\documentclass{article}

\usepackage[utf8]{inputenc}
\usepackage[T1]{fontenc}
\usepackage[english]{babel}
\usepackage{csquotes}
\usepackage[margin=2.5cm]{geometry}

\usepackage{amsmath,amsfonts,amsthm,dsfont,bm,stmaryrd,amssymb}
\usepackage[colorlinks=true,linkcolor=blue,citecolor=red]{hyperref}
\usepackage{cleveref}

\usepackage{comment}
\usepackage{enumerate}
\usepackage{hhline}
\usepackage{multirow}

\theoremstyle{plain}
\newtheorem{thm}{Theorem}[section]
\newtheorem*{thmstar}{Theorem}
\newtheorem{cor}[thm]{Corollary}
\newtheorem{prop}[thm]{Proposition}
\newtheorem{lm}[thm]{Lemma}

\theoremstyle{definition}
\newtheorem{df}[thm]{Definition}
\newtheorem{algo}[thm]{Algorithm}
\newtheorem{Rk}[thm]{Remark}

\makeatletter
\newenvironment{pf}[1][\proofname]{\par
    \pushQED{\qed}%
    \normalfont \topsep6\p@\@plus6\p@\relax
    \trivlist
    \item\relax
        {\bfseries
    #1\@addpunct{:}}\hspace\labelsep\ignorespaces
}{%
    \popQED\endtrivlist\@endpefalse
}
\makeatother

\numberwithin{equation}{section}

\newcommand{\ent}[2]{\llbracket #1,#2 \rrbracket}
\renewcommand{\geq}{\geqslant}
\renewcommand{\leq}{\leqslant}

\renewcommand{\bar}[1]{\overline{#1}}
\renewcommand{\tilde}[1]{\widetilde{#1}}

\newcommand{\Z}{\mathbb{Z}}
\newcommand{\Q}{\mathbb{Q}}

\newcommand{\C}{\mathbb{C}}
\newcommand{\F}{\mathbb{F}}	
\renewcommand{\S}{\mathfrak{S}} 
\newcommand{\A}{\mathfrak{A}} 
\renewcommand{\P}{\mathbb{P}}

\DeclareMathOperator{\GL}{\mathrm{GL}}
\DeclareMathOperator{\PGL}{\mathrm{PGL}}
\DeclareMathOperator{\SL}{\mathrm{SL}}
\DeclareMathOperator{\PSL}{\mathrm{PSL}}
\DeclareMathOperator{\Tr}{\mathrm{Tr}}

\DeclareMathOperator{\lcm}{\mathrm{lcm}}

\DeclareMathOperator{\Mod}{M}
\DeclareMathOperator{\Smod}{S}
\DeclareMathOperator{\Gal}{\mathrm{Gal}}
\newcommand{\Id}{\mathrm{Id}}

\renewcommand{\Im}{\mathrm{Im}}

\usepackage[backend=biber]{biblatex}
\title{Explicit Small Image Theorems for Residual Modular Representations}
\author{Baptiste Peaucelle}
\date{}

\AtEndDocument{\bigskip{\footnotesize%
  Université Clermont Auvergne, CNRS, LMBP, F-63000 CLERMONT-FERRAND, FRANCE \par 
  E-mail address: \texttt{baptiste.peaucelle@uca.fr}
}}

\begin{document}
\maketitle
\begin{abstract}
Let $\bar{\rho}_{f,\lambda}$ be the residual Galois representation attached to a newform $f$ and a prime ideal $\lambda$ in the integer ring of its coefficient field. In this paper, we prove explicit bounds for the residue characteristic of the prime ideals $\lambda$ such that $\bar{\rho}_{f,\lambda}$ is exceptional, that is reducible, of projective dihedral image, or of projective image isomorphic to $\A_{4}$, $\S_{4}$ or $\A_{5}$. We also develop explicit criteria to check the reducibility of $\bar{\rho}_{f,\lambda}$, leading to an algorithm that compute the exact set of such $\lambda$. We have implemented this algorithm in PARI/GP. Along the way, we construct lifts of Katz' $\theta$ operator in character zero, and we prove a new Sturm bound theorem.
\end{abstract}
\section*{Introduction}
\subsection*{The setup}

Let $f$ be a newform of weight $k \geq 2$, level $N \geq 1$, and character $\varepsilon : \left(\Z/N\Z\right)^{\times} \to \C^{\times}$. We denote by $K_{f}$ the number field generated by the Fourier coefficients of $f$, and by $\mathcal{O}_{f}$ its ring of integers. Given a prime number $\ell$, let us write $\rho_{f,\ell}$ for the $\ell$-adic representation attached to $f$ by Deligne:
\[\rho_{f,\ell} : G_{\Q} \to \GL_{2}(K_{f} \otimes \Q_{\ell}),\]
where $G_{\Q}$ denotes the Galois group of an algebraic closure $\bar{\Q}$ of $\Q$. Choosing a stable lattice for $\rho_{f,\ell}$, we can assume it takes values in $\GL_{2}(\mathcal{O}_{f} \otimes \Z_{\ell})$. For a prime ideal $\lambda$ of $\mathcal{O}_{f}$, let us denote by $\mathcal{O}_{f,\lambda}$ the completion of $\mathcal{O}_{f}$ at $\lambda$. The decomposition $\mathcal{O}_{f} \otimes \Z_{\ell} = \prod\limits_{\lambda \mid \ell} \mathcal{O}_{f,\lambda}$, where the product runs over all the prime ideals in $\mathcal{O}_{f}$ above $\ell$, gives rise to a Galois representation $\rho_{f,\lambda}$ with values in $\GL_{2}(\mathcal{O}_{f,\lambda})$. Reducing $\rho_{f,\lambda}$ modulo the maximal ideal of $\mathcal{O}_{f,\lambda}$ and taking the semi-simplification, we finally get a residual representation $\bar{\rho}_{f,\lambda} : G_{\Q} \to \GL_{2}(\F_{\lambda})$ where $\F_{\lambda}$ denotes the residue field of $\mathcal{O}_{f,\lambda}$. Notice that $\bar{\rho}_{f,\lambda}$ does not depend on the choice of the stable lattice. In this article, we will be interested in the following result which has been proved by Ribet in 1985 \cite{Ribet85}:

\begin{thmstar}[Ribet, 1985]\label{Ribet}
For all but finitely many $\lambda$ the representation $\bar{\rho}_{f,\lambda}$ is irreducible. Furthermore, if $f$ is not a form with complex multiplication (see \cref{CM}), then for all but finitely many $\lambda$, the order of the image of $\bar{\rho}_{f,\lambda}$ is divisible by the residue characteristic of $\lambda$.
\end{thmstar}
For simplicity, we shall say that $\lambda$ is exceptional if one of the conclusions of Ribet's theorem does not hold. Let us denote by $\P\bar{\rho}_{f,\lambda}$ the projectivisation of $\bar{\rho}_{f,\lambda}$, that is the composition of $\bar{\rho}_{f,\lambda}$ with the canonical projection map from $\GL_{2}(\F_{\lambda})$ to $\PGL_{2}(\F_{\lambda})$. According to the classification of the subgroups of $\PGL_{2}(\F_{\lambda})$ and $\PSL_{2}(\F_{\lambda})$ \cite[\nopp II.8.27]{Huppert67}, if $\lambda$ is exceptional, then one of the following is true:
\begin{enumerate}
\item the representation $\bar{\rho}_{f,\lambda}$ is reducible;\label{red}
\item the image of $\P\bar{\rho}_{f,\lambda}$ is a dihedral group $D_{2n}$ with $\ell \nmid n$;\label{dih}
\item the image of $\P\bar{\rho}_{f,\lambda}$ is isomorphic to $\A_{4}$, $\S_{4}$ or $\A_{5}$.\label{exp}
\end{enumerate}
Therefore, Ribet's theorem is equivalent to showing that the above three cases occur finitely many times (under the assumption that $f$ is non-CM in the second case).

The theorem generalised results of Ribet from 1975 for $N = 1$ \cite{Ribet75}, which already extended those of Serre and Swinnerton-Dyer from 1973 for $N = 1$ and $K_{f} = \Q$ \cite{Serre73,Swinnerton-Dyer73}. Although the first result of Ribet from 1975 gave an explicit description of the prime ideals for which the associated representation is reducible, it was no more the case in 1985. The second part of the theorem was already ineffective in 1975. The first step in making Ribet's result effective has been accomplished by Billerey and Dieulefait in 2014 \cite{Billerey14}. Assuming the character of $f$ is trivial, they gave explicit criteria for the residue characteristics $\ell$ of $\lambda$ in terms of $k$ and $N$, for $\bar{\rho}_{f,\lambda}$ to be reducible. In the two other cases, they gave explicit bounds for $\ell$, in terms of $k$ and $N$.

\subsection*{Statement of results}

In this paper, we extend the results of \cite{Billerey14} to all newforms of arbitrary weight, level and character. The argument used to deal with the third case given in \cite{Billerey14} can be applied without modifications to a form with non-trivial character. We recall their result for the sake of completeness.
\begin{thmstar}
If the image of $\P\bar{\rho}_{f,\lambda}$ is isomorphic to $\A_{4}$, $\S_{4}$ or $\A_{5}$, then either $\ell \mid N$ or $\ell \leq 4k-3$.
\end{thmstar}
We now focus on the reducible and dihedral cases. The main idea to deal with these cases is to translate them into congruences between modular forms. In the dihedral case, we get a congruence between twists of $f$, while in the reducible case we deal with congruences involving Eisenstein series. To deduce a bound from a congruence between two modular forms, the so-called Sturm bounds are usually used. We generalise them to deal with modular forms of arbitrary weight, level and character (see \cref{Sturm_theta}). In the dihedral projective image case, the bound we get for the residue characteristic is the following.
\begin{thmstar}[\cref{Dihthm}]\label{Dihthm_Intro}
Assume $\bar{\rho}_{f,\lambda}$ has dihedral projective image. If $N = 1$ then we have $\ell \leq k$ or $\ell \in \{2k - 1,2k-3\}$. Else, if $N \geq 2$ and $f$ does not have complex multiplication, then we have
\[\ell \leq \left(2N^{\frac{k-1}{2}}\right)^{[K_{f}:\Q]} \times \max\left(
    \left(\frac{k}{3}\left(2\log\log(N)+2.4\right)\right)^{\frac{k-1}{2}},
    \left(\frac{5}{2}N^{\frac{k-1}{2}}\right)
\right)^{[K_{f}:\Q]}.\]
\end{thmstar}
This result gives us indeed an upper bound for $\ell$ in terms of $N$ and $k$ because $[K_{f}:\Q]$ can be bounded by the dimension of the $\C$-vector space generated by the newforms of weight $k$, level $N$ and character $\varepsilon$ (see for example \cite{Martin05}).

The most challenging case is when $\bar{\rho}_{f,\lambda}$ is reducible. In this situation we generalise the results of \cite[Section 2.]{Billerey14} to newforms of arbitrary character. The restriction on the character in \cite{Billerey14} was mainly due to a partial knowledge of the constant term of Eisenstein series at arbitrary cusps. This computation has been done in full generality in \cite{Billerey18}, allowing us to generalise their result. The following theorem then follows from combining this technical result with a detail study of modular reducible representations, hence extending the strategy used for the proof of \cite[Theorem 2.7]{Billerey14}.
\begin{thmstar}[\cref{red_bound}]\label{red_boundIntro}
Assume $\bar{\rho}_{f,\lambda}$ to be reducible. Then one of following holds:
\begin{enumerate}
\item $\ell \leq k+1$;
\item $\ell \mid N\varphi(N)$;
\item there exists a prime-to-$\ell$ primitive Dirichlet character $\eta$ of conductor $\mathfrak{c}_{0} \mid N$ and such that $\eta(-1) = (-1)^{k}$ and $\ell$ divides the algebraic norm of one of the non-zero following quantities:
\begin{enumerate}
\item $p^{k} - \eta(p)$ for a prime number $p \mid N$;
\item the $k$-th Bernoulli number $B_{k,\eta}$ attached to $\eta$ (see \cref{Ber_df}).
\end{enumerate}
\end{enumerate}
\end{thmstar}
The precise study of reducible modular representations used in the proof of the previous theorem is the main novelty of the present paper. The basic question we consider is as follows: How to characterise the reducibility of $\bar{\rho}_{f,\lambda}$ by a finite number of explicit congruences? We give two answers to this question. A general one that applies without any restriction on $\ell$ or $f$, and, under some assumptions on $\ell$, a second one for which the number of congruences to check is independent from $\ell$. A weaker form of our first main result in this direction is the following.

\begin{thmstar}[\cref{red_thm}]
The following are equivalent:
\begin{enumerate}
\item $\bar{\rho}_{f,\lambda}$ is reducible;
\item Let $\mathfrak{L}$ be a place of $\bar{\Q}$ above $\lambda$. There exist two primitive Dirichlet characters $\varepsilon_{1}$, $\varepsilon_{2}$ of conductor $\mathfrak{c}_{1}$, $\mathfrak{c}_{2}$ respectively such that $\mathfrak{c}_{1}\mathfrak{c}_{2} \mid N$, and two integers $m_{1}$, $m_{2}$ such that $0 \leq m_{1} \leq m_{2} \leq \ell -2$ and $\bar{\chi}_{\ell}^{m_{1}+m_{2}}\varepsilon_{1}\varepsilon_{2} \equiv \bar{\chi}_{\ell}^{k-1}\varepsilon \pmod{\mathfrak{L}}$. Define 
\[\tilde{k} = \left\{\begin{array}{ll}
3 + \max(k,m_{2}+2m_{1}+1) & \text{if } \ell \mid N\\
\ell + 5 + \max(k,m_{2}+\ell m_{1} + 1) & \text{if } \ell \nmid N
\end{array}\right..\]
For every prime number $p \leq \frac{N \tilde{k}}{3} \prod\limits_{\underset{q \text{ prime}}{q \mid 2N}} \left(1+\frac{1}{q}\right)$ and not dividing $2\ell$, we have
\begin{itemize}
\item $p \nmid N$ and $a_{p}(f) \equiv p^{m_{1}}\varepsilon_{1}(p) + p^{m_{2}}\varepsilon_{2}(p) \pmod{\mathfrak{L}}$;
\item $p \mid N$ and $a_{p}(f) \equiv p^{m_{1}}b_{p} \pmod{\mathfrak{L}}$ for some $b_{p} \in \{0,\varepsilon_{1}(p),p^{m_{2}-m_{1}}\varepsilon_{2}(p)\}$.
\end{itemize}
\end{enumerate}
When this holds, we moreover have $\bar{\rho}_{f,\lambda} \cong \bar{\chi}_{\ell}^{m_{1}}\bar{\varepsilon_{1}} \oplus \bar{\chi}_{\ell}^{m_{2}}\bar{\varepsilon_{2}}$.
\end{thmstar}
Notice that this result applies with no assumption on $f$ and $\ell$. In particular, it can be used to check the reducibility of $\bar{\rho}_{f,\lambda}$ for any given $\lambda$, including the ones which residue characteristic is small compared to the weight, or divides the level. Such restrictions appear for instance in the work of Anni \cite[Algorithm 7.2.4]{Anni13}, where the author develops a different, “bottom-up” approach, towards these questions in the context of modular forms “à la Katz”.

In the previous statement, the number of congruences to be satisfied in order to prove the reducibility of $\bar{\rho}_{f,\lambda}$ depends not only on $N$, $k$ and $\varepsilon$, but also on $\ell$. Under some assumptions on $\ell$, we have been able to remove this dependency in the bound. A weaker form of our second main result can be stated as follows.

\begin{thmstar}[\Cref{bigred_thm}]\label{bigred_thmIntro}
Assume $\ell > k+1$ and $\ell \nmid N\varphi(N)$, where $\varphi$ denotes the Euler totient function. The following are equivalent:
\begin{enumerate}
\item $\bar{\rho}_{f,\lambda}$ is reducible;
\item Let $\mathfrak{L}$ be a place of $\bar{\Q}$ above $\lambda$. There exist two Dirichlet characters $\varepsilon_{1}$, $\varepsilon_{2}$ of conductor $\mathfrak{c}_{1}$, $\mathfrak{c}_{2}$ respectively such that $\mathfrak{c}_{1}\mathfrak{c}_{2} \mid N$, and $\varepsilon_{1}\varepsilon_{2} = \varepsilon$. For all odd primes $p \leq \frac{Nk}{3} \prod\limits_{\underset{q \text{ prime}}{q \mid 2N}} \left(1 + \frac{1}{q}\right)$, we have
\begin{itemize}
\item $p \nmid N$ and $a_{p}(f) \equiv \varepsilon_{1}(p) + p^{k-1}\varepsilon_{2}(p) \pmod{\mathfrak{L}}$;
\item $p \mid N$ and $a_{p}(f) \equiv b_{p} \pmod{\mathfrak{L}}$ for some $b_{p} \in \{0,\varepsilon_{1}(p),p^{k-1}\varepsilon_{2}(p)\}$.
\end{itemize}
\end{enumerate}
When this holds, we moreover have $\bar{\rho}_{f,\lambda} \cong \bar{\varepsilon_{1}} \oplus \bar{\chi}_{\ell}^{k-1}\bar{\varepsilon_{2}}$.
\end{thmstar}
We stress the fact that according to these latter two results, proving the reducibility of $\bar{\rho}_{f,\lambda}$ requires checking almost $Nk^{2}\log\log(N)$ congruences (and even $Nk\log\log(N)$ for the primes $\ell$ satisfying the assumptions $\ell > k+1$ and $\ell \nmid N\varphi(N)$). Notice that the $\log\log(N)$ part comes from the upper-bound we prove in \cref{Sturm_upperbound}.

To achieve such bound, we extensively use the local description of $\bar{\rho}_{f,\lambda}$ at the bad prime numbers (\textit{i.e.} the prime numbers dividing $N$), together with generalised Sturm bounds theorems and an appropriate use of degeneracy maps between modular forms spaces of various levels. Having a sharp bound is especially important from a computational point of view. Indeed, those results also provide us with a algorithm that explicitly computes the exact set of $\lambda$ such that $\bar{\rho}_{f,\lambda}$ is reducible. We have implemented such an algorithm in PARI/GP \cite{PARI}.

\subsection*{Organisation of the paper}

In the first section of the paper we collect various results on the local description of residual modular representations. The second section is devoted to recall basic facts about Dirichlet characters and Eisenstein series. In \cref{sec_prelim}, we develop various tools that will play a central role in the proof of the main results. In paragraph \ref{sec_theta}, we construct the so-called theta operators in characteristic zero, that we then use to generalise Sturm bounds results in paragraph \ref{sec_Sturm}. In paragraph \ref{sec_raise}, we introduce various other operators acting on modular forms spaces that allow us to modify some Fourier coefficients of a given form. The two next sections are devoted to prove the results involving reducible residual modular representations and dihedral representations respectively. Finally, we explain in the last section how to translate our theoretic results on reducibility into an algorithm that takes a modular form as input, and outputs the exact list of prime ideals $\lambda$ such that $\bar{\rho}_{f,\lambda}$ is reducible. We then illustrate this algorithm on numerical examples.

\subsection*{Acknowledgements}

The author really wants to thank Nicolas Billerey for his great advice during the redaction of this paper. The author also wishes to thank Samuele Anni for stimulating discussions at Grenoble, and Bill Allombert and Karim Belabas for the fruitful exchanges about PARI/GP and the numerical implementation of this work.

\subsection*{General notations}

In the whole paper, for two positive integers $k$ and $N$, and a Dirichlet character $\varepsilon$ modulo $N$, the notations $\Mod_{k}(N,\varepsilon)$, $\Smod_{k}(N,\varepsilon)$ and $\Smod_{k}^{\mathrm{new}}(N,\varepsilon)$ stand for the complex vector spaces of modular forms, cuspidal modular forms, and new modular forms of weight $k$, level $N$, and character $\varepsilon$ respectively. We also denote by $\Mod(N)$ the graded algebra of modular forms of level $N$. For any modular form $g$ in $\Mod(N)$, we always write $\sum\limits_{n = 0}^{\infty} a_{n}(g) q^{n}$ for its $q$-expansion at infinity, with $q = e^{2i\pi z}$ and $z$ in the complex upper-half plane.

\section{Background on residual modular Galois representations}\label{sec_BackGal}
We fix for this section a newform $f$ of weight $k \geq 2$, level $N \geq 1$, and character $\varepsilon$. Let $\ell$ be a prime number and let $\lambda$ be a prime ideal of $\mathcal{O}_{f}$ above $\ell$. We recall some results about the shape of $\bar{\rho}_{f,\lambda}$ restricted to a decomposition subgroup. In the following, given a prime number $p$ we denote by $G_{p}$ a decomposition subgroup of $G_{\Q}$ at $p$, and by $I_{p}$ its inertia subgroup. For some $\lambda$-integer $x$, $\mu_{p}(x)$ stands for the unique unramified character of $G_{p}$ sending a Frobenius element at $p$ to the reduction of $x$ modulo $\lambda$. We will simply write $\mu(x)$ when no confusion on $p$ can arise. We finally write $\bar{\chi}_{\ell}$ for the cyclotomic character modulo $\ell$.

We begin by the local description of $\bar{\rho}_{f,\lambda}$ at $\ell$ that has been established by Deligne and Fontaine. Their result is the following.

\begin{prop}[Deligne-Fontaine, {\cite[\nopp 2.4]{Edixhoven92}}]\label{DF}
Assume $2 \leq k \leq \ell + 1$ and $\ell \nmid N$.
\begin{itemize}
\item If $f$ is ordinary at $\lambda$ (that is if $a_{\ell}(f) \not\equiv 0 \pmod{\lambda}$), then we have
\[\bar{\rho}_{f,\lambda|G_{\ell}} \cong \begin{pmatrix}
\bar{\chi}_{\ell}^{k-1}\mu\left(\frac{\varepsilon(\ell)}{a_{\ell}(f)}\right) & \star\\
0 & \mu\left(a_{\ell}(f)\right)
\end{pmatrix},\]
\item If not, then $\bar{\rho}_{f,\lambda|G_{\ell}}$ is irreducible, and we have
\[\bar{\rho}_{f,\lambda|I_{\ell}} \cong \begin{pmatrix}
\psi^{k-1} & 0\\
0 & \psi'^{k-1}
\end{pmatrix}.\]
\end{itemize}
Where $\{\psi,\psi'\} = \{\psi,\psi^{\ell}\}$ is the set of fundamental characters of level $2$ (see \cite[\nopp 2.4]{Edixhoven92}).
\end{prop}

For the primes $p$ different from $\ell$ and dividing $N$, the shape of $\bar{\rho}_{f,\lambda|G_{p}}$ has been computed by Langlands and compiled in \cite[Proposition 2.8]{Loeffler12}.
\begin{prop}\label{LW}
Let $p \neq \ell$ be a prime dividing $N$ and let $\mathfrak{c}$ be the conductor of $\varepsilon$. We denote by $v_{p}$ the $p$-adic valuation.
\begin{itemize}
\item If $v_{p}(N) = 1$ and $v_{p}(\mathfrak{c}) = 0$, then we have
\[\bar{\rho}_{f,\lambda|G_{p}} \cong \begin{pmatrix}
\mu(a_{p}(f)) \bar{\chi}_{\ell} & \star\\
0 & \mu(a_{p}(f)) 
\end{pmatrix};\]
\item If $v_{p}(N) = v_{p}(\mathfrak{c})$, then $a_{p}(f)$ is a unit in $\mathcal{O}_{f,\lambda}$ and we have
\[\bar{\rho}_{f,\lambda|G_{p}} \cong \mu\left(a_{p}(f)\right) \oplus \mu\left(a_{p}(f)^{-1}\right) \bar{\chi}_{\ell}^{k-1} \bar{\varepsilon_{|G_{p}}},\]
where $\bar{\varepsilon_{|G_{p}}}$ stands for the reduction modulo $\lambda$ of the restriction of $\varepsilon$ to $G_{p}$.
\end{itemize}
\end{prop}
\begin{pf}
From \cite[Theorem 4.6.17]{Miyake06}, we have in the second case $|a_{p}(f)|^{2} = p^{k-1}$. Therefore, it is indeed invertible in $\mathcal{O}_{f,\lambda}$ because $p \neq \ell$.

The only thing to prove is that the hypothesis of \cite[Proposition 2.8]{Loeffler12} holds in our cases, namely that $f$ is $p$-primitive in the terminology of \cite[Definition 2.7]{Loeffler12}. To do so, we use \cite[Theorem]{Loeffler15}. We recall a direct consequence of this result: Define $u = \min\left(\left\lfloor \frac{v_{p}(N)}{2}\right\rfloor,v_{p}(N)-v_{p}(\mathfrak{c})\right)$. If $u = 0$, then $f$ is $p$-primitive. We easily check that in our two cases we have $u = 0$.
\end{pf}

Since the work of Carayol \cite[Théorème (A)]{Carayol86}, it is well-known that the value of the prime-to-$\ell$ part of the Artin conductor of $\rho_{f,\lambda}$ is equal to the prime-to-$\ell$ part of $N$. The study of the behaviour of the Artin conductor after reduction (\textit{i.e.} the one of $\bar{\rho}_{f,\lambda}$) has then been established independently by Carayol \cite{Carayol89} and Livné \cite{Livne89}. Here is their result.

\begin{prop}[Carayol-Livné]\label{CL}
Let $N(\bar{\rho}_{f,\lambda})$ be the prime-to-$\ell$ part of the Artin conductor of $\bar{\rho}_{f,\lambda}$. We then have
\[N(\bar{\rho}_{f,\lambda}) \mid N.\]
Moreover, write $e_{p} := v_{p}(N) - v_{p}(N(\bar{\rho}_{f,\lambda}))$ for a prime $p \neq \ell$. If $e_{p} > 0$, then we have $v_{p}(N)-v_{N(\bar{\rho}_{f,\lambda})} \in \{1,2\}$.
\end{prop}
\section{Background on Dirichlet characters and Eisenstein series}\label{sec_BackMF}

\subsection{Generalised Bernoulli numbers and Gauß sums}\label{sec_char}
Let $\varepsilon$ be a primitive Dirichlet character of conductor $\mathfrak{c}$. We recall the definition and properties of the Gauß sums and generalised Bernoulli numbers attached to $\varepsilon$.

\begin{df}\label{Ber_df}
The Bernoulli numbers $(B_{m,\varepsilon})_{m \geq 0}$ attached to $\varepsilon$ are defined by the following generating series:
\[\sum_{n=1}^{\mathfrak{c}} \varepsilon(n) \frac{te^{nt}}{e^{\mathfrak{c}t} - 1} = \sum_{m = 0}^{\infty} B_{m,\varepsilon} \frac{t^{m}}{m!}.\]
In particular, when $\varepsilon$ is odd, we have $B_{1,\varepsilon} = \frac{1}{\mathfrak{c}} \sum\limits_{n = 1}^{\mathfrak{c}-1} n\varepsilon(n)$, and when $\varepsilon$ is both even and non-trivial, we have $B_{2,\varepsilon} = \frac{1}{\mathfrak{c}} \sum\limits_{n = 1}^{\mathfrak{c}-1} n^{2}\varepsilon(n)$.
\end{df}

\begin{Rk}
If $\varepsilon = \mathds{1}$ is the trivial character modulo $1$, we get the classical Bernoulli numbers except when $m = 1$, in which case we have $B_{1,\mathds{1}} = \frac{1}{2} = -B_{1}$.
\end{Rk}
We state below the main properties of the Bernoulli numbers. First, we exactly know when the Bernoulli numbers vanish (see \cite[Theorem 3.3.4]{Miyake06} for a proof).

\begin{prop}\label{Ber_vanishing}
We have $B_{m,\varepsilon} = 0$ if and only if $\varepsilon(-1) \neq (-1)^{m}$.
\end{prop}
Secondly, the behaviour of the Bernoulli numbers after reduction modulo a prime ideal has been studied by Van-Staudt \cite{Staudt40} in the case $\varepsilon = \mathds{1}$, and by Carlitz \cite[Theorem 1]{Carlitz59} in the case $\varepsilon \neq \mathds{1}$. We summarise their results in the following proposition.

\begin{prop}\label{Ber_divisors}
Let $m$ be a positive integer.
\begin{enumerate}
\item Let $\ell$ be a prime number. If $\ell-1$ divides $m$, then we have the congruence $\ell B_{m,\mathds{1}} \equiv -1 \pmod{\ell}$. Otherwise, $\frac{B_{m,\mathds{1}}}{m}$ is $\ell$-integral and its reduction modulo $\ell$ only depends on $m$ modulo $\ell-1$. In particular, the denominator of $B_{m,\mathds{1}}$ is equal to $\prod\limits_{\underset{\ell-1 \mid m}{\ell \text{ prime}}} \ell$.
\item For $\varepsilon \neq \mathds{1}$, write $\frac{B_{m,\varepsilon}}{m} = \mathfrak{N}\mathfrak{D}^{-1}$, with $\mathfrak{N}$ and $\mathfrak{D}$ two coprime ideals of $\Z[\varepsilon]$, the ring spanned by the image of $\varepsilon$. If the conductor of $\varepsilon$ admits at least two distinct prime factors, then $\mathfrak{D} = 1$. Otherwise, if the conductor of $\varepsilon$ is a power of a prime number $\ell$, then $\mathfrak{D}$ contains only prime ideals above $\ell$.
\end{enumerate}
\end{prop}

Another classical quantity attached to Dirichlet characters is the Gauß sum. We recall its definition and properties below.
\begin{df}
The Gauß sum attached to $\varepsilon$ is
\[W(\varepsilon) = \sum\limits_{n=1}^{\mathfrak{c}} \varepsilon(n) e^{\frac{2i\pi n}{\mathfrak{c}}}.\]
\end{df}
One can find the following result in \cite[Lemma 2.1.]{Billerey14}.
\begin{prop}\label{W_divisors}
The prime divisors of the algebraic norm of $W(\varepsilon)$ are those of $\mathfrak{c}$.
\end{prop}

\subsection{Teichmüller lifts}\label{sec_Tlifts}
We present here the behaviour of roots of unity after reduction modulo $\mathfrak{L}$.
\begin{lm}\label{Rootmod}
Let $n$ be a positive integer and let $\zeta$ be a primitive $n$-th root of unity in $\bar{\Q}$. Let $\ell$ a prime number and let $\mathfrak{L}$ be a place of $\bar{\Q}$ above $\ell$. We have $\zeta \equiv 1 \pmod{\mathfrak{L}}$ if and only if $n$ is a power of $\ell$. In particular, a Dirichlet character is trivial modulo $\mathfrak{L}$ if and only if it has order a power of $\ell$.
\end{lm}
\begin{pf}
According to \cite[Proposition 3.5.4.]{Cohen07}, the algebraic norm of $1-\zeta$ over $\Q(\zeta)$ is equal to:
\[\left\{\begin{array}{ll}
0 & \text{if } n = 1;\\
q & \text{if } n = q^{r} \text{ with } q \text{ prime and } r \geq 1;\\
1 & \text{otherwise}.
\end{array}\right.\]
Thus, if $n$ is not an $\ell$-power, then $\ell$ does not divide the norm of $\zeta-1$ and we have $\zeta \not\equiv 1 \pmod{\mathfrak{L}}$. Assume $n = \ell^{r}$, $r \geq 1$. We then have
\[\ell\Z[\zeta] = (1-\zeta)^{\ell^{r-1}(\ell-1)}\Z[\zeta].\]
Thus, the only prime ideal above $\ell$ in $\Z[\zeta]$ is $(1-\zeta)\Z[\zeta]$ and we therefore have $\zeta \equiv 1 \pmod{\mathfrak{L}}$.

For the second part of the lemma, let $\varepsilon$ be a Dirichlet character modulo $N$. From above, $\varepsilon$ is trivial modulo $\mathfrak{L}$ if and only if for all $x \in \left(\Z/N\Z\right)^{\times}$, $\varepsilon(x)$ is a root of unity of order a power of $\ell$. This can happen if and only if $\varepsilon$ has order a power of $\ell$.
\end{pf}
\Cref{Rootmod} implies that the kernel of the reduction modulo $\mathfrak{L}$ from the group of all roots of unity to $\bar{\F}_{\ell}^{\times}$, is the subgroup of primitive roots of unity of order a power of $\ell$. In particular, the restriction of this map to the subgroup of roots of unity of order prime to $\ell$ is injective. Moreover, because the subgroup of roots of unity of order $\ell^{n}-1$ maps to $\F_{\ell^{n}}^{\times}$, it is onto and therefore an isomorphism. The inverse map
\[T_{\mathfrak{L}}: \bar{\F}_{\ell}^{\times} \to \{\zeta \in \C^{\times}, \gcd(\ell,\mathrm{ord}(\zeta)) = 1\},\]
is the so-called Teichmüller lift with respect to the place $\mathfrak{L}$. It will allow us to lift multiplicatively characters modulo $\mathfrak{L}$ to Dirichlet characters of prime-to-$\ell$ order. Moreover, the Dirichlet characters that arise this way are exactly the ones that have the same conductor as a Dirichlet character and as a character modulo $\mathfrak{L}$.

\subsection{Eisenstein series}\label{sec_Eisenstein}
Let $k$ be a positive integer and let $\varepsilon_{1}$, $\varepsilon_{2}$ be two Dirichlet characters modulo $\mathfrak{c}_{1}$ and $\mathfrak{c}_{2}$ respectively, such that $\varepsilon_{1}\varepsilon_{2}(-1) = (-1)^{k}$. Moreover, if $k = 2$ and $\varepsilon_{1}$, $\varepsilon_{2}$ are both trivial, then assume $\mathfrak{c}_{1} = 1$ and $\mathfrak{c}_{2}$ is a prime number, otherwise, assume $\varepsilon_{1}$ and $\varepsilon_{2}$ are primitive. For a complex number $z$ in the complex upper-half plane $\mathcal{H}$, consider the following $q$-expansion:
\begin{equation}\label{qexp}
E_{k}^{\varepsilon_{1},\varepsilon_{2}}(z) := C + \sum_{n=1}^{\infty} \sigma_{k-1}^{\varepsilon_{1},\varepsilon_{2}}(n) q^{n},
\end{equation}
with $\sigma_{r}^{\varepsilon_{1},\varepsilon_{2}}(n) = \sum\limits_{0 < d \mid n} d^{r}\varepsilon_{1}\left(\frac{n}{d}\right)\varepsilon_{2}(d)$ for any $r \geq 0$ and
\begin{equation}\label{cst}
C = \left\{\begin{array}{cl}
0 & \left|\begin{array}{ll}
&\text{if } k \geq 2 \text{ and } \varepsilon_{1} \neq \mathds{1},\\
\text{or} &\text{if } k = 1 \text{ and } \varepsilon_{1}, \varepsilon_{2} \text{ are both non-trivial};
\end{array}\right.\\
\dfrac{1}{24}(\mathfrak{c}_{2}-1) & \text{if } k = 2 \text{ and } \varepsilon_{1}, \varepsilon_{2} \text{ both trivial};\\
-\dfrac{B_{k,\varepsilon_{1}\varepsilon_{2}}}{2k} & \text{otherwise}.
\end{array}\right.
\end{equation}
The following result is proved in \cite[Theorem 4.7.1]{Miyake06} and \cite[(4.7.16)]{Miyake06}.

\begin{prop}\label{Eisenstein}
The $q$-series $E_{k}^{\varepsilon_{1},\varepsilon_{2}}$ defines a modular form of weight $k$, level $\mathfrak{c}_{1}\mathfrak{c}_{2}$ and character $\varepsilon_{1}\varepsilon_{2}$. It is a normalised eigenform for all the Hecke operators at level $\mathfrak{c}_{1}\mathfrak{c}_{2}$.
\end{prop}
For $\varepsilon_{1}=\varepsilon_{2}=\mathds{1}$, the definition of the series $E_{k}^{\mathds{1},\mathds{1}}$ agrees with the definition of the classical Eisenstein series of weight $k$. We simply write it $E_{k}$ in this case. For $k = 2$, we denote by $E_{2}$ the $q$-series
\[E_{2}(z) = -\frac{1}{24} + \sum\limits_{n = 1}^{\infty} \Bigg(\sum\limits_{0 < d \mid n} d\Bigg) q^{n}.\]
Note that this formula defines a holomorphic function on $\mathcal{H}$, but $E_{2}$ is not modular.

In the case $k \geq 2$ and $\varepsilon_{1}$, $\varepsilon_{2}$ primitive, the behaviour of the constant coefficient of $E_{k}^{\varepsilon_{1},\varepsilon_{2}}$ at a cusp of $\Gamma_{1}(N)$ has been computed in \cite[Proposition 4]{Billerey18}. It states the following:

\begin{prop}\label{Eisenstein_cst}
Assume $k \geq 2$ and $\varepsilon_{1}$, $\varepsilon_{2}$ are primitive. Let $M$ be a positive integer and let $\gamma = \begin{pmatrix} u & \beta\\v & \delta\end{pmatrix} \in \SL_{2}(\Z)$. Put $\tilde{v} := \frac{v}{\gcd(v,M)}$ and $\tilde{M} := \frac{M}{\gcd(v,M)}$. We define
\[\Upsilon_{k}^{\varepsilon_{1},\varepsilon_{2}}(\gamma,M) := \lim\limits_{\Im(z) \to +\infty} (E_{k}^{\varepsilon_{1},\varepsilon_{2}}(M\cdot)|_{k}\gamma)(z),\]
where we denote by $|_{k}$ the classical slash action of weight $k$.
 
If $\mathfrak{c}_{2} \nmid \tilde{v}$ then $\Upsilon_{k}^{\varepsilon_{1},\varepsilon_{2}}(\gamma,M) = 0$. Otherwise, if $\mathfrak{c}_{2} \mid \tilde{v}$ then $\Upsilon_{k}^{\varepsilon_{1},\varepsilon_{2}}(\gamma,M) \neq 0 \Leftrightarrow \gcd\left(\mathfrak{c}_{1},\frac{\tilde{v}}{\mathfrak{c}_{2}}\right) = 1$. In this case, we moreover have
\[\Upsilon_{k}^{\varepsilon_{1},\varepsilon_{2}}(\gamma,M) = -\varepsilon_{2}^{-1}\left(\tilde{M}u\right)\varepsilon_{1}\left(-\frac{\tilde{v}}{\mathfrak{c}_{2}}\right) \frac{W((\varepsilon_{1}\varepsilon_{2}^{-1})_{0})}{W(\varepsilon_{2}^{-1})}
\frac{B_{k,(\varepsilon_{1}^{-1}\varepsilon_{2})_{0}}}{2k}
\left(\frac{\mathfrak{c}_{2}}{\tilde{M}\mathfrak{c}_{0}}\right)^{k} \prod_{p \mid \mathfrak{c}_{1}\mathfrak{c}_{2}} \left(1 - \frac{(\varepsilon_{1}\varepsilon_{2}^{-1})_{0}(p)}{p^{k}}\right),\]
where $\chi_{0}$ denotes the primitive character associated to a Dirichlet character $\chi$, and $\mathfrak{c}_{0}$ the conductor of $\varepsilon_{1}^{-1}\varepsilon_{2}$.
\end{prop}
The proof of \cite{Billerey18} is only given in the cases $k \geq 3$ and $k = 2$, $\varepsilon_{1}$, $\varepsilon_{2}$ non-trivial, but we easily see that the proof is still valid in the case $k = 2$, $\varepsilon_{1} = \varepsilon_{2} = \mathds{1}$.

\section{Preliminary results on modular forms}\label{sec_prelim}
\subsection{Theta operators}\label{sec_theta}
We fix for this paragraph a place $\mathfrak{L}$ of $\bar{\Q}$. Consider the operator $\theta$ acting on the space of holomorphic functions on $\mathcal{H}$ by $\frac{1}{2i\pi}\frac{\mathrm{d}}{\mathrm{d}z} = q\frac{\mathrm{d}}{\mathrm{d}q}$. On $q$-expansions, this operator maps $\sum\limits_{n \geq 0} a_{n}q^{n}$ to $\sum\limits_{n \geq 0} na_{n}q^{n}$. It is well-known that if $g$ is a modular form, then $\theta g$ is no longer modular (see for example \cite[Chapter 5]{Zagier08}). However, Swinnerton-Dyer and Serre \cite[\nopp § 3]{Swinnerton-Dyer73} have proved that for $\ell \geq 5$ and a level $1$ modular form $g$ with $\mathfrak{L}$-integral Fourier coefficients, one can construct a level $1$ form with $\mathfrak{L}$-integral Fourier coefficients and which Fourier coefficients are congruent modulo $\mathfrak{L}$ to the one of $\theta g$. More generally, Katz \cite{Katz77} has proved using his geometric theory of modular forms that there is an operator on the space of modular forms with coefficients in an algebraic closure of $\F_{\ell}$, which action on the $q$-expansions is the same as the one of $\theta$. For our purposes, the main drawbacks of this latter approach is that Katz' modular forms modulo $\mathfrak{L}$ do not always lift in characteristic $0$, and have by essence a prime-to-$\ell$ level. To remedy this, we will construct for any given level $N \geq 1$ and place $\mathfrak{L}$, an operator $\tilde{\theta}$ acting on $\Mod(N)$, stabilising the subspace of forms with $\mathfrak{L}$-integral Fourier coefficients, and such that for every modular form $g$ with $\mathfrak{L}$-integral coefficients we have
\[\tilde{\theta}g \equiv \theta g \pmod{\mathfrak{L}},\]
meaning that $a_{n}\big(\tilde{\theta}g\big)$ and $n a_{n}(g)$ are congruent modulo $\mathfrak{L}$ for all $n$.

The main tool we will use in the construction of $\tilde{\theta}$ is the Rankin-Cohen bracket, introduced by Cohen in \cite[Corollary 7.2]{Cohen75}. We recall its definition and properties below.

\begin{prop}[Rankin-Cohen bracket]\label{RCbracket}
Let $g$ and $h$ be two modular forms of weight $k_{g}$ and $k_{h}$, level $N_{g}$ and $N_{h}$, and character $\varepsilon_{g}$ and $\varepsilon_{h}$ respectively. The Rankin-Cohen bracket of $g$ and $h$ is
\[[g,h] := k_{g}g \theta h - k_{h} h \theta g.\]
It is a modular form of weight $k_{g} + k_{h} + 2$, level $\lcm(N_{g},N_{h})$ and character $\varepsilon_{g}\varepsilon_{h}$. Moreover, if both $g$ and $h$ have their Fourier coefficients in a ring $R$, then so has $[g,h]$.
\end{prop}
Let $N$ be a positive integer. For a prime number $p$, we denote by $T_{p}^{N}$ the $p$-th Hecke operator acting on $\Mod(N)$. Recall that a modular form $g \in \Mod_{k}(N,\varepsilon)$ is an eigenform for $T_{p}^{N}$ modulo $\mathfrak{L}$ with eigenvalue $a_{p} \in \bar{\F}_{\ell}$ in the sense of \cite[\nopp §6(b)]{Deligne74} if $g$ has $\mathfrak{L}$-integral Fourier coefficients, and if $T_{p}^{N}g$ is congruent to $a_{p} g$ modulo $\mathfrak{L}$. If $g$ is moreover normalised modulo $\mathfrak{L}$, that is if $a_{1}(g) \equiv 1 \pmod{\mathfrak{L}}$, then $g$ is an eigenform for $T_{p}^{N}$ modulo $\mathfrak{L}$ if and only if for all integer $n \geq 0$ prime to $p$, and all $\alpha \geq 1$, we have
\[\left\{\begin{array}{l}
a_{np^{\alpha}}(g) \equiv a_{n}(g)a_{p^{\alpha}}(g) \pmod{\mathfrak{L}};\\
a_{p^{\alpha+1}}(g) \equiv a_{p}(g)a_{p^{\alpha}}(g) - p^{k-1}\varepsilon(p)a_{p^{\alpha-1}}(g) \pmod{\mathfrak{L}}.
\end{array}\right.\]
The eigenvalue is then moreover the reduction of $a_{p}(g)$ modulo $\mathfrak{L}$. The following lemma is the central result that shows how to construct an operator $\tilde{\theta}$ satisfying the properties described above, using Rankin-Cohen brackets.

\begin{lm}\label{theta_lm}
Let $k_{A}$ be a positive integer, and let $\chi_{A}$ be a Dirichlet character modulo $N$. Let $A \in \Mod_{k_{A}}(N,\chi_{A})$ be such that $A$ and $\frac{1}{k_{A}}\theta A$ have $\mathfrak{L}$-integral Fourier coefficients and satisfy
\begin{equation}\label{theta_cond}
A \equiv 1 \pmod{\mathfrak{L}}, \quad \frac{1}{k_{A}}\theta A \equiv 0 \pmod{\mathfrak{L}}, \quad \text{and} \quad \chi_{A} \equiv \bar{\chi}_{\ell}^{-k_{A}} \pmod{\mathfrak{L}}.
\end{equation}
Then, we have a well-defined operator $\tilde{\theta}_{A}$ on $\Mod(N)$ given by $\tilde{\theta}_{A}g = -\frac{1}{k_{A}}[g,A]$. For every $g \in \Mod_{k}(N,\varepsilon)$ with $\mathfrak{L}$-integral Fourier coefficients, this operator satisfies the following properties.
\begin{itemize}
\item $\tilde{\theta}_{A}g \in \Mod_{k+k_{A}+2}(N,\varepsilon\chi_{A})$ and has $\mathfrak{L}$-integral Fourier coefficients;
\item $\tilde{\theta}_{A}g \equiv \theta g \pmod{\mathfrak{L}}$;
\item Moreover, if for some prime number $p$, $g$ is a normalised eigenform for $T_{p}^{N}$ modulo $\mathfrak{L}$ then $\tilde{\theta}_{A}g$ is also a normalised eigenform for $T_{p}^{N}$ modulo $\mathfrak{L}$, with eigenvalue $pa_{p}(g)$.
\end{itemize}
In the following, when there will be no confusion on the form $A$, we shall write $\tilde{\theta}$ for the operator $\tilde{\theta}_{A}$.
\end{lm}
\begin{pf}
According to \cref{RCbracket} above, this is clear that $\tilde{\theta}_{A}$ is a well-defined operator and that $\tilde{\theta}_{A}g$ has the announced weight, level, and character. Furthermore, we have $\tilde{\theta}_{A}g = -\frac{k}{k_{A}}g\theta A + A\theta g$. Therefore, from the assumption, if $g$ has $\mathfrak{L}$-integral Fourier coefficients, then so does $\tilde{\theta}_{A}g$ and we have
\[\tilde{\theta}_{A}g = A\theta g - \frac{k}{k_{A}}g \theta A \equiv \theta g \pmod{\mathfrak{L}}.\]
Assume $g$ is a normalised eigenform for $T_{p}^{N}$ modulo $\mathfrak{L}$. Then, $\tilde{\theta}_{A}g$ is also normalised modulo $\mathfrak{L}$ because we have $a_{1}\left(\tilde{\theta}_{A}g\right) \equiv 1 \times a_{1}(g) \equiv 1 \pmod{\mathfrak{L}}$. Let $n \geq 0$ be prime to $p$, and $\alpha \geq 1$. We have
\[a_{np^{\alpha}}\left(\tilde{\theta}_{A}g\right) \equiv np^{\alpha}a_{np^{\alpha}}(g) \equiv \left(na_{n}(g)\right)\left(p^{\alpha}a_{p^{\alpha}}(g)\right) \equiv a_{n}\left(\tilde{\theta}_{A}g\right)a_{p^{\alpha}}\left(\tilde{\theta}_{A}g\right) \pmod{\mathfrak{L}},\]
\begin{align*}
a_{p^{\alpha+1}}\left(\tilde{\theta}_{A}g\right) \equiv p^{\alpha + 1}a_{p^{\alpha + 1}}(g) &\equiv p^{\alpha+1}\left(a_{p}(g)a_{p^{\alpha}}(g) - p^{k-1}\varepsilon(p)a_{p^{\alpha-1}}(g)\right) \pmod{\mathfrak{L}}\\
	&\equiv a_{p}\left(\tilde{\theta}_{A}g\right)a_{p^{\alpha}}\left(\tilde{\theta}_{A}g\right) - p^{k+1}\varepsilon(p)a_{p^{\alpha-1}}\left(\tilde{\theta}_{A}g\right) \pmod{\mathfrak{L}}.
\end{align*}
If $p \mid N\ell$, we have $p^{k+1}\varepsilon(p) \equiv 0 \equiv p^{(k+k_{A}+2)-1}\varepsilon\chi_{A}(p) \pmod{\mathfrak{L}}$. Otherwise, we have $p^{k_{A}}\chi_{A}(p) \equiv 1 \pmod{\mathfrak{L}}$ by assumption, and we again have $p^{k+1}\varepsilon(p) \equiv p^{(k+k_{A}+2)-1}\varepsilon\chi_{A}(p) \pmod{\mathfrak{L}}$. As desired, we therefore have
\[a_{p^{\alpha+1}}\left(\tilde{\theta}_{A}g\right) \equiv a_{p}\left(\tilde{\theta}_{A}g\right)a_{p^{\alpha}}\left(\tilde{\theta}_{A}g\right) - p^{(k+k_{A}+2)-1}\varepsilon\chi_{A}(p)a_{p^{\alpha-1}}\left(\tilde{\theta}_{A}g\right) \pmod{\mathfrak{L}}.\]
\end{pf}
\begin{Rk}
When $\ell$ does not divide $N$, the reduction of $A$ modulo $\mathfrak{L}$ is the so-called Katz' Hasse invariant.
\end{Rk}
The rest of this paragraph is devoted to construct, for each level $N$ and place $\mathfrak{L}$, a form $A$ that satisfies the hypotheses of \cref{theta_lm} that we will constantly use. Among all possible forms, the ones presented in table \ref{theta_table} are those we found with the smallest weight. Notice that if we have a form $A$ of level $M$ satisfying the hypotheses of \cref{theta_lm} for a given place $\mathfrak{L}$, this form also satisfies the hypotheses of \cref{theta_lm} at the multiple-of-$M$ levels. We will use this fact to consider the smallest set of level possible.

\subsubsection{Theta operators in characteristic greater than 3}\label{theta5}

The following proposition was already known to Swinnerton-Dyer in \cite[Theorem 2]{Swinnerton-Dyer73}.

\begin{prop}\label{theta5_coprime}
Assume $\ell \geq 5$. The form $A := -2\ell E_{\ell-1} \in \Mod_{\ell-1}(1,\mathds{1})$ satisfies the hypotheses of \cref{theta_lm} for any level $N$.
\end{prop}
\begin{pf}
Since $\ell \geq 5$, $A$ is well-defined and the constant coefficient of $A$ is equal to $\frac{\ell B_{\ell-1,\mathds{1}}}{\ell-1}$. From \cref{Ber_divisors}, it is $\mathfrak{L}$-integral and congruent to $1$ modulo $\mathfrak{L}$. Moreover, because $E_{\ell-1}$ has integral coefficients, away from the constant one, it follows that $A$ and $-\frac{1}{k_{A}}\theta A$ have $\mathfrak{L}$-integral Fourier coefficients and that $A \equiv 1 \pmod{\mathfrak{L}}$ and $\frac{1}{k_{A}}\theta A \equiv 0 \pmod{\mathfrak{L}}$. Finally we have $\bar{\chi}_{\ell}^{-k_{A}} = \bar{\chi}_{\ell}^{1-\ell} \equiv \mathds{1} \pmod{\mathfrak{L}}$.
\end{pf}

If the level $N$ is divisible by $\ell$, the situation is in fact much more pleasant for us, in the sense that we can find a form with $k_{A} = 1$. We find a record of the following fact in \cite[(2.1) Theorem]{Ribet94}.

\begin{prop}\label{theta5_level}
Assume $\ell \geq 5$ and $\ell \mid N$. Let $\chi_{\mathfrak{L}}$ be the Teichmüller lift of $\bar{\chi}_{\ell}$ with respect to the place $\mathfrak{L}$, viewed as a primitive Dirichlet modulo $\ell$. The form $A := 2\ell E_{1}^{\mathds{1},\chi_{\mathfrak{L}}^{-1}} \in \Mod_{1}(\ell,\chi_{\mathfrak{L}}^{-1})$ satisfies the hypotheses of \cref{theta_lm}.
\end{prop}
\begin{pf}
The form $A$ is well-defined because $\chi_{\mathfrak{L}}^{-1}$ is an odd character. Indeed, we have $\chi_{\mathfrak{L}}^{-1}(-1) \equiv \bar{\chi}_{\ell}^{-1}(-1) \equiv -1 \pmod{\mathfrak{L}}$, and because $\ell$ is odd, this lifts to $\chi_{\mathfrak{L}}^{-1}(-1) = -1$. The constant term of $A$ is equal to $-\ell B_{1,\chi_{\mathfrak{L}}^{-1}} = -\sum\limits_{i = 1}^{\ell-1} i\chi_{\mathfrak{L}}^{-1}(i)$ which is $\mathfrak{L}$-integral. Therefore, because $\chi_{\mathfrak{L}}$ induces the identity modulo $\mathfrak{L}$, this coefficient is congruent to $1$ modulo $\mathfrak{L}$, and because $E_{1}^{\mathds{1},\chi_{\mathfrak{L}}^{-1}}$ has integral coefficients away from the constant one, $A$ and $\frac{1}{k_{A}}\theta A$ have $\mathfrak{L}$-integral Fourier coefficients. Moreover, we also get $A \equiv 1 \pmod{\mathfrak{L}}$ and $\frac{1}{k_{A}}\theta A \equiv 0 \pmod{\mathfrak{L}}$. Finally, by definition we have $\bar{\chi}_{\ell}^{-k_{A}} = \bar{\chi}_{\ell}^{-1} \equiv \chi_{A} \pmod{\mathfrak{L}}$. Thus, $A$ satisfies the hypotheses of \cref{theta_lm}.
\end{pf}
This finishes the case $\ell \geq 5$. For $\ell \leq 3$, the two previous constructions do no always give well-defined modular forms. We present in the two next paragraph specific constructions in the cases $\ell = 2$ and $\ell = 3$.

\subsubsection{Theta operators in characteristic 2}\label{theta2}

For $\ell = 2$, the most favourable case is when $4 \mid N$. The following construction is very analogous to the one of \cref{theta5_level}.

\begin{prop}\label{theta2_even}
Assume $\ell = 2$ and $N$ is divisible by $4$. Let $\chi_{4}$ be the only non-trivial Dirichlet character modulo $4$. The form $A := 4E_{1}^{\mathds{1},\chi_{4}} \in \Mod_{1}(4,\chi_{4})$ satisfies the hypotheses of \cref{theta_lm}.
\end{prop}
\begin{pf}
The form $A$ is well-defined because $\chi_{4}$ is odd. Moreover, the constant coefficient of $A$ is equal to $-2B_{1,\chi_{4}} = -\frac{1}{2} \left(1\chi_{4}(1)+3\chi_{4}(3)\right) = 1$. Therefore, the constant coefficient of $A$ is equal to $1$ and because $E_{1}^{\mathds{1},\chi_{4}}$ has integral coefficients away from the constant one, $A$ and $\frac{1}{k_{A}}\theta A$ have $\mathfrak{L}$-integral coefficients. Moreover, we also get $A \equiv 1 \pmod{\mathfrak{L}}$ and $\frac{1}{k_{A}}\theta A \equiv 0 \pmod{\mathfrak{L}}$. Finally, it is straightforward that $\chi_{4}$ is trivial modulo $\mathfrak{L}$, as is the cyclotomic character modulo $2$.
\end{pf}

The next favourable case is when $N$ admits at least one odd prime divisor. The following result was inspired by \cite[Appendix A.]{Meier17}. As it has never been published, we prove it for the sake of completeness.

\begin{prop}\label{theta2_odd}
Assume $\ell = 2$ and $N$ has an odd prime divisor. Let $p$ be the least odd prime divisor of $N$, and let $\chi_{N}$ be a Dirichlet character modulo $p$ of order $2^{m}$, the greatest power of $2$ dividing $p-1$. Let $\zeta$ be any $2^{m}$-th root of unity. The form $A := (\zeta-1)E_{1}^{\mathds{1},\chi_{N}} \in \Mod_{1}(p,\chi_{N})$ satisfies the hypotheses of \cref{theta_lm}.
\end{prop}
\begin{pf}
Let $g$ be an integer generating $\left(\Z/p\Z\right)^{\times}$ and such that $\chi_{N}(g) = \zeta$. Because $\zeta$ is a root of unity of order $2^{m}$, we have $\chi_{N}(-1) = \chi_{N}(g)^{\frac{p-1}{2}} = -1$. Therefore, $\chi_{N}$ is odd and $A$ is well-defined and its constant coefficient is equal to $\frac{1-\zeta}{2}B_{1,\chi_{N}} = \frac{1-\zeta}{2p}\sum\limits_{a = 1}^{p-1} a\chi_{N}(a)$.

For $i \in \ent{0}{2^{m-1}-1}$, we have $\chi_{N}(g^{i}) = \zeta^{i}$, and $\chi_{N}(-g^{i}) = -\zeta^{i} = \zeta^{i+2^{m-1}} = \chi_{N}(g^{i+2^{m-1}})$. Therefore, the set $\left\{\pm g^{i}, i \in \ent{0}{2^{m-1}-1}\right\}$ is a set of representatives of $\left(\Z/p\Z\right)^{\times}/\mathrm{Ker}(\chi_{N})$. For an integer $x$, we write $[x]$ for the only integer between $0$ and $p-1$ that is congruent to $x$ modulo $p$. We then have
\begin{align*}
\frac{1-\zeta}{2}B_{1,\chi_{N}} &= \frac{1-\zeta}{2p} \sum\limits_{e \in \mathrm{Ker}(\chi_{N})}\sum\limits_{i = 0}^{2^{m-1}-1} \left(\left[eg^{i}\right]\chi_{N}(g^{i}) + \left[-eg^{i}\right]\chi_{N}(-g^{i})\right)\\
    &= \frac{1-\zeta}{2p}\sum\limits_{i = 0}^{2^{m-1}-1} \sum\limits_{e \in \mathrm{Ker}(\chi_{N})} \left(\left[eg^{i}\right]\zeta^{i}  - \left(p - \left[eg^{i}\right]\right)\zeta^{i} \right)\\
    &= \frac{1-\zeta}{2p} \sum\limits_{i = 0}^{2^{m-1}} \zeta^{i} \left(-p \cdot \#\mathrm{Ker}(\chi_{N}) + 2\sum\limits_{e \in \mathrm{Ker}(\chi_{N})} \left[eg^{i}\right]\right)\\
    &= -\frac{p-1}{2^{m}} + (1-\zeta)\left(\frac{1}{p}\sum\limits_{i = 0}^{2^{m-1}-1}\zeta^{i}\sum\limits_{e \in \mathrm{Ker}(\chi_{N})} \left[eg^{i}\right]\right).
\end{align*}
The term inside the parentheses is $\mathfrak{L}$-integral and $\frac{p-1}{2^{m}}$ is odd. Moreover, the only prime ideal above $2$ in the ring $\mathcal{O}_{\Q(\zeta)} = \Z[\zeta]$ is $(1-\zeta)\Z[\zeta]$. Therefore, we have $\frac{1-\zeta}{2}B_{1,\chi_{N}} \equiv 1 \pmod{\mathfrak{L}}$. Because the non-constant Fourier coefficients of $E_{1}^{\mathds{1},\chi_{N}}$ are integral, $A$ and $\frac{1}{k_{A}}\theta A$ have $\mathfrak{L}$-integral coefficients, and we get $A \equiv 1 \pmod{\mathfrak{L}}$ and $\frac{1}{k_{A}}\theta A \equiv 0 \pmod{\mathfrak{L}}$. Finally, from \cref{Rootmod}, because $\chi_{N}$ has order a power of $2$, it is trivial modulo $\mathfrak{L}$ as well as the cyclotomic character modulo $2$. This finishes the proof.
\end{pf}

We are left with the cases $N = 1$ and $N = 2$. There is no modular form of weight 1 of these levels, so we have to look at bigger weights in order to construct the form $A$. For level $2$, we show that weight 2 suffices.

\begin{prop}\label{theta2_2}
Assume that $\ell = 2$ and $N = 2$. Let $\mathds{1}_{(2)}$ be the trivial character modulo $2$. The modular form $A := 24E_{2}^{\mathds{1},\mathds{1}_{(2)}} \in \Mod_{2}(2,\mathds{1}_{(2)})$ satisfies the hypotheses of \cref{theta_lm}.
\end{prop}
\begin{pf}
The constant coefficient of $A$ is equal to $1$ and $E_{2}^{\mathds{1},\mathds{1}_{(2)}}$ has integral coefficients away from the constant one. Therefore, the forms $A$ and $\frac{1}{k_{A}}\theta A$ have both $\mathfrak{L}$-integral Fourier coefficients, and we have $A \equiv 1 \pmod{\mathfrak{L}}$ and $\frac{1}{k_{A}}\theta A \equiv 0 \pmod{\mathfrak{L}}$. Finally, the character of $A$ is trivial modulo $\mathfrak{L}$ as well as the cyclotomic character modulo $2$.
\end{pf}

For $N = 1$, the weight needs to be at least $4$ and we have the following result.

\begin{prop}\label{theta2_1}
Assume $\ell = 2$ and $N = 1$. The form $A := 240E_{4} \in \Mod_{4}(1,\mathds{1})$ satisfies the hypotheses of \cref{theta_lm}.
\end{prop}
\begin{pf}
The constant coefficient of $A$ is equal to $1$ and the non-constant Fourier coefficients of $E_{4}$ are integer. Therefore, $A$ and $\frac{1}{k_{A}}\theta A$ have integer coefficients and we have $A \equiv 1 \pmod{\mathfrak{L}}$, and $\frac{1}{k_{A}}\theta A \equiv 0 \pmod{\mathfrak{L}}$. The character of $A$ is again trivial, as well as the cyclotomic character modulo $2$.
\end{pf}

\subsubsection{Theta operators in characteristic 3}\label{theta3}

For $N$ divisible by $3$, the form of \cref{theta5_level} is still appropriate.
\begin{prop}\label{theta3_level}
Assume $\ell = 3$ and $N$ is divisible by $3$. Let $\chi_{3}$ be the unique non-trivial Dirichlet character modulo $3$. The form $A := 6E_{1}^{\mathds{1},\chi_{3}} \in \Mod_{1}(3,\chi_{3})$ satisfies the hypotheses of \cref{theta_lm}.
\end{prop}
\begin{pf}
We have $\chi_{3} \equiv \bar{\chi}_{\ell}^{-1} \pmod{\mathfrak{L}}$ and the proof is exactly the same as the one of \cref{theta5_level}.
\end{pf}

For the levels containing a prime divisor congruent to $2$ modulo $3$, we can still consider an Eisenstein series for the form $A$.

\begin{prop}\label{theta3_2mod3}
Assume $\ell = 3$ and $N$ has a prime divisor congruent to $2$ modulo $3$. Let $p$ be the least such prime divisor, and let $\mathds{1}_{(p)}$ be the trivial Dirichlet character modulo $p$. The form $A := \frac{24}{p-1}E_{2}^{\mathds{1},\mathds{1}_{(p)}} \in \Mod_{2}(p,\mathds{1}_{(p)})$ satisfies the hypotheses of \cref{theta_lm}.
\end{prop}
\begin{pf}
The constant coefficient of $A$ is equal to $1$. Moreover, $E_{2}^{\mathds{1},\mathds{1}_{(p)}}$ has integral Fourier coefficients away from the constant one. Because, $p$ is congruent to $2$ modulo $3$, $\frac{24}{p-1}$ is $0$ modulo $\mathfrak{L}$. Therefore $A$ and $\frac{1}{k_{A}}\theta A$ have $\mathfrak{L}$-integral Fourier coefficients, and we have $A \equiv 1 \pmod{\mathfrak{L}}$ and $\frac{1}{k_{A}}\theta A \equiv 0 \pmod{\mathfrak{L}}$. Finally, the character of $A$ is trivial, and we have $\bar{\chi}_{3}^{-2} = \mathds{1}$.
\end{pf}

The remaining cases are the levels containing only prime factors that are congruent to 1 modulo 3. For the levels divisible by a prime $p$ congruent to 4 modulo 9 (that is if 3 divides $p-1$ only once), we found the following construction.

\begin{prop}\label{theta3_4mod9}
Assume $\ell = 3$ and $N$ has a prime divisor congruent to $4$ modulo $9$. Let $p$ be the least such prime divisor of $N$, and let $\chi^{N}$ be a Dirichlet character modulo $p$ of order 3. The modular form $A := \frac{3}{p-1} \left(E_{2}^{\mathds{1},\chi^{N}} - E_{2}^{\chi^{N},\mathds{1}}\right) \in \Mod_{2}\left(p,\chi^{N}\right)$ satisfies the hypotheses of \cref{theta_lm}.
\end{prop}
\begin{pf}
First notice that $\chi^{N}$ indeed exists as $\left(\Z/p\Z\right)^{\times}$ is a cyclic group of order $p-1$ that is divisible by $3$. Moreover, $\chi^{N}$ has order $3$. Therefore, it is trivial modulo $\mathfrak{L}$ and even, and the two Eisenstein series $E_{2}^{\mathds{1},\chi^{N}}$ and $E_{2}^{\chi^{N},\mathds{1}}$ exist.

The constant coefficient of $A$ is equal to $\frac{3}{4(1-p)}B_{2,\chi^{N}}$ which is $\mathfrak{L}$-integral by \cref{Ber_divisors}. We have
\[\frac{3}{4(1-p)} B_{2,\chi^{N}} = \frac{3}{4p(1-p)} \sum\limits_{a = 1}^{p-1} a^{2}\chi^{N}(a) \equiv \frac{3}{1-p} \sum\limits_{a = 1}^{p-1} a^{2} \equiv \frac{3}{1-p} \frac{p(p-1)(2p-1)}{6} \equiv 1 \pmod{\mathfrak{L}}.\]
Therefore, the constant coefficient of $A$ is $1$ modulo $\mathfrak{L}$, and because the non-constant Fourier coefficients of $E_{2}^{\mathds{1},\chi^{N}}$ and $E_{2}^{\chi^{N},\mathds{1}}$ are integral, $A$ and $\frac{1}{k_{A}}\theta A$ have $\mathfrak{L}$-integral coefficients. The weight $k_{A}$ is invertible modulo $3$, it consequently suffices to prove that $A \equiv 1 \pmod{\mathfrak{L}}$ to conclude.

The forms $E_{2}^{\mathds{1},\chi^{N}}$ and $E_{2}^{\chi^{N},\mathds{1}}$ are both normalised eigenforms for all the Hecke operators at level $p$, and have the same weight and character. Thus, it is enough to prove that $a_{r}\left(E_{2}^{\mathds{1},\chi^{N}}\right) \equiv a_{r}\left(E_{2}^{\chi^{N},\mathds{1}}\right) \pmod{\mathfrak{L}}$ for all prime numbers $r$. This last congruence is straightforward, because we have
\[a_{r}\left(E_{2}^{\mathds{1},\chi^{N}}\right) = 1+r\chi^{N}(r) \equiv \chi^{N}(r) + r = a_{r}\left(E_{2}^{\chi^{N},\mathds{1}}\right) \pmod{\mathfrak{L}}.\]
\end{pf}

It only remains the levels containing only primes congruent to 1 modulo 9. We found no general way to express the modular form $A$ as form of weight $2$. Using computations in PARI/GP, we look for a modular form of level $p \equiv 1 \pmod{9}$ satisfying the hypotheses of \cref{theta_lm} for $p$ up to $1000$. We always find a form except for $p \in \{307,379,433,487,523,613,631,757,811,829,991\}$, \textit{i.e.} we found 16 forms out of the 27 we were looking for. It can be proved that such a modular form cannot be expressed as a linear combination of forms in the Eisenstein space, meaning that one has necessarily to consider cusp forms to construct $A$. To fill this gap anyway, we can still consider the modular form $A := 240E_{4}$ as in the case of \cref{theta2_1}.

\begin{prop}\label{theta3_bad}
Assume $\ell = 3$, and $N$ contains only prime factors congruent to $1$ modulo $9$. The modular form $A := 240E_{4}$ satisfies the hypotheses of \cref{theta_lm}.
\end{prop}
\begin{pf}
The constant coefficient of $A$ is equal to $1$ and the non-constant Fourier coefficients of $E_{4}$ are integral. Therefore, the forms $A$ and $\frac{1}{k_{A}}\theta A$ have both $\mathfrak{L}$-integral Fourier coefficients, and we have $A \equiv 1 \pmod{\mathfrak{L}}$ and $\frac{1}{k_{A}}\theta A \equiv 0 \pmod{\mathfrak{L}}$. Finally, the character of $A$ is trivial and we have $\bar{\chi}_{\ell}^{-k_{A}} = \bar{\chi}_{\ell}^{-4} \equiv \mathds{1} \pmod{\mathfrak{L}}$.
\end{pf}
We have compiled in \cref{theta_table} the definition of $A$ depending on $\ell$ and $N$. When multiple definitions were possible, we have taken the one with the least weight among all the possible forms. The third column corresponds to the proposition where the properties of the form have been proved.
\begin{table}[!ht]
\caption{}
\begin{minipage}{.45\linewidth}
\begin{center}
\begin{tabular}{|c|c|c|}
\hline
$\ell \geq 5$ & Form $A$ & Proposition\\
\hline
$\ell \nmid N$ & $\Big.-2\ell E_{\ell-1}\Big.$ & \cref{theta5_coprime}\\
\hline
$\ell \mid N$ & $\bigg.2\ell E_{1}^{\mathds{1},\chi_{\mathfrak{L}}^{-1}}\bigg.$ & \cref{theta5_level}\\
\hline
\end{tabular}
\end{center}
\end{minipage}
\begin{minipage}{.5\linewidth}
\begin{center}
\begin{tabular}{|c|c|c|}
\hline
$\ell = 2$ & Form $A$ & Proposition\\
\hline
$4 \mid N$ & $\Big.4E_{1}^{\mathds{1},\chi_{4}}\Big.$ & \cref{theta2_even}\\
\hline
$N \geq 3$ and $4 \nmid N$ & $\Big.(\zeta-1)E_{1}^{\mathds{1},\chi_{N}}\Big.$ & \cref{theta2_odd}\\
\hline
$N = 2$ & $\Big.24E_{2}^{\mathds{1},\mathds{1}_{(2)}}\Big.$ & \cref{theta2_2}\\
\hline
$N = 1$ & $\Big.240E_{4}$ & \cref{theta2_1} \\
\hline
\end{tabular}
\end{center}
\end{minipage}

\begin{center}
\begin{tabular}{|c|c|c|}
\hline
$\ell = 3$ & Form $A$ & Proposition\\
\hline
$\ell \mid N$ & $\bigg.6E_{1}^{\mathds{1},\chi_{3}}\bigg.$ & \cref{theta3_level}\\
\hline
$\ell \nmid N$ and $N$ has a prime & \multirow{2}{*}{$\Big.\frac{24}{p-1} E_{2}^{\mathds{1},\mathds{1}_{(p)}}\Big.$} & \multirow{2}{*}{\cref{theta3_2mod3}}\\
factor $q \equiv 2 \pmod{3}$ & & \\
\hline
$\forall d \mid N$, $d \equiv 1 \pmod{3}$ and $N$ has & \multirow{2}{*}{$\frac{3}{p-1}\left(E_{2}^{\mathds{1},\chi^{N}}-E_{2}^{\chi^{N},\mathds{1}}\right)$} & \multirow{2}{*}{\cref{theta3_4mod9}}\\
a prime factor $p \equiv 4 \pmod{9}$ & & \\
\hline
$\forall p \mid N$, $p \equiv 1 \pmod{9}$ & $\Big.240E_{4}$ & \cref{theta3_bad}\\
\hline
\end{tabular}
\end{center}
\label{theta_table}
\end{table}

Looking at the various results above, we state the following definition that will be useful in the proofs of the next paragraph.
\begin{df}\label{bad_df}
We say a pair $(\ell,N)$ is bad, if we have one of the following.
\begin{itemize}
\item $\ell = 2$ and $N = 1$;
\item $\ell = 3$ and all the prime factors of $N$ are congruent to $1$ modulo $9$.
\end{itemize}
\end{df}

\begin{Rk}
When $(\ell,N)$ is bad, the modular $-504E_{6}$ is also congruent to $1$ modulo $\mathfrak{L}$. Its weight is greater than the one of \cref{theta_table} but we will have to use it in the proof of \cref{Sturm_theta} in the next section.
\end{Rk}

\subsection{Sturm bounds}\label{sec_Sturm}
A Sturm bound for a space of modular forms is an upper bound on the number of leading coefficients that characterise a form of this space. Equivalently, it is the maximal number of zero leading coefficients that a non-zero form of this space can have. The study of such bounds has first been made by Sturm \cite{Sturm87} and was later generalised among others by Murty \cite{Murty97}. The same kind of bounds exists if we look at modular forms modulo a prime ideal -- and are in fact the same as the first one. In the next lemma we give a slight improvement of Murty's result for modular forms of same weight. We then state a more general result for modular forms of any weight and level.

For all the paragraph, we fix a prime number $\ell$ and a place  $\mathfrak{L}$ of $\bar{\Q}$ above $\ell$.

\begin{lm}\label{Sturm}
Let $f$, $g$ be two modular forms of same weight $k \geq 0$, level $N_{f}$, $N_{g}$ and character $\varepsilon_{f}$, $\varepsilon_{g}$ respectively. Let $N$ be the lcm of $N_{f}$ and $N_{g}$. Assume that $f$ and $g$ have both $\mathfrak{L}$-integral Fourier coefficients and that $\varepsilon_{f} \equiv \varepsilon_{g} \pmod{\mathfrak{L}}$.

If $a_{n}(f) \equiv a_{n}(g) \pmod{\mathfrak{L}}$ for every integer $n \leq \frac{kN}{12} \prod\limits_{\underset{p \text{ prime}}{p \mid N}} \left(1 + \frac{1}{p}\right)$, then $f \equiv g \pmod{\mathfrak{L}}$.
\end{lm}

\begin{pf}
We follow substantially the proof of Murty of \cite[Paragraph 5]{Murty97}. Consider $\phi = f-g$. As explained in \cite{Murty97}, for $\gamma \in \SL_{2}(\Z)$ there is an element $A_{\gamma} \in \bar{\Q}^{\times}$ such that the modular form $A_{\gamma} \phi|_{k}\gamma$ has $\mathfrak{L}$-integral coefficients and is not congruent to $0$ modulo $\mathfrak{L}$.

We write $m := [\SL_{2}(\Z):\Gamma_{0}(N)] = N\prod\limits_{p \mid N} \left(1+\frac{1}{p}\right)$ and consider a system of representatives $(\gamma_{i})_{i \in \ent{1}{m}}$ of right cosets of $\Gamma_{0}(N)$ in $\SL_{2}(\Z)$. We can further assume $\gamma_{1} = I_{2}$, the identity matrix. Also choose a set $(\tau_{j})_{j \in \ent{1}{\varphi(N)}}$ of representatives of $\Gamma_{1}(N)$ in $\Gamma_{0}(N)$. We then have,
\[\SL_{2}(\Z) = \bigcup_{i = 1}^{m} \Gamma_{0}(N) \gamma_{i} = \bigcup_{i = 1}^{m} \bigcup_{j = 1}^{\varphi(N)} \Gamma_{1}(N) \tau_{j} \gamma_{i}.\]
Taking the norm function of $\phi$ according to this system of representatives, we get
\[F := \prod_{i = 1}^{m} \prod_{j = 1}^{\varphi(N)} A_{\tau_{j}\gamma_{i}} \phi|_{k} \tau_{j}\gamma_{i} \in \Mod_{km\varphi(N)}(\SL_{2}(\Z)).\]
The form $F$ has $\mathfrak{L}$-integral Fourier coefficients by construction. We want to show that the vanishing order of $F$ modulo $\mathfrak{L}$ is at least $\frac{km\varphi(N)}{12}$. For $i = 1$ and $j \in \ent{1}{\varphi(N)}$, we have
\[A_{\tau_{j}\gamma_{1}}\phi|_{k}\tau_{j}\gamma_{1} = A_{\tau_{j}}\phi|_{k}\tau_{j} = A_{\tau_{j}}\left(\varepsilon_{f}(\tau_{j})f - \varepsilon_{g}(\tau_{j})g\right) \equiv A_{\tau_{j}}\varepsilon_{f}(\tau_{j}) \phi \pmod{\mathfrak{L}}.\]
By assumption, the vanishing order at infinity of $\phi$ modulo $\mathfrak{L}$ is at least equal to $\frac{km}{12}$. We conclude that the modular form $\prod\limits_{j=1}^{\varphi(N)} A_{\tau_{j}} \varphi|_{k}\tau_{j}$ has a zero modulo $\mathfrak{L}$ at infinity of order at least equal to $\frac{km\varphi(N)}{12}$, and therefore the same is true for $F$. Applying Sturm's theorem for level $1$ modular forms \cite[Theorem 5]{Murty97}, $F$ must vanish modulo $\mathfrak{L}$ and by construction of the coefficients $A_{\tau_{j}\gamma_{i}}$, $\phi$ must also be trivial modulo $\mathfrak{L}$.
\end{pf}

The following proposition generalises the previous lemma to modular forms of arbitrary weights and levels. The proof uses extensively the construction of theta operators given in \cref{sec_theta}.
\begin{prop}\label{Sturm_theta}
Let $f$, $g$ be two modular forms of weights $k_{f}$, $k_{g} \geq 0$, levels $N_{f}$, $N_{g} \geq 1$ and characters $\varepsilon_{f}$, $\varepsilon_{g}$ respectively. Let $m_{f}$, $m_{g}$ be two non-negative integers. Assume that $f$ and $g$ have both $\mathfrak{L}$-integral Fourier coefficients and that $\bar{\chi}_{\ell}^{k_{f}+2m_{f}}\varepsilon_{f} \equiv \bar{\chi}_{\ell}^{k_{g}+2m_{g}}\varepsilon_{g} \pmod{\mathfrak{L}}$. Let $N$ be the lcm of $N_{f}$ and $N_{g}$, and define
\[a = \left\{\begin{array}{ll}
4 & \text{if } \left|\begin{array}{l}
    k_{f} + 2m_{f} \equiv k_{g} + 2m_{g}+2 \pmod{4}\\
    \text{and } (\ell,N) \text{ is bad}
\end{array}\right.\\
0 & \text{otherwise}
\end{array}\right. b = \left\{\begin{array}{ll}
4 & \text{if } \ell = 2 \text{ and } N = 2\\
3 & \text{if } \ell \mid N \text{ and } (\ell,N) \neq (2,2)\\
6 & \text{if } (\ell,N) \text{ is bad}\\
\ell + 1 & \text{otherwise}
\end{array}\right.\]
\[k = a + \max(k_{f}+bm_{f},k_{g}+bm_{g}),\]
where "bad" refers to \cref{bad_df}.

If for every $n \leq \frac{Nk}{12} \prod\limits_{\underset{p \text{ prime}}{p \mid N}} \left(1+\frac{1}{p}\right)$, we have $n^{m_{f}}a_{n}(f) \equiv n^{m_{g}}a_{n}(g) \pmod{\mathfrak{L}}$ (with $0^{0} = 1$), then this holds for all integers $n \geq 0$.
\end{prop}

\begin{pf}
For the whole proof, we write $A$ for the modular form associated with $\mathfrak{L}$ and $N$ constructed in \cref{sec_theta}. According to \cref{theta_table}, it has weight $b-2$ and level $N$.  We write $\chi_{A}$ for the character of $A$ and $B(N,k) := \frac{Nk}{12} \prod\limits_{\underset{p \text{ prime}}{p \mid N}} \left(1+\frac{1}{p}\right)$.

Assume without loss of generality that $k_{f}+bm_{f} \leq k_{g} + bm_{g}$. We first prove that, assuming $b-2$ divides $k_{g}-k_{f}+b(m_{g}-m_{f})$, we have $n^{m_{f}}a_{n}(f) \equiv n^{m_{g}}a_{n}(g) \pmod{\mathfrak{L}}$ for all non-negative integers $n$ if these congruences hold for $n \leq B(N,k_{g}+bm_{g})$. By applying \cref{theta_lm} recursively, we have
\[\tilde{\theta}^{m_{f}}f \in \Mod_{k_{f}+bm_{f}}(N,\varepsilon_{f}\chi_{A}^{m_{f}}) \quad \text{and} \quad \tilde{\theta}^{m_{g}}g \in  \Mod_{k_{g}+bm_{g}}(N,\varepsilon_{g}\chi_{A}^{m_{g}}).\]
We cannot apply \cref{Sturm} to $\tilde{\theta}^{m_{f}}f$ and $\tilde{\theta}^{m_{g}}g$ directly since they don't have the same weight. However, the functions $A^{\frac{k_{g}-k_{f}+b(m_{g}-m_{f})}{b-2}}\tilde{\theta}^{m_{f}}f$ and $\tilde{\theta}^{m_{g}}g$ are well-defined modular forms by assumption. They have the same weight $k_{g}+bm_{g}$, the same level $N$ and characters $\varepsilon_{f}\chi_{A}^{m_{f}+\frac{k_{g}-k_{f}+b(m_{g}-m_{f})}{b-2}}$ and $\varepsilon_{g}\chi_{A}^{m_{g}}$ respectively. Moreover, by \cref{theta_lm} again, we have $\chi_{A} \equiv \bar{\chi}_{\ell}^{2-b} \pmod{\mathfrak{L}}$. By the assumption on the characters we get
\begin{align*}
\varepsilon_{f}\chi_{A}^{m_{f}+\frac{k_{g}-k_{f}+b(m_{g}-m_{f})}{b-2}} &\equiv \varepsilon_{f} \bar{\chi}_{\ell}^{-(b-2)\left(m_{f}+\frac{k_{g}-k_{f}+b(m_{g}-m_{f})}{b-2}\right)} \pmod{\mathfrak{L}}\\
    &\equiv \bar{\chi}_{\ell}^{k_{f}+2m_{f}}\varepsilon_{f} \cdot \bar{\chi}_{\ell}^{-(k_{g}+bm_{g})} \pmod{\mathfrak{L}}\\
    &\equiv \bar{\chi}_{\ell}^{k_{g}+2m_{g}}\varepsilon_{g} \cdot \bar{\chi}_{\ell}^{-(k_{g}+bm_{g})} \equiv \varepsilon_{g}\chi_{A}^{m_{g}} \pmod{\mathfrak{L}}.
\end{align*}
Therefore, the assumptions of \cref{Sturm} are satisfied for these two modular forms. Since $A$ reduces to $1$ modulo $\mathfrak{L}$, we get that if the coefficients of $\tilde{\theta}^{m_{f}}f$ and $\tilde{\theta}^{m_{g}}g$ are congruent up to the $B(N,k_{g}+bm_{g})$-th one, then $\tilde{\theta}^{m_{f}}f$ and $\tilde{\theta}^{m_{g}}g$ are congruent modulo $\mathfrak{L}$ by \cref{Sturm}.

We now look at the hypothesis $b-2 \mid k_{g}-k_{f}+b(m_{g}-m_{f})$. We claim that if $(\ell,N)$ is not bad, then it is always satisfied. We have three cases: (i) $\ell = N = 2$, (ii) $\ell \mid N$ and $(\ell,N) \neq (2,2)$, (iii) $\ell \nmid N$ and $(\ell,N)$ is not bad.
\begin{enumerate}[(i)]
\item If $\ell = N = 2$, then $b-2 = 2$, and $k_{f} \equiv k_{g} \equiv 0 \pmod{2}$, because the weight of a modular form of level $2$ is necessarily even. Thus, $k_{g}-k_{f}+4(m_{g}-m_{f})$ is divisible by $b-2$.
\item If $\ell \mid N$, then $b-2 = 1$ and there is nothing to prove.
\item If $\ell \nmid N$ and $(\ell,N)$ is not bad, we have $b-2 = \ell - 1$. Because $\ell \nmid N$, $\varepsilon_{f}$ and $\varepsilon_{g}$ are unramified at $\ell$. From the assumption $\bar{\chi}_{\ell}^{k_{f}+2m_{f}}\varepsilon_{f} \equiv \bar{\chi}_{\ell}^{k_{g}+2m_{g}}\varepsilon_{g} \pmod{\mathfrak{L}}$, we get that $k_{g}-k_{f}+2(m_{g}-m_{f}) \equiv 0 \pmod{\ell-1}$, hence $b-2 \mid k_{g}-k_{f}+b(m_{g}-m_{f})$.
\end{enumerate}
Therefore, when $(\ell,N)$ is not bad, the proposition is proved because we have $a = 0$ and $k = k_{g} + bm_{g}$.

From now on, assume that $(\ell,N)$ is bad. By definition, we have $b = 6$, and either $(\ell,N) = (2,1)$, or $\ell = 3$ and $\ell \nmid N$. Let us first prove that we have $k_{g}-k_{f}+6(m_{g}-m_{f}) \equiv 0 \pmod{2}$ (\textit{i.e.} $k_{g} \equiv k_{f} \pmod{2}$). When $(\ell,N) = (2,1)$, it is true because the weights $k_{f}$ and $k_{g}$ are both even. When $\ell = 3$ and $\ell \nmid N$, the hypothesis on the characters again implies that $k_{f}+2m_{f} \equiv k_{g}+2m_{g} \pmod{2}$ and the conclusion follows.

If the even number $k_{g}-k_{f}+6(m_{g}-m_{f})$ is divisible by $4 = b-2$, then by definition we have $a = 0$ and $k = k_{g}+bm_{g}$. The result follows as before in this case. Otherwise, we have $4 \mid k_{g}-k_{f}+6(m_{g}-m_{f})-2$ and $a = 4$. Write $A_{4} := 240E_{4}$ and $A_{6} := -504E_{6}$. We have seen in \cref{sec_theta}, that both $A_{4}$ and $A_{6}$ are congruent to $1$ modulo $\mathfrak{L}$. We set
\[f' := A_{6}f \quad \text{and} \quad g' = A_{4}g.\]
Then $f'$ and $g'$ are modular forms with $\mathfrak{L}$-integral Fourier coefficients of weight $k_{f'} = k_{f}+6$, $k_{g'} = k_{g} + 4$, level $N$ and character $\varepsilon_{f}$, $\varepsilon_{g}$ respectively. Since $\bar{\chi}_{\ell}^{2}$ is trivial for $\ell = 2$, $3$, the congruence
\[\bar{\chi}_{\ell}^{k_{f'}+bm_{f}}\varepsilon_{f} \equiv \bar{\chi}_{\ell}^{k_{g'}+bm_{g}}\varepsilon_{g} \pmod{\mathfrak{L}}\]
is satisfied. Moreover, we have $k_{f'}+bm_{f} \leq k_{g'}+bm_{g}$, and $b-2 = 4$ divides $k_{g'}-k_{f'}+b(m_{g}-m_{f})$. According to the discussion at the beginning of the proof, we therefore get the desired result since $f'$, $g'$ reduce to $f$, $g$ respectively and $k_{g'}+bm_{g} = a + k_{g} + bm_{g} = k$.
\end{pf}
\begin{Rk}
Notice that \cref{Sturm} corresponds to the special case $m_{f} = m_{g} = 0$ and $k_{f} = k_{g} = k$. Moreover, in practice we can always take $m_{f} \in \{0,1\}$ and $m_{g} \in \ent{0}{\ell-1}$.
\end{Rk}
In \cref{sec_red} we will mainly deal with eigenforms. It is well-known that the knowledge of the Fourier coefficients of prime index and of the constant coefficient characterises such forms. We can therefore simplify \cref{Sturm_theta} and get the following corollary.

\begin{cor}\label{Sturm_eigen}
Let $f$, $g$ be as in \cref{Sturm_theta} and define also $N$ and $k$ similarly. Assume further that $f$ and $g$ are normalised eigenforms for the Hecke operators at level $N$ modulo $\mathfrak{L}$ of prime index less than $\frac{Nk}{12} \prod\limits_{\underset{p \text{ prime}}{p \mid N}} \left(1+\frac{1}{p}\right)$ (and different from $\ell$ if $m_{f}$, $m_{g} \geq 1$).

If $0^{m_{f}}a_{0}(f) \equiv 0^{m_{g}}a_{0}(g) \pmod{\mathfrak{L}}$ (with $0^{0} = 1$) and if for every prime number $p \leq \frac{Nk}{12} \prod\limits_{\underset{p \text{ prime}}{p \mid N}}\left(1+\frac{1}{p}\right)$ we have $p^{m_{f}}a_{p}(f) \equiv p^{m_{g}}a_{p}(g) \pmod{\mathfrak{L}}$, then we have $n^{m_{f}}a_{n}(f) \equiv n^{m_{g}}a_{n}(g) \pmod{\mathfrak{L}}$ for every non-negative integer $n$.
\end{cor}

We now give an upper-bound for the product appearing in the Sturm bound. We use a technique of Kraus \cite{Kraus95} to get a slightly better bound than the one suggested by Serre in Kraus' article.
\begin{lm}\label{Sturm_upperbound}
Let $n$ be an integer greater or equal to $2$, we have:
\[\prod_{\underset{p \text{ prime}}{p \mid n}} \left(1 + \frac{1}{p}\right) \leq 2\log\log(n) + 2.4.\]
\end{lm}

\begin{pf}
We first split the product in two parts: $\prod\limits_{\underset{p \text{ prime}}{p \mid n}} \left(1 + \frac{1}{p}\right) = P(n)Q(n)$ with
\[P(n) = \prod_{\underset{p > \log n}{p \mid n}} \left(1 + \frac{1}{p}\right)\quad \text{and}\quad Q(n) = \prod_{\underset{p \leq \log n}{p \mid n}} \left(1 + \frac{1}{p}\right).\]
Let $m$ be the number of primes $p$ dividing $n$ and being greater than $\log n$. As $n \geq \log(n)^{m}$, we get $m \leq \frac{\log n}{\log\log n}$. Thus,

\begin{equation}\label{P}
P(n) \leq \exp\left(\frac{\log n}{\log\log n} \log\left(1 + \frac{1}{\log n}\right)\right) \leq \exp\left(\frac{1}{\log\log n}\right).
\end{equation}
Applying \cite[(3.27)]{Rosser62}, we get an upper bound for $Q$:

\begin{equation}\label{Q}
Q(n) \leq \prod_{p \leq \log n} \left(1 - \frac{1}{p}\right)^{-1} < e^{\gamma} \log\log(n) \left(1 - \frac{1}{(\log \log n)^{2}}\right)^{-1},
\end{equation}
where $\gamma$ is the Euler-Mascheroni constant.

Putting (\ref{P}) and (\ref{Q}) together we have
\[\prod_{\underset{p \text{ prime}}{p \mid n}} \left(1 + \frac{1}{p}\right) \leq e^{\gamma} \log\log(n) \exp\left(\frac{1}{\log\log n}\right)\left(1 - \frac{1}{(\log \log n)^{2}}\right)^{-1}.\]
The function $x \mapsto e^{\gamma+x}(1-x^{2})^{-1}$ is bounded by $2$ for $x \in [0,0.1]$. Therefore, the lemma holds for all integers $n \geq \exp(\exp(10))$. For $n$ between $2$ and $\exp\exp(10)$, we first notice that we only have to deal with square-free integers. Then, among the square-free integers having $k$ prime factors, it suffices to only check the lemma for $n_{k} = \prod\limits_{i=1}^{k} p_{i}$, $p_{i}$ being the $i$-th prime number. The greatest $k$ such that $n_{k} \leq \exp\exp(10)$ is $2486$, and we have checked the lemma with a computer for all those $n_{k}$.
\end{pf}

\subsection{Modifying modular forms}\label{sec_raise}
In this paragraph we discuss a way to construct from a given eigenform, another eigenform with slightly different Fourier coefficients but with a bigger level. It will be crucial in \cref{sec_red}.

Let $\mathcal{O}(\mathcal{H})$ be the space of holomorphic functions on the complex upper-half plane. For an integer $n \geq 1$ and a complex number $b$, we define two operators $V_{n}$ and $S_{n}(b)$ on $\mathcal{O}(\mathcal{H})$ by
\[V_{n}: \left\{\begin{array}{rcl}
\mathcal{O}(\mathcal{H}) & \longrightarrow & \mathcal{O}(\mathcal{H})\\
h & \longmapsto & \left(z \mapsto h(nz)\right)
\end{array}\right. \quad \text{and} \quad S_{n}(b) : \left\{\begin{array}{rcl}
\mathcal{O}(\mathcal{H}) & \longrightarrow & \mathcal{O}(\mathcal{H})\\
h & \longmapsto & h - b V_{n}h
\end{array}\right..\]
For a prime number $p$, we denote by $U_{p}$ the operator which action on Fourier expansions is given by
\[U_{p}\left(\sum\limits_{n = 0}^{\infty} a_{n}q^{n}\right) := \sum\limits_{n = 0}^{\infty}  a_{np} q^{n}.\]
We recall the following facts about the operators $U_{p}$ and $V_{p}$: for any primes $p$ and $r$, the operators $V_{p}$ and $V_{r}$ commute and the image of $\Mod_{k}(M,\varepsilon)$ by $V_{p}$ is $\Mod_{k}(Mp,\varepsilon)$. Letting $V_{p}$ act on $q$-expansions, it commutes with $U_{r}$ for $r \neq p$ and satisfies $U_{p}V_{p} = \Id$. Moreover, $T_{p}^{M}$ decomposes on the space $\Mod_{k}(M,\varepsilon)$ as
\[T_{p}^{M} = U_{p} + p^{k-1}\varepsilon(p) V_{p}.\]

From now on, consider a modular form $g$ of weight $k \geq 1$, level $M \geq 1$ and character $\varepsilon$ that is a normalised eigenform for all the Hecke operators at level $M$. For any prime number $p$, we denote by $\alpha_{p}$, $\beta_{p}$ the roots of the Hecke polynomial $X^{2} - a_{p}(g)X + p^{k-1}\varepsilon(p)$.

\begin{lm}\label{Spb}
Let $p$ be a prime number and let $b \in \{\alpha_{p},\beta_{p}\}$. The function $S_{p}(b)g$ is a modular form of same weight and character as $g$ and of level $Mp^{n_{p}}$ with
\[n_{p} = \left\{\begin{array}{ll}
1 & \text{if } b \neq 0;\\
0 & \text{if } b = 0.
\end{array}\right.\]
It is a normalised eigenform for all the Hecke operators at level $Mp^{n_{p}}$, and for any prime $r$ we have
\[a_{r}\left(S_{p}(b)g\right) = \left\{\begin{array}{ll}
a_{r}(g) & \text{if } r \neq p;\\
a_{p}(g)-b & \text{if } r = p.
\end{array}\right..\]
Moreover, if $g$ has Fourier coefficients in a ring $R$, then those of $S_{p}(b)g$ lie in the ring $R(b)$.
\end{lm}
\begin{pf}
If $b = 0$, then there is nothing to prove as $S_{p}(0)g = g$. Assume $b \neq 0$. Because both $g$ and $V_{p}g$ are modular forms of weight $k$, level $Mp$ and character $\varepsilon$, it is also the case for $S_{p}(b)g$. Let us compute the action of the Hecke operators at level $Mp$ on $S_{p}(b)g$.

Let $r$ be a prime number different from $p$. The operator $T_{r}^{Mp}$ has the same action as $T_{r}^{M}$. Thus, because the operators $V_{p}$ and $T_{r}^{M}$ commute, we have
\[T_{r}^{Mp}S_{p}(b)g = T_{r}^{M}g-bV_{p}T_{r}^{M}g = a_{r}(g)g - ba_{r}(g)V_{p}g = a_{r}(g)S_{p}(b)g.\]
For $r = p$, we have $T_{p}^{Mp}g = U_{p}g = T_{p}^{M}g-p^{k-1}\varepsilon(p)V_{p}g = a_{p}(g)g - p^{k-1}\varepsilon(p)V_{p}g$. It gives
\[T_{p}^{Mp}S_{p}(b)g = \left(a_{p}(g)g - p^{k-1}\varepsilon(p)V_{p}g\right) - bU_{p}V_{p}g = \left(a_{p}(g)-b\right)g - p^{k-1}\varepsilon(p)V_{p}g.\]
As $b$ is a root of $X^{2} - a_{p}(g)X + p^{k-1}\varepsilon(p)$, it satisfies $b(a_{p}(g)-b) = p^{k-1}\varepsilon(p)$. We finally get
\[T_{p}^{Mp}S_{p}(b)g = \left(a_{p}(g)-b\right)g - \left(a_{p}(g)-b\right)bV_{p}g = \left(a_{p}(g)-b\right)S_{p}g.\]
The form $S_{p}(b)g$ is thus a normalised eigenform for the all Hecke operators at level $Mp$. The fact about the ring of Fourier coefficients of $S_{p}(b)g$ is now straightforward.
\end{pf}
We now apply this result to construct from the eigenform $g$, an eigenform which $p$-th Fourier coefficient is a chosen number $b$ in $\{\alpha_{p},\beta_{p},0\}$.

\begin{prop}\label{raise_once}
Let $p$ be a prime number and let $b \in \{\alpha_{p}, \beta_{p}, 0\}$. Define
\[\left\{\begin{array}{lll}
g_{p}^{b} = g &\text{ and } n_{p} = 
0, & \text{if } b = a_{p}(g);\\
g_{p}^{b} = S_{p}(a_{p}(g)-b)g &\text{ and } n_{p} = 1, & \text{if } b \neq a_{p}(g) \text{ and } b \in \{\alpha_{p}, \beta_{p}\};\\
g_{p}^{b} =S_{p}(\alpha_{p}) \circ S_{p}(\beta_{p})g &\text{ and }n_{p} = 2, & \text{if } b \neq a_{p}(g) \text{ and } b \notin \{\alpha_{p},\beta_{p}\}.
\end{array}\right.\]
Then, $g_{p}^{b}$ is a modular form of same weight and character as $g$ and of level $Mp^{n_{p}}$. It is a normalised eigenform for all the Hecke operators at level $Mp^{n_{p}}$, and for any prime $r$ we have
\[a_{r}\left(g_{p}^{b}\right) = \left\{\begin{array}{ll}
a_{r}(g) & \text{if } r \neq p;\\
b & \text{if } r = p.
\end{array}\right.\]
Moreover, if $g$ has Fourier coefficients in a ring $R$, then those of $g_{p}^{b}$ lie in the ring $R(b)$.
\end{prop}
\begin{pf}
In the first two cases, we have $g_{p}^{b} = S_{p}(a_{p}(g)-b)g$. \Cref{Spb} gives directly the result. In the third case, we necessarily have $b = 0$ and $\alpha_{p}$, $\beta_{p}$ non-zero. From \cref{Spb} apply to $g$ and $\beta_{p}$, the $p$-th Hecke polynomial of $S_{p}(\beta_{p})g$ is $X^{2}-\alpha_{p}X$, which $\alpha_{p}$ is a root. We can then apply \cref{Spb} to $S_{p}(\beta_{p})g$ and $\alpha_{p}$ to conclude. Finally, the calculation
\begin{align*}
S_{p}(\alpha_{p}) \circ S_{p}(\beta_{p})g &= \left(g-\beta_{p}V_{p}g\right) - \alpha_{p}V_{p}\left(g-\beta_{p}V_{p}g\right)\\
	&= g - (\alpha_{p}+\beta_{p})V_{p}g + \alpha_{p}\beta_{p}V_{p}^{2}g\\
	&= g - a_{p}(g)V_{p}g + p^{k-1}\varepsilon(p)V_{p}^{2}g,
\end{align*}
proves that the Fourier coefficients of $g_{p}^{b}$ lie in the same ring as $g$, because the values of the character of an eigenform always lie in the ring of coefficient of it (see \cite[Corollary (3.1)]{Ribet75}).
\end{pf}

\begin{Rk}\label{raise_Rk}
Notice that the modular form $g_{p}^{b}$ is always of the shape $P(V_{p})g$ with $P = 1 - \varepsilon_{p}X + \delta_{p}X^{2}$ and $(\varepsilon_{p},\delta_{p}) \in \{(\alpha_{p},0),(\beta_{p},0),(a_{p}(g),p^{k-1}\varepsilon(p))\}$.
\end{Rk}

For any prime number $p$ and $b_{p} \in \{0,\alpha_{p},\beta_{p}\}$, define
\[S_{p}^{b_{p}} = \left\{\begin{array}{ll}
\Id & \text{if } b_{p} = a_{p}(g);\\
\Id - (a_{p}(g)-b_{p})V_{p} & \text{if } b_{p} \neq a_{p}(g) \text{ and } b_{p} \in \{\alpha_{p},\beta_{p}\};\\
\Id - a_{p}(g)V_{p} + p^{k-1}\varepsilon(p)V_{p}^{2} & \text{if } b_{p} \neq a_{p}(g) \text{ and } b_{p} \notin \{\alpha_{p},\beta_{p}\},
\end{array}\right.\]
so that we have $g_{p}^{b_{p}} = S_{p}^{b_{p}}g$. By \cref{raise_once}, applying $S_{p}^{b_{p}}$ to $g$ only modifies the Fourier coefficients of index divisible by $p$. Moreover, it gives us a modular form that is still a normalised eigenform for the whole Hecke algebra at its level. It means that for another prime $r$ and $b_{r} \in \{0,\alpha_{r},\beta_{r}\}$, the modular forms $\big(g^{b_{p}}_{p}\big)_{r}^{b_{r}}$ and $\left(g_{r}^{b_{r}}\right)_{p}^{b_{p}}$ are well-defined and equal to $S_{p}^{b_{p}}S_{r}^{b_{r}}g = S_{r}^{b_{r}}S_{p}^{b_{p}}g$.

For any finite set of primes $\bm{\mathrm{P}}$ and any $\bm{\mathrm{b}} \in \prod\limits_{p \in \bm{\mathrm{P}}} \{0,\alpha_{p},\beta_{p}\}$, we define
\[g_{\bm{\mathrm{P}}}^{\bm{\mathrm{b}}} := \prod_{p \in \bm{\mathrm{P}}} S_{p}^{b_{p}}g.\]
With the notations of \cref{raise_once}, we deduce the following result.

\begin{cor}\label{raise_several}
The function $g_{\bm{\mathrm{P}}}^{\bm{\mathrm{b}}}$ is a modular form of same weight and character as $g$ and of level $M \prod\limits_{p \in \bm{\mathrm{P}}} p^{n_{p}}$. It is a normalised eigenform for all the Hecke operators at level $M \prod\limits_{p \in \bm{\mathrm{P}}} p^{n_{p}}$, and for any prime $r$ we have
\[a_{r}\left(g_{\bm{\mathrm{P}}}^{\bm{\mathrm{b}}}\right) = \left\{\begin{array}{ll}
a_{r}(g) & \text{if } r \notin \bm{\mathrm{P}};\\
b_{r} & \text{if } r \in \bm{\mathrm{P}}.
\end{array}\right.\]
Moreover, if $g$ has Fourier coefficients in a ring $R$, then those of $g_{\bm{\mathrm{P}}}^{\bm{\mathrm{b}}}$ lie in the ring $R(\bm{\mathrm{b}})$.
\end{cor}

There is another function that we can modify with the operator $S_{p}(b)$ and get a modular form: the Eisenstein series $E_{2}$. Indeed, an easy computation shows that for any prime $p$, we have
\[S_{p}(p)E_{2} = E_{2}^{\mathds{1},\mathds{1}_{(p)}}.\]
In particular, the form $S_{p}(p)E_{2} = S_{p}(a_{p}(E_{2})-1)E_{2}$ is a normalised eigenform of weight $2$, level $p$, trivial character, and for any prime $r$, its $r$-th Fourier coefficient is equal to $r+1 = a_{r}(E_{2})$ if $r \neq p$, and $1$ if $r = p$. Moreover, the Hecke polynomial at $p$ of $E_{2}^{\mathds{1},\mathds{1}_{(p)}}$ is $X(X-1)$. Thus, $S_{p}(1)E_{2}^{\mathds{1},\mathds{1}_{(p)}} = S_{p}(1) \circ S_{p}(p)E_{2}$ is a normalised eigenform of weight $2$, trivial character and level $p^{2}$ and we have $a_{p}(S_{p}(1)\circ S_{p}(p)E_{2}) = 0$. We have proved the following result.

\begin{prop}\label{raise_onceE2}
Let $p$ be any prime number and $b \in \{p,0\}$. Define
\[\left\{\begin{array}{lll}
(E_{2})^{b}_{p} = S_{p}(p)E_{2} & \text{and } n_{p} = 1 & \text{if } b = 1;\\
(E_{2})^{b}_{p} = S_{p}(1) \circ S_{p}(p)E_{2} & \text{and } n_{p} = 2 & \text{if } b = 0.
\end{array}\right.\]
The function $(E_{2})_{p}^{b}$ is a modular form of weight $2$, level $p^{n_{p}}$, and trivial character. It is a normalised eigenform for all the Hecke operators at level $p^{n_{p}}$, and for any prime $r$ we have
\[a_{r}\left((E_{2})_{p}^{b}\right) = \left\{\begin{array}{ll}
r+1 & \text{if } r \neq p;\\
b & \text{if } r = p.
\end{array}\right.\]
Moreover, all the Fourier coefficients of $\left(E_{2}\right)_{p}^{b}$ are integers, except maybe the constant one that is rational.
\end{prop}
We can then state a result of the shape of \cref{raise_several} for $E_{2}$.
\begin{cor}\label{raise_severalE2}
Let $\bm{\mathrm{P}}$ be a finite set of prime numbers and let $\bm{\mathrm{b}} \in \prod\limits_{p \in \bm{\mathrm{P}}} \{0,1,p\} \setminus (1)_{p \in \bm{\mathrm{P}}}$. There is a modular form $\left(E_{2}\right)_{\bm{\mathrm{P}}}^{\bm{\mathrm{b}}}$ of weight $2$, level $\prod\limits_{p \in \bm{\mathrm{P}}} p^{n_{p}}$, and trivial character. It is a normalised eigenform for all the Hecke operators at its level, and for any prime $r$ we have
\[a_{r}\left(\left(E_{2}\right)_{\bm{\mathrm{P}}}^{\bm{\mathrm{b}}}\right) = \left\{\begin{array}{ll}
r+1 & \text{if } r \notin \bm{\mathrm{P}};\\
b_{r} & \text{if } r \in \bm{\mathrm{P}}.
\end{array}\right.\]
Moreover, all the Fourier coefficients of $\left(E_{2}\right)_{\bm{\mathrm{P}}}^{\bm{\mathrm{b}}}$ are integers, except maybe the constant one that is rational.
\end{cor}

We finally give a result on the constant coefficient of an Eisenstein series that has been modified with \cref{raise_several}.
\begin{prop}\label{cst_coef}
Let $k \geq 2$, let $\varepsilon_{1}$, $\varepsilon_{2}$ be two primitive Dirichlet characters. Let $\bm{\mathrm{P}}$ be a finite set of prime numbers and let $\bm{\mathrm{b}} := (b_{p}) \in \prod\limits_{p \in \bm{\mathrm{P}}} \{0,\varepsilon_{1}(p),p^{k-1}\varepsilon_{2}(p)\}$, different from $(1)_{p \in \bm{\mathrm{P}}}$ if $(k,\varepsilon_{1},\varepsilon_{2}) = (2,\mathds{1},\mathds{1})$. Then the constant coefficient of $\left(E_{k}^{\varepsilon_{1},\varepsilon_{2}}\right)_{\bm{\mathrm{P}}}^{\bm{\mathrm{b}}}$ is equal to
\[\left\{\begin{array}{ll}
0 & \text{if } \varepsilon_{1} \neq \mathds{1};\\
-\frac{B_{k,\varepsilon_{2}}}{2k} \prod\limits_{p \in \bm{\mathrm{P}}} b_{p}(b_{p}-p^{k-1}\varepsilon_{2}(p)) & \text{if } \varepsilon_{1} = \mathds{1}.
\end{array}\right.\]
\end{prop}
\begin{pf}
First, if $\varepsilon_{1} \neq \mathds{1}$, then the constant coefficient of $E_{k}^{\varepsilon_{1},\varepsilon_{2}}$ is trivial by \eqref{qexp}. Assume $\varepsilon_{1} = \mathds{1}$. Then the modular form $\left(E_{k}^{\varepsilon_{1},\varepsilon_{2}}\right)_{\bm{\mathrm{P}}}^{\bm{\mathrm{b}}}$ is equal to
\[\prod\limits_{p \in \bm{\mathrm{P}}} \left(\Id - \varepsilon_{p} V_{p} + \delta_{p} V_{p}^{2}\right) E_{k}^{\varepsilon_{1},\varepsilon_{2}},\]
where
\begin{equation}\label{epsdeltap}
(\varepsilon_{p},\delta_{p}) = \left\{\begin{array}{ll}
(1+p^{k-1}\varepsilon_{2}(p),p^{k-1}\varepsilon_{2}(p)) & \text{if } b_{p} = 0;\\
(1,0) & \text{if } b_{p} = p^{k-1}\varepsilon_{2}(p);\\
(p^{k-1}\varepsilon_{2}(p),0) & \text{if } b_{p} = 1.
\end{array}\right.
\end{equation}
Therefore, the constant coefficient is equal to $-\frac{B_{k,\varepsilon_{2}}}{2k} \prod\limits_{p \in \bm{\mathrm{P}}} (1-\varepsilon_{p}+\delta_{p})$. A straightforward computation gives that $1-\varepsilon_{p}+\delta_{p}$ is equal to $0$ if $b_{p} \in \{0,p^{k-1}\varepsilon_{2}(p)\}$, and to $1-p^{k-1}\varepsilon_{2}(p)$ if $b_{p} = 1$. Therefore, if one of the $b_{p}$'s is equal to $0$ or $p^{k-1}\varepsilon_{2}(p)$, then the constant coefficient is equal to
\[0 = -\frac{B_{k,\varepsilon_{2}}}{2k} \prod\limits_{p \in \bm{\mathrm{P}}} b_{p}(b_{p}-p^{k-1}\varepsilon_{2}(p)).\]
Else, if all the $b_{p}$'s are equal to $1$, then the constant coefficient is equal to
\[-\frac{B_{k,\varepsilon_{2}}}{2k} \prod\limits_{p \in \bm{\mathrm{P}}} (1-p^{k-1}\varepsilon_{2}(p)) = -\frac{B_{k,\varepsilon_{2}}}{2k} \prod\limits_{p \in \bm{\mathrm{P}}} b_{p}(b_{p}-p^{k-1}\varepsilon_{2}(p)).\]
\end{pf}

\begin{prop}\label{Upsilon}
Let $k$, $\varepsilon_{1}$, $\varepsilon_{2}$, $\bm{\mathrm{P}}$ and $\bm{\mathrm{b}}$ be as in \cref{cst_coef}. Let $\mathfrak{c}_{1}$, $\mathfrak{c}_{2}$ be the conductors of $\varepsilon_{1}$ and $\varepsilon_{2}$ respectively. Then the constant coefficient of $\left(E_{k}^{\varepsilon_{1},\varepsilon_{2}}\right)_{\bm{\mathrm{P}}}^{\bm{\mathrm{b}}}$ at the cusp $\frac{1}{\mathfrak{c}_{2}}$ is equal to
\begin{align*}
-\varepsilon_{1}(-1) \frac{W((\varepsilon_{1}\varepsilon_{2}^{-1})_{0})}{W(\varepsilon_{2}^{-1})} \frac{B_{k,(\varepsilon_{1}^{-1}\varepsilon_{2})_{0}}}{2k} &\left(\frac{\mathfrak{c}_{2}}{\mathfrak{c}_{0}}\right)^{k} \prod_{p \mid \mathfrak{c}_{1}\mathfrak{c}_{2}} \left(1 - \frac{(\varepsilon_{1}\varepsilon_{2}^{-1})_{0}(p)}{p^{k}}\right)\\
&\times \prod\limits_{b_{p} \neq \varepsilon_{1}(p)} \left(1-\frac{\varepsilon_{1}\varepsilon_{2}^{-1}(p)}{p^{k}}\right)\prod\limits_{b_{p} \neq p^{k-1}\varepsilon_{2}(p)} \left(1-\frac{1}{p}\right).
\end{align*}
\end{prop}
\begin{pf}
Let $\gamma := \begin{pmatrix} 1 & 0\\ \mathfrak{c}_{2} & 1\end{pmatrix}$ be an element of $\SL_{2}(\Z)$ such that $\gamma\infty = \frac{1}{\mathfrak{c}_{2}}$. Write the modular form $\left(E_{k}^{\varepsilon_{1},\varepsilon_{2}}\right)_{\bm{\mathrm{P}}}^{\bm{\mathrm{b}}}$ as
\[\prod\limits_{p \in \bm{\mathrm{P}}} \left(\Id - \varepsilon_{p} V_{p} + \delta_{p} V_{p}^{2}\right) E_{k}^{\varepsilon_{1},\varepsilon_{2}},\]
with $\varepsilon_{p}$ and $\delta_{p}$ defined by \eqref{epsdeltap}. By \cref{Eisenstein_cst}, for an integer $M$, the constant coefficient of $\left(V_{M}E_{k}^{\varepsilon_{1},\varepsilon_{2}}\right)|_{k} \gamma$ is non-zero if and only if $M$ and $\mathfrak{c}_{2}$ are coprime. Under this assumption, we have, with the notations of \cref{Eisenstein_cst},
\begin{equation}\label{upsilon}
\Upsilon_{k}^{\varepsilon_{1},\varepsilon_{2}}(\gamma,M) = \left[-\varepsilon_{1}(-1) \frac{W((\varepsilon_{1}\varepsilon_{2}^{-1})_{0})}{W(\varepsilon_{2}^{-1})} \frac{B_{k,(\varepsilon_{1}^{-1}\varepsilon_{2})_{0}}}{2k} \left(\frac{\mathfrak{c}_{2}}{\mathfrak{c}_{0}}\right)^{k} \prod_{p \mid \mathfrak{c}_{1}\mathfrak{c}_{2}} \left(1 - \frac{(\varepsilon_{1}\varepsilon_{2}^{-1})_{0}(p)}{p^{k}}\right)\right] \frac{\varepsilon_{2}^{-1}(M)}{M^{k}}.
\end{equation}
The expression in brackets is independent from $M$, let write it $\mathbf{D}$. Notice that if $M$ is not coprime to $\mathfrak{c}_{2}$, the formula still holds, as $\varepsilon_{2}(M) = 0$. Define
\[P := \prod_{p \in \bm{\mathrm{P}}} \left(1 - \varepsilon_{p}X_{p} + \delta_{p}X_{p}^{2}\right) \in \C\left[(X_{p})_{p \in \bm{\mathrm{P}}}\right].\]
As \eqref{upsilon} is fully multiplicative in $M$, the constant coefficient of $\left(E_{k}^{\varepsilon_{1},\varepsilon_{2}}\right)_{\bm{\mathrm{P}}}^{\bm{\mathrm{b}}}$ is then equal to
\[\lim\limits_{\Im(z) \to +\infty} P((V_{p})_{p \in \bm{\mathrm{P}}})E_{k}^{\varepsilon_{1},\varepsilon_{2}}|_{k}\gamma(z) = \mathbf{D} \cdot P\left(\left(\frac{\varepsilon_{2}^{-1}(p)}{p^{k}}\right)_{p \in \bm{\mathrm{P}}}\right).\]
We just have to compute the value of $P\left(\left(\frac{\varepsilon_{2}^{-1}(p)}{p^{k}}\right)_{p \in \bm{\mathrm{P}}}\right)$ to conclude. Let $P_{p}(X_{p}) = 1-\varepsilon_{p}X_{p}+\delta_{p}X_{p}^{2}$, so that we have $P = \prod\limits_{p \in \mathbf{P}} P_{p}(X_{p})$. A straightforward calculation gives us the value of $P_{p}\left(\frac{\varepsilon_{2}^{-1}(p)}{p^{k}}\right)$:
\[P_{p}\left(\frac{\varepsilon_{2}^{-1}(p)}{p^{k}}\right) = \left\{\begin{array}{ll}
\left(1 - \frac{\varepsilon_{1}\varepsilon_{2}^{-1}(p)}{p^{k}}\right)\left(1 - \frac{1}{p}\right) & \text{if } (\varepsilon_{p},\delta_{p}) = (\varepsilon_{1}(p)+p^{k-1}\varepsilon_{2}(p),p^{k-1}\varepsilon_{1}\varepsilon_{2}(p));\\
1 - \frac{\varepsilon_{1}\varepsilon_{2}^{-1}(p)}{p^{k}} & \text{if } (\varepsilon_{p},\delta_{p}) = (\varepsilon_{1}(p),0);\\
1-\frac{1}{p} & \text{if } (\varepsilon_{p},\delta_{p}) = (p^{k-1}\varepsilon_{2}(p),0).
\end{array}\right.\]
\end{pf}

\section{Reducible Galois representations}\label{sec_red}

\subsection{Reducible Galois representations and Eisenstein series}\label{sec_redEisenstein}
Before dealing with the reducibility of a residual representation attached to a newform, we examine the general case of a reducible, residual, semi-simple, odd Galois representation. We begin by a definition that will be enlightened in the following proposition. Throughout this paragraph, we fix a prime number $\ell$ and a place $\mathfrak{L}$ of $\bar{\Q}$ above $\ell$.

\begin{df}\label{Rl}
Define $R(\ell)$ as the set of quadruplets $(\varepsilon_{1},\varepsilon_{2},m_{1},m_{2})$ consisting of two primitive Dirichlet characters $\varepsilon_{1}$ and $\varepsilon_{2}$ that are of prime-to-$\ell$ order and unramified at $\ell$, and of two integers $m_{1}$, $m_{2}$ such that $0 \leq m_{1} \leq m_{2} < \ell - 1$, and $\varepsilon_{1}\varepsilon_{2}(-1) \equiv (-1)^{m_{1}+m_{2}+1} \pmod{\ell}$.
\end{df}

\begin{prop}\label{red_params0}
Let $\bar{\rho} : G_{\Q} \to \GL_{2}(\bar{\F}_{\ell})$ be a semi-simple, odd representation. Then, $\bar{\rho}$ is reducible if and only if there exists $(\varepsilon_{1},\varepsilon_{2},m_{1},m_{2}) \in R(\ell)$ such that $\bar{\rho} \cong \bar{\chi}_{\ell}^{m_{1}}\bar{\varepsilon_{1}} \oplus \bar{\chi}_{\ell}^{m_{2}} \bar{\varepsilon_{2}}$, where $\bar{\varepsilon_{1}}$ and $\bar{\varepsilon_{2}}$ correspond to the reduction modulo $\mathfrak{L}$ of $\varepsilon_{1}$ and $\varepsilon_{2}$ respectively.
\end{prop}

\begin{pf}
If $\bar{\rho}$ is reducible, then we can decompose it as $\varphi_{1}\bar{\chi}_{\ell}^{m_{1}} \oplus \varphi_{2}\bar{\chi}_{\ell}^{m_{2}}$, with $\varphi_{1}$, $\varphi_{2}$ two characters modulo $\mathfrak{L}$ unramified at $\ell$, and $0 \leq m_{1} \leq m_{2} < \ell-1$. Let $\varepsilon_{i}$ be the Teichmüller lift of $\varphi_{i}$ with respect to $\mathfrak{L}$ (see \cref{sec_Tlifts}). By construction, the characters $\varepsilon_{1}$, $\varepsilon_{2}$ are primitive and have prime-to-$\ell$ order. Finally, if $c$ is a complex conjugation, then we have
\[-1 = \det(\bar{\rho}(c)) = \bar{\varepsilon_{1}\varepsilon_{2}}(-1)\bar{\chi}_{\ell}(c)^{m_{1}+m_{2}} = (-1)^{m_{1}+m_{2}}\bar{\varepsilon_{1}\varepsilon_{2}}(-1).\]
Therefore we have $(\varepsilon_{1},\varepsilon_{2},m_{1},m_{2}) \in R(\ell)$.
\end{pf}

Our goal is now to construct from an element $(\varepsilon_{1},\varepsilon_{2},m_{1},m_{2})$ of $R(\ell)$, various modular forms which we can attach a Galois representation modulo $\mathfrak{L}$ that is isomorphic to $\bar{\chi}_{\ell}^{m_{1}}\bar{\varepsilon_{1}} \oplus \bar{\chi}_{\ell}^{m_{2}}\bar{\varepsilon_{2}}$. According to \cite[Théorème 6.7]{Deligne74}, for a modular form $g \in \Mod_{k_{g}}(N_{g},\varepsilon_{g})$ to admit a Galois representation modulo $\mathfrak{L}$, it suffices that
\begin{itemize}
\item $g$ has $\mathfrak{L}$-integral Fourier coefficients;
\item $g$ is a normalised eigenform for the Hecke operators $T_{p}^{N_{g}}$ modulo $\mathfrak{L}$ of prime index $p$ not dividing $\ell N_{g}$.
\end{itemize}
Moreover, this representation $\bar{\rho}_{g,\mathfrak{L}}$ is characterised up to isomorphism by the properties: $\det(\bar{\rho}_{g,\mathfrak{L}}) = \bar{\chi}_{\ell}^{k_{g}-1} \bar{\varepsilon_{g}}$, $\bar{\rho}_{g,\mathfrak{L}}$ is unramified outside $N_{g}\ell$, and for any prime $p \nmid N_{g}\ell$, we have $\Tr(\bar{\rho}_{g,\mathfrak{L}}(\mathrm{Frob}_{p})) = a_{p}(g) \pmod{\mathfrak{L}}$. In particular by \cite[Lemme 3.2.]{Deligne74}, to have $\bar{\rho}_{g,\mathfrak{L}} \cong \bar{\chi}_{\ell}^{m_{1}}\bar{\varepsilon_{1}} \oplus \bar{\chi}_{\ell}^{m_{2}}\bar{\varepsilon_{2}}$, it suffices to have $\bar{\chi}_{\ell}^{k_{g}-1}\bar{\varepsilon_{g}} = \bar{\chi}_{\ell}^{m_{1}+m_{2}}\bar{\varepsilon_{1}}\bar{\varepsilon_{2}}$ and for any prime $p \nmid N_{g}\mathfrak{c}_{1}\mathfrak{c}_{2}\ell$ in a set of density one, $a_{p}(g) \equiv p^{m_{1}}\varepsilon_{1}(p) + p^{m_{2}}\varepsilon_{2}(p) \pmod{\mathfrak{L}}$, where $\mathfrak{c}_{i}$ is the conductor of $\varepsilon_{i}$.

For $(\varepsilon_{1},\varepsilon_{2},m_{1},m_{2})$ in $R(\ell)$, define
\begin{equation}\label{kE0}
k' := \left\{\begin{array}{ll}
m_{2}-m_{1}+1 & \text{if } \ell > 2;\\
2 & \text{if } \ell = 2,
\end{array}\right. \quad \text{and} \quad E_{0} := E_{k'}^{\varepsilon_{1},\varepsilon_{2}}.
\end{equation}

\begin{prop}\label{E0}
We have $\varepsilon_{1}\varepsilon_{2}(-1) = (-1)^{k'}$. In particular, $E_{0}$ is well-defined and modular if and only if $(k',\varepsilon_{1},\varepsilon_{2}) \neq (2,\mathds{1},\mathds{1})$, in which case $E_{0}$ is a normalised eigenform of weight $k'$, level $\mathfrak{c}_{1}\mathfrak{c}_{2}$, and character $\varepsilon_{1}\varepsilon_{2}$. Moreover, for any prime number $p$ we have $a_{p}(E_{0}) = \varepsilon_{1}(p) + p^{k'-1}\varepsilon_{2}(p)$ in any case.
\end{prop}
\begin{pf}
If $\ell = 2$, then $\varepsilon_{1}$ and $\varepsilon_{2}$ are even and we have $\varepsilon_{1}\varepsilon_{2}(-1) = 1 = (-1)^{k'}$. Otherwise, we have
\[\varepsilon_{1}\varepsilon_{2}(-1) = (-1)^{m_{2}+m_{1}+1} = (-1)^{m_{2}-m_{1}+1} = (-1)^{k'}.\]
The rest of the proposition follows from \cref{Eisenstein} and \eqref{qexp}.
\end{pf}
To be sure that we can associate a Galois representation modulo $\mathfrak{L}$ to $E_{0}$, we study the integrality of its Fourier coefficients. The following lemma states when the coefficients of $E_{0}$ may not be $\mathfrak{L}$-integral.

\begin{lm}\label{integral}
Assume $(k',\varepsilon_{1},\varepsilon_{2}) \neq (2,\mathds{1},\mathds{1})$. The Fourier coefficients of $E_{0}$ are $\mathfrak{L}$-integral unless perhaps in the following cases:
\begin{itemize}
\item $\ell = 2$, $\varepsilon_{1} = \mathds{1}$ and $\varepsilon_{2} \neq \mathds{1}$;
\item $\ell \geq 5$, $\varepsilon_{1} = \varepsilon_{2} = \mathds{1}$, and $(m_{1},m_{2}) = (0,\ell-2)$.
\end{itemize}
\end{lm}
\begin{pf}
Apart from the constant one, the coefficients of $E_{0}$ are all algebraic integers. We therefore only need to focus on the constant Fourier coefficient $a_{0}$ of $E_{0}$.

In the case $\ell \neq 2$, if $(\varepsilon_{1},\varepsilon_{2}) \neq (\mathds{1},\mathds{1})$, then $a_{0}$ is always $\mathfrak{L}$-integral by \cref{Ber_divisors}, because $\varepsilon_{1}$ and $\varepsilon_{2}$ are unramified at $\ell$. If $\varepsilon_{1} = \varepsilon_{2} = \mathds{1}$, then $a_{0} = -\frac{1}{2k'}B_{k'}$. By \cref{Ber_divisors} again, if $(m_{1},m_{2}) \neq (0,\ell-2)$, then $a_{0}$ is always $\mathfrak{L}$-integral. If $(m_{1},m_{2}) = (0,\ell-2)$, then $a_{0}$ is always not $\mathfrak{L}$-integral. Notice moreover that we must have $\ell \neq 3$, because otherwise $k' = 2$ and $\varepsilon_{1} = \varepsilon_{2} = \mathds{1}$, which is excluded.

Assume $\ell = 2$, hence $k' = 2$. If $\varepsilon_{1} \neq \mathds{1}$, then as before we have $a_{0} = 0$. Else, if $\varepsilon_{1} = \mathds{1}$, then $a_{0} = -\frac{B_{2,\varepsilon_{2}}}{4}$ which may not be $\mathfrak{L}$-integral. Moreover, we must have $\varepsilon_{2} \neq \mathds{1}$ because otherwise we would have $(k',\varepsilon_{1},\varepsilon_{2}) = (2,\mathds{1},\mathds{1})$.
\end{pf}
We can now construct from $(\varepsilon_{1},\varepsilon_{2},m_{1},m_{2}) \in R_{\ell}$ (and hence $k_{0}$ and $E_{0}$), a modular form with the properties discussed above. With the notations of \cref{raise_once} and \cref{raise_onceE2}, define
\begin{equation}\label{rE}
\left\{\begin{array}{lll}
r := 4 & \text{and } E := \left(E_{0}\right)_{2}^{0} & \text{if }\left|\begin{array}{ll}
    & \text{we are in one of the cases listed in \cref{integral}}\\
    \text{or} & (k',\varepsilon_{1},\varepsilon_{2}) = (2,\mathds{1},\mathds{1});
    \end{array}\right.\\
r := 1 & \text{and } E := E_{0} & \text{otherwise}.
\end{array}\right.
\end{equation}
The following proposition sums up the properties of $E$.
\begin{prop}\label{E}
The function $E$ is a modular form of weight $k'$, level $M := \lcm(\mathfrak{c}_{1}\mathfrak{c}_{2},r)$, and character $\varepsilon_{1}\varepsilon_{2}$. It is a normalised eigenform for all the Hecke operators at level $M$, all its Fourier coefficients are $\mathfrak{L}$-integral, and for any prime $p \nmid r$, we have $a_{p}(E) = \varepsilon_{1}(p)+p^{k'-1}\varepsilon_{2}(p)$.
\end{prop}
\begin{pf}
The only thing to prove it that $M$ is indeed the level of $E$. The rest of the proposition then follows from \cref{raise_once}, \cref{raise_onceE2}, \cref{integral} and \cref{E0}. If $r = 1$, the level of $E$ is equal to $\mathfrak{c}_{1}\mathfrak{c}_{2} = \lcm(\mathfrak{c}_{1}\mathfrak{c}_{2},r)$. Assume $r = 4$. We then always have $\mathfrak{c}_{1} = 1$ and either $\varepsilon_{2} = \mathds{1}$ or $\ell = 2$. In the first case, $\mathfrak{c}_{2} = 1$ and $\varepsilon_{2}(2) = 1 \neq 0$. In the second case, $\mathfrak{c}_{2}$ is odd because prime to $\ell$. Thus, we have $\varepsilon_{2}(2) \neq 0$. In every case, the level of $E$ is equal to $4\mathfrak{c}_{1}\mathfrak{c}_{2} = \lcm(\mathfrak{c}_{1}\mathfrak{c}_{2},4)$.
\end{pf}
The final step to have a modular form that has a Galois representation modulo $\mathfrak{L}$ that is isomorphic to $\bar{\chi}_{\ell}^{m_{1}}\bar{\varepsilon_{1}} \oplus \bar{\chi}_{\ell}^{m_{2}}\bar{\varepsilon_{2}}$ is to apply the operator $\tilde{\theta}$ that we have constructed is \cref{sec_theta}. 

\begin{prop}\label{red_forms}
Let $\bar{\rho} : G_{\Q} \to \GL_{2}(\bar{\F}_{\ell})$ be a residual, semi-simple, odd Galois representation. Then $\bar{\rho}$ is reducible if and only if there exist $(\varepsilon_{1},\varepsilon_{2},m_{1},m_{2}) \in R(\ell)$ such that $\bar{\rho} \cong \bar{\rho}_{\tilde{\theta}^{m_{1}}E,\mathfrak{L}}$, where $E$ is the Eisenstein series associated to $(\varepsilon_{1},\varepsilon_{2},m_{1},m_{2})$ in \eqref{rE}.

Moreover, if $\bm{\mathrm{P}}$ is a finite set of prime numbers, and $\bm{\mathrm{b}} \in \prod\limits_{p \in \bm{\mathrm{P}}} \{0,\varepsilon_{1}(p),p^{k'-1}\varepsilon_{2}(p)\}$, then we also have $\bar{\rho} \cong \bar{\rho}_{\tilde{\theta}^{m_{1}}E_{\bm{\mathrm{P}}}^{\bm{\mathrm{b}}},\mathfrak{L}}$, where $E_{\bm{\mathrm{P}}}^{\bm{\mathrm{b}}}$ is defined in \cref{raise_several}.
\end{prop}
\begin{pf}
From \cref{red_params0}, we only have to prove that for a quadruplet $(\varepsilon_{1},\varepsilon_{2},m_{1},m_{2}) \in R(\ell)$, we have $\bar{\rho}_{\tilde{\theta}^{m_{1}}E,\mathfrak{L}} \cong \bar{\chi}_{\ell}^{m_{1}}\bar{\varepsilon_{1}} \oplus \bar{\chi}_{\ell}^{m_{2}}\bar{\varepsilon_{2}}$. It follows from \cref{E} and \cref{theta_lm} that the representation $\bar{\rho}_{\tilde{\theta}^{m_{1}}E,\mathfrak{L}}$ is well-defined. Moreover, from \cref{theta_lm} again, the form $\tilde{\theta}_{A}^{m_{1}}E$ is of weight $k'+m_{1}(k_{A}+2)$, level $M$, and character $\varepsilon_{1}\varepsilon_{2}\chi_{A}^{m_{1}}$, the value of $k_{A}$ and $\chi_{A}$ depending on the level $M$ and the place $\mathfrak{L}$ (see \cref{theta_table}). Recall that in any case we have $\chi_{A} \equiv \bar{\chi}_{\ell}^{-k_{A}} \pmod{\mathfrak{L}}$ and $k' \equiv m_{2}-m_{1}+1 \pmod{\ell-1}$. Therefore, we have \[\bar{\chi}_{\ell}^{(k'+m_{1}(k_{A}+2))-1}\chi_{A}^{m_{1}} \varepsilon_{1}\varepsilon_{2} \equiv \bar{\chi}_{\ell}^{m_{2}-m_{1}+1+m_{1}(k_{A}+2)-1-m_{1}k_{A}} \varepsilon_{1}\varepsilon_{2} \equiv \bar{\chi}_{\ell}^{m_{1}+m_{2}}\varepsilon_{1}\varepsilon_{2} \pmod{\mathfrak{L}}.\]
Finally, from \cref{theta_lm}, \cref{E}, and again the congruence $k' \equiv m_{2}-m_{1}+1 \pmod{\ell-1}$, we have for any prime number $p \nmid r\ell$, 
\[a_{p}\left(\tilde{\theta}^{m_{1}}E\right) \equiv p^{m_{1}}a_{p}(E) \equiv p^{m_{1}}\varepsilon_{1}(p) + p^{k'+m_{1}-1}\varepsilon_{2}(p) \equiv p^{m_{1}}\varepsilon_{1}(p) + p^{m_{2}}\varepsilon_{2}(p) \pmod{\mathfrak{L}}.\]
By the discussion at the beginning of the paragraph, we get $\bar{\rho}_{\tilde{\theta}^{m_{1}}E,\mathfrak{L}} \cong \bar{\chi}_{\ell}^{m_{1}}\bar{\varepsilon_{1}} \oplus \bar{\chi}_{\ell}^{m_{2}}\bar{\varepsilon_{2}}$.

By \cref{raise_several} and \cref{theta_lm}, the modular form $\tilde{\theta}^{m_{1}}E_{\mathbf{P}}^{\mathbf{b}}$ still has $\mathfrak{L}$-integral Fourier coefficients, and is again an eigenform for all the Hecke operators modulo $\mathfrak{L}$ at its level. Moreover, by the same computations as above, the determinant of $\bar{\rho}_{\tilde{\theta}^{m_{1}}E_{\mathbf{P}}^{\mathbf{b}}}$ is still equal to $\bar{\chi}_{\ell}^{m_{1}+m_{2}}\bar{\varepsilon_{1}}\bar{\varepsilon_{2}}$, and for any prime number $p \nmid r\ell$, $p \notin \mathbf{P}$, we have $a_{p}\left(\tilde{\theta}^{m_{1}}E_{\mathbf{P}}^{\mathbf{b}}\right) \equiv p^{m_{1}}a_{p}\left(E_{\mathbf{P}}^{\mathbf{b}}\right) \equiv p^{m_{1}}\varepsilon_{1}(p)+p^{m_{2}}\varepsilon_{2}(p) \pmod{\mathfrak{L}}$. Therefore, we also have $\bar{\rho} \cong \bar{\rho}_{\tilde{\theta}^{m_{1}}E_{\mathbf{P}}^{\mathbf{b}}}$.
\end{pf}
\begin{Rk}
Be careful that the operator $\tilde{\theta}$ we apply to $E_{\bf{P}}^{\bf{b}}$ is not strictly the same as the one we apply to $E$, the levels of these two forms not being the same.
\end{Rk}

\subsection{General study of modular reducible representations}\label{sec_rednew}
Let $f = q + \sum\limits_{n = 2}^{\infty} a_{n}(f)q^{n}$ be a newform of weight $k \geq 2$, level $N \geq 1$, and character $\varepsilon$ of conductor $\mathfrak{c}$. Let $K_{f}$ be the number field generated by $(a_{n}(f))_{n \geq 2}$ and let $\lambda$ be a prime ideal of the ring of integers of $K_{f}$ above a prime number $\ell$. As in the previous paragraph, we begin with the definition of a set that corresponds to the possible reductions of $\bar{\rho}_{f,\lambda}$.

\begin{df}\label{RNkeps}
Let $\mathfrak{L}$ be a place of $\bar{\Q}$ above $\ell$. Define the set $R_{N,k,\varepsilon}(\mathfrak{L})$ as the subset of $R(\ell)$ (see \cref{Rl}) consisting of the quadruples $(\varepsilon_{1},\varepsilon_{2},m_{1},m_{2})$ such that $\bar{\chi}_{\ell}^{m_{1}+m_{2}} \bar{\varepsilon_{1}}\bar{\varepsilon_{2}} = \bar{\chi}_{\ell}^{k-1} \bar{\varepsilon}$, and for every prime $p \neq \ell$, we have $v_{p}\left(\frac{N}{\mathfrak{c}_{1}\mathfrak{c}_{2}}\right) \in \{0,1,2\}$, where $\bar{\cdot}$ denotes the reduction modulo $\mathfrak{L}$, and $\mathfrak{c}_{i}$ is the conductor of $\varepsilon_{i}$.

Notice that in particular if $(\varepsilon_{1},\varepsilon_{2},m_{1},m_{2}) \in R_{N,k,\varepsilon}(\mathfrak{L})$, then $\mathfrak{c}_{1}\mathfrak{c}_{2} \mid N$. The set $R_{N,k,\varepsilon}(\mathfrak{L})$ is therefore finite.
\end{df}
\begin{Rk}
We will see later that the set $R_{N,k,\varepsilon}(\mathfrak{L})$ depends in fact only on $\mathfrak{L}\cap\Q(\varepsilon)$ (and obviously on $N$, $k$ and $\varepsilon$). For now, this dependency will not matter and we postpone this proof to \cref{sec_Num}.
\end{Rk}
From now on, assume that $\mathfrak{L}$ is a place of $\bar{\Q}$ extending $\lambda$.
\begin{prop}\label{red_params}
The representation $\bar{\rho}_{f,\lambda}$ is reducible if and only if there exists $(\varepsilon_{1},\varepsilon_{2},m_{1},m_{2}) \in R_{N,k,\varepsilon}(\mathfrak{L})$ such that $\bar{\rho}_{f,\lambda} \cong \bar{\chi}_{\ell}^{m_{1}}\bar{\varepsilon_{1}} \oplus \bar{\chi}_{\ell}^{m_{2}}\bar{\varepsilon_{2}}$, where $\bar{\varepsilon_{i}}$ denotes the reduction of $\varepsilon_{i}$ modulo $\mathfrak{L}$.
\end{prop}
\begin{pf}
The representation $\bar{\rho}_{f,\lambda}$ is semi-simple and odd. Therefore, by \cref{red_params0}, it is reducible if and only if there exist a place $\mathfrak{L}$ above $\ell$ and $(\varepsilon_{1},\varepsilon_{2},m_{1},m_{2}) \in R(\ell)$ such that $\bar{\rho}_{f,\lambda} \cong \bar{\chi}_{\ell}^{m_{1}}\bar{\varepsilon_{1}} \oplus \bar{\chi}_{\ell}^{m_{2}}\bar{\varepsilon_{2}}$. We just have to check that in this case we have $(\varepsilon_{1},\varepsilon_{2},m_{1},m_{2}) \in R_{N,k,\varepsilon}(\mathfrak{L})$. The determinant of $\bar{\rho}_{f,\lambda}$ is $\bar{\chi}_{\ell}^{k-1}\bar{\varepsilon}$. Therefore, we have $\bar{\chi}_{\ell}^{m_{1}+m_{2}} \bar{\varepsilon_{1}}\bar{\varepsilon_{2}} = \bar{\chi}_{\ell}^{k-1} \bar{\varepsilon}$. Moreover, the prime-to-$\ell$ part of the Artin conductor of $\bar{\chi}_{\ell}^{m_{1}}\bar{\varepsilon_{1}} \oplus \bar{\chi}_{\ell}^{m_{2}}\bar{\varepsilon_{2}}$ is equal to $\mathfrak{c}_{1}\mathfrak{c}_{2}$. By \cref{CL}, we necessarily have $v_{p}\left(\frac{N}{\mathfrak{c}_{1}\mathfrak{c}_{2}}\right) \in \{0,1,2\}$ for all primes $p \nmid N$, $p \neq \ell$.
\end{pf}

In regard of \cref{red_forms}, this result is equivalent to saying that the prime index coefficients of $f$ are all but finitely many congruent modulo $\mathfrak{L}$ to those of any of the forms $\tilde{\theta}^{m_{1}}E_{\mathbf{P}}^{\mathbf{b}}$ described in paragraph \ref{sec_redEisenstein}. We will prove that there exist in fact $\mathbf{P}$ and $\mathbf{b}$ such that $f$ is congruent to $\tilde{\theta}^{m_{1}}E_{\mathbf{P}}^{\mathbf{b}}$ except maybe at the primes dividing $\ell$ and $r$. The following result is the key step in this direction. It uses in a crucial way the local description of $\bar{\rho}_{f,\lambda}$ at the bad prime numbers (see \cref{sec_BackGal}).

\begin{lm}\label{bp}
If the representation $\bar{\rho}_{f,\lambda}$ is reducible, then there exists $(\varepsilon_{1},\varepsilon_{2},m_{1},m_{2}) \in R_{N,k,\varepsilon}(\mathfrak{L})$ such that for any prime number $p \neq \ell$, we have
\begin{equation}\label{red_cond}
a_{p}(f) \equiv \left\{\begin{array}{ll}
p^{m_{1}}\varepsilon_{1}(p) + p^{m_{2}}\varepsilon_{2}(p) \pmod{\mathfrak{L}} & \text{if } p \nmid N;\\
p^{m_{1}}b_{p} \pmod{\mathfrak{L}} & \text{if } p \mid N \text{, for some } b_{p} \in \{0,\varepsilon_{1}(p),p^{m_{2}-m_{1}}\varepsilon_{2}(p)\}.
\end{array}\right.
\end{equation}
Conversely, if for some $(\varepsilon_{1},\varepsilon_{2},m_{1},m_{2}) \in R_{N,k,\varepsilon}(\mathfrak{L})$, those congruences hold for every prime $p$ in a set of density $1$, then we have $\bar{\rho}_{f,\lambda} \cong \bar{\chi}_{\ell}^{m_{1}}\bar{\varepsilon_{1}} \oplus \bar{\chi}_{\ell}^{m_{2}}\bar{\varepsilon_{2}}$.
\end{lm}
\begin{pf}
Let us first prove the second statement. Write $\bar{\rho} := \bar{\chi}_{\ell}^{m_{1}}\bar{\varepsilon_{1}}\oplus\bar{\chi}_{\ell}^{m_{2}}\bar{\varepsilon_{2}}$. By construction, the determinants of $\bar{\rho}_{f,\lambda}$ and $\bar{\rho}$ agree. Moreover, by assumption for any prime number $p \nmid N\ell$ in a set of density $1$, we have
\[\Tr\left(\bar{\rho}_{f,\lambda}(\mathrm{Frob}_{p})\right) \equiv a_{p}(f) \equiv p^{m_{1}}\varepsilon_{1}(p)+p^{m_{2}}\varepsilon_{2}(p) \equiv \Tr\left(\bar{\rho}(\mathrm{Frob}_{p})\right) \pmod{\mathfrak{L}}.\]
By Brauer-Nesbitt theorem \cite[Lemma 3.2.]{Deligne74}, $\bar{\rho}_{f,\lambda}$ must be isomorphic to $\bar{\rho}$ and is thus reducible.

We now prove the first statement. Assume $\bar{\rho}_{f,\lambda}$ to be reducible. \Cref{red_params} gives us the existence of $(\varepsilon_{1},\varepsilon_{2},m_{1},m_{2}) \in R_{N,k,\varepsilon}(\mathfrak{L})$, such that $\bar{\rho}_{f,\lambda} \cong \bar{\chi}_{\ell}^{m_{1}}\bar{\varepsilon_{1}}\oplus\bar{\chi}_{\ell}^{m_{2}}\bar{\varepsilon_{2}}$. For any prime $p \nmid N\ell$, taking the trace at a Frobenuis at $p$ gives the congruence $a_{p}(f) \equiv p^{m_{1}}\varepsilon_{1}(p)+p^{m_{2}}\varepsilon_{2}(p) \pmod{\mathfrak{L}}$. Let us now consider a prime $p \mid N$ and different from $\ell$. We treat $3$ cases separately:
\begin{enumerate}[(i)]
\item If $v_{p}(N) \geq 2$ and $v_{p}(N) > v_{p}(\mathfrak{c})$, we know that $a_{p}(f) = 0$ (see \cite[Theorem 4.6.17]{Miyake06}). Hence, we have $a_{p}(f) \equiv p^{m_{1}}b_{p} \pmod{\mathfrak{L}}$ with $b_{p} = 0$.

\item If $v_{p}(N) = 1$ and $v_{p}(\mathfrak{c}) = 0$, then $\varepsilon$ is unramified at $p$. Moreover, by the congruence $\bar{\chi}_{\ell}^{k-1}\varepsilon \equiv \bar{\chi}_{\ell}^{m_{1}+m_{2}}\varepsilon_{1}\varepsilon_{2} \pmod{\mathfrak{L}}$, the character $\varepsilon_{1}\varepsilon_{2}$ is then also unramified at $p$. Therefore, $p \mid \mathfrak{c}_{1}$ if and only if $p \mid \mathfrak{c}_{2}$, and because $v_{p}(N) = 1$ and $\mathfrak{c}_{1}\mathfrak{c}_{2} \mid N$, we deduce that $p$ does not divide $\mathfrak{c}_{1}\mathfrak{c}_{2}$. Thus, we are in the first case of \cref{LW} and we get an equality of sets of characters of $G_{p}$:
\[\left\{\mu\left(a_{p}(f)\right),\mu\left(a_{p}(f)\right) \bar{\chi}_{\ell}\right\} = \left\{\bar{\chi}_{\ell}^{m_{1}}\bar{\varepsilon_{1}},\bar{\chi}_{\ell}^{m_{2}}\bar{\varepsilon_{2}}\right\}.\]
There are two cases to look at:
\begin{itemize}
\item If $\mu\left(a_{p}(f)\right) = \bar{\chi}_{\ell}^{m_{1}}\bar{\varepsilon_{1}}$, then $a_{p}(f) \equiv \varepsilon_{1}(p) p^{m_{1}} \pmod{\mathfrak{L}}$. In this case, we define $b_{p} = \varepsilon_{1}(p)$.
\item If $\mu\left(a_{p}(f)\right) = \bar{\chi}_{\ell}^{m_{2}}\bar{\varepsilon_{2}}$, then $a_{p}(f) \equiv \varepsilon_{2}(p) p^{m_{2}} \pmod{\mathfrak{L}}$. We put $b_{p} = p^{m_{2}-m_{1}}\varepsilon_{2}(p)$.
\end{itemize}
In both cases we have $a_{p}(f) \equiv p^{m_{1}}b_{p} \pmod{\mathfrak{L}}$ with $b_{p} \in \{\varepsilon_{1}(p),p^{m_{2}-m_{1}}\varepsilon_{2}(p)\}$.

\item Finally, if $v_{p}(N) = v_{p}(\mathfrak{c})$, we are in the second case of \cref{LW} and we get the equality 
\[\left\{\mu\left(a_{p}(f)\right),\mu\left(a_{p}(f)^{-1}\right) \bar{\varepsilon}_{|G_{p}} \bar{\chi}_{\ell}^{k-1}\right\} = \left\{\bar{\chi}_{\ell}^{m_{1}}\bar{\varepsilon_{1}}, \bar{\chi}_{\ell}^{m_{2}}\bar{\varepsilon_{2}}\right\},\]
We again have two cases to consider:
\begin{itemize}
\item If $\mu\left(a_{p}(f)\right) = \bar{\chi}_{\ell}^{m_{1}}\bar{\varepsilon_{1}}$, then $a_{p}(f) \equiv \varepsilon_{1}(p) p^{m_{1}} \pmod{\mathfrak{L}}$. Let us put $b_{p} = \varepsilon_{1}(p)$.
\item If $\mu\left(a_{p}(f)\right) = \bar{\chi}_{\ell}^{m_{2}}\bar{\varepsilon_{1}}$, then $a_{p}(f) \equiv \varepsilon_{2}(p) p^{m_{2}} \pmod{\mathfrak{L}}$. We put $b_{p} = p^{m_{2}-m_{1}} \varepsilon_{2}(p)$.
\end{itemize}
In both cases, we again have $a_{p}(f) \equiv p^{m_{1}}b_{p} \pmod{\mathfrak{L}}$ with $b_{p} \in \{\varepsilon_{1}(p),p^{m_{2}-m_{1}}\varepsilon_{2}(p)\}$.
\end{enumerate}
\end{pf}

For $(\varepsilon_{1},\varepsilon_{2},m_{1},m_{2}) \in R_{N,k,\varepsilon}(\mathfrak{L})$ and consider $k'$, $r$ and $E$ as defined in \eqref{kE0} and \eqref{rE} respectively. With the notations of \cref{raise_once}, we define
\begin{equation}\label{Nf'}
\left\{\begin{array}{lll}
f' := f_{2}^{0} & \text{and } N' := \left\{\begin{array}{ll}
	N & \text{if } 2 \mid N \text{ and } a_{2}(f) = 0;\\
	2N & \text{if } 2 \mid N \text{ and } a_{2}(f) \neq 0;\\
	4N & \text{if } 2 \nmid N,
\end{array}\right. & \text{if } r = 4;\\
f' := f & \text{and } N' := N, & \text{if } r = 1.
\end{array}\right.
\end{equation}

\begin{prop}\label{f'}
The form $f'$ is a normalised eigenform of weight $k$, level $N'$, and character $\varepsilon$. Its Fourier coefficients are $\mathfrak{L}$-integral and if a prime $p$ divides $r$, then $a_{p}(f') = a_{p}(E) = 0$. Moreover, the level $M := \lcm(\mathfrak{c}_{1}\mathfrak{c}_{2},r)$ of $E$ (see \cref{E}) always divides $N'$, and if $\ell = 2$, then $N' \geq 3$.
\end{prop}
\begin{pf}
The only facts that do not follow directly from \cref{raise_once} are those on the level $N'$. First, if $r = 1$, then by construction $N' = N$ and $M = \mathfrak{c}_{1}\mathfrak{c}_{2} \mid N$. Assume $r = 4$. If $2 \nmid N$, then $N' = 4N$ is divisible by $r$. If $2 \mid N$ and $a_{2}(f) \neq 0$, then $N' = 2N$ is also divisible by $r$. Finally, if $2 \mid N$ and $a_{2}(f) = 0$, then by \cite[Theorem 4.6.17]{Miyake06}, $v_{2}(N) \geq 2$ and $N' = N$ is again divisible by $r$. In every case, we have $M \mid N'$.

Finally, assume $\ell = 2$. If $\varepsilon_{1} = \mathds{1}$, then we necessarily have $r = 4$ and $N' \geq 4$. Otherwise, if $\varepsilon_{1} \neq \mathds{1}$, we then have $N' \geq \mathfrak{c}_{1} \geq 3$ because there is no non-trivial primitive character of conductor less than $3$.
\end{pf}

It follows from \cref{bp} that if $\bar{\rho}_{f,\lambda}$ is reducible, then there exist $\mathbf{P}$ and $\mathbf{b}$ such that for all primes $p$, we have $pa_{p}(f') \equiv p^{m_{1}+1}a_{p}\left(E_{\mathbf{P}}^{\mathbf{b}}\right) \pmod{\mathfrak{L}}$. In order to use \cref{Sturm_eigen} and make these infinite set of congruences equivalent to a finite one, we now need to control the level of $E_{\mathbf{P}}^{\mathbf{b}}$.
\begin{prop}\label{np}
Let $(\varepsilon_{1},\varepsilon_{2},m_{1},m_{2}) \in R_{N,k,\varepsilon}(\mathfrak{L})$, let $p \neq \ell$ be any prime number dividing $N$, and let $b_{p} \in \{0,\varepsilon_{1}(p),p^{k'-1}\varepsilon_{2}(p)\}$.

If we have a congruence $a_{p}(f) \equiv p^{m_{1}}b_{p} \pmod{\mathfrak{L}}$, then we have
\[1 \leq v_{p}(\mathfrak{c}_{1}\mathfrak{c}_{2}) + n_{p} \leq v_{p}(N),\]
where $n_{p}$ is defined as in \cref{raise_once} with respect to $g = E$ and $b_{p}$. In particular, those inequalities are independent of the choice of $b_{p}$.
\end{prop}
\begin{pf}
First, we always have $v_{p}(\mathfrak{c}_{1}\mathfrak{c}_{2})+n_{p} \geq 1$, because $n_{p} = 0$ only if $b_{p} = a_{p}(E) = \varepsilon_{1}(p)+p^{k'-1}\varepsilon_{2}(p)$, which implies that $v_{p}(\mathfrak{c}_{1}\mathfrak{c}_{2}) \geq 1$.

Next, we claim the following:
\[b_{p} = 0 \text{ if and only if } v_{p}(\mathfrak{c}) < v_{p}(N) \text{ and } v_{p}(N) \geq 2.\]
Indeed, we have $b_{p} = 0$ if and only if $a_{p}(f) \equiv 0 \pmod{\mathfrak{L}}$. Moreover, by \cite[Theorem 4.6.17]{Miyake06} we have either $v_{p}(\mathfrak{c}) < v_{p}(N)$, $v_{p}(N) \geq 2$ and $a_{p}(f) = 0$, or $|a_{p}(f)|^{2} = p^{s}$ with $s \geq 0$. Therefore, we have $a_{p}(f) \equiv 0 \pmod{\mathfrak{L}}$ if and only if $v_{p}(\mathfrak{c}) < v_{p}(N)$ and $v_{p}(N) \geq 2$.

We now prove that $v_{p}(\mathfrak{c}_{1}\mathfrak{c}_{2})+n_{p} \leq v_{p}(N)$. If $n_{p} = 0$, it follows from the fact that $\mathfrak{c}_{1}\mathfrak{c}_{2} \mid N$. If $n_{p} = 2$, then from \cref{raise_once}, we must have $b_{p} = 0 \notin \{\varepsilon_{1}(p),p^{k'-1}\varepsilon_{2}(p)\}$. Therefore, $p \nmid \mathfrak{c}_{1}\mathfrak{c}_{2}$ and from the discussion above we have $v_{p}(N) \geq 2 = v_{p}(\mathfrak{c}_{1}\mathfrak{c}_{2})+n_{p}$.

Assume finally that $n_{p} = 1$. We then have $b_{p} \neq \varepsilon_{1}(p) + p^{k'-1}\varepsilon_{2}(p)$ and $b_{p} \in \{\varepsilon_{1}(p),p^{k'-1}\varepsilon_{2}(p)\}$. Therefore, $p$ does not divide both $\mathfrak{c}_{1}$ and $\mathfrak{c}_{2}$. If $p \nmid \mathfrak{c}_{1}\mathfrak{c}_{2}$, we have $v_{p}(\mathfrak{c}_{1}\mathfrak{c}_{2})+n_{p} = 1 \leq v_{p}(N)$. Otherwise, assume that $p \mid \mathfrak{c}_{1}$ and $p \nmid \mathfrak{c}_{2}$. We then necessarily have $b_{p} = 0$ and from the discussion above we get $v_{p}(\mathfrak{c}) < v_{p}(N)$. Looking at the $p$-part of the Artin conductor of both side of the equality $\bar{\chi}_{\ell}^{k-1}\bar{\varepsilon} = \bar{\chi}_{\ell}^{m_{1}+m_{2}}\bar{\varepsilon_{1}}\bar{\varepsilon_{2}}$, we get $v_{p}(\bar{\mathfrak{c}}) = v_{p}(\bar{\mathfrak{c}_{1}})$ where $\bar{\mathfrak{c}}$ and $\bar{\mathfrak{c}_{1}}$ denote the conductors of $\bar{\varepsilon}$ and $\bar{\varepsilon_{1}}$ respectively. Because $\varepsilon_{1}$ has prime-to-$\ell$ order, we have $\bar{\mathfrak{c}_{1}} = \mathfrak{c}_{1}$. On the other sides we always have $\bar{\mathfrak{c}} \mid \mathfrak{c}$. Therefore, we have $v_{p}(\mathfrak{c}_{1}\mathfrak{c}_{2}) = v_{p}(\mathfrak{c}_{1}) \leq v_{p}(\mathfrak{c})$. Hence, we get $v_{p}(\mathfrak{c}_{1}\mathfrak{c}_{2})+n_{p} \leq v_{p}(N)$. The case $p \mid \mathfrak{c}_{2}$ and $p \nmid \mathfrak{c}_{1}$ is exactly the same.
\end{pf}

\begin{cor}\label{E'}
Let $(\varepsilon_{1},\varepsilon_{2},m_{1},m_{2}) \in R_{N,k,\varepsilon}(\mathfrak{L})$. Define $k'$, $r$ and $E$ as in \eqref{kE0} and \eqref{rE} respectively. Consider $\bm{\mathrm{P}} \subseteq \left\{p \text{ prime, } p \mid N, p \nmid r\ell\right\}$ and $\bm{\mathrm{b}} := (b_{p})_{p \in \bm{\mathrm{P}}} \in \prod\limits_{p \in \bm{\mathrm{P}}} \{0,\varepsilon_{1}(p),p^{k'-1}\varepsilon_{2}(p)\}$ such that for all $p \in \bm{\mathrm{P}}$, we have $p^{m_{1}}b_{p} \equiv a_{p}(f) \pmod{\mathfrak{L}}$.

The modular form $E' := E_{\bm{\mathrm{P}}}^{\bm{\mathrm{b}}}$ is of weight $k'$, character $\varepsilon_{1}\varepsilon_{2}$, and its level divides $N'$. It has $\mathfrak{L}$-integral Fourier coefficients and for every prime $p$ such that either $p \nmid N\ell$ or $p \in \bm{\mathrm{P}} \cup \{r\}$, $E'$ is a normalised eigenform for the Hecke operator $T_{p}^{N'}$.
\end{cor}
\begin{pf}
From \cref{raise_several}, the form $E'$ is a normalised eigenform for all the Hecke operators at its level $N_{E'} := M \prod\limits_{p \in \bm{\mathrm{P}}} p^{n_{p}}$. Moreover, the action of $T_{p}^{N'}$ and $T_{p}^{N_{E'}}$ on $E'$ are the same if $p$ divides both $N'$ and $N_{E'}$ or none of them. If $p \nmid N\ell$, then $p \nmid N_{E'}$. If $p \in \mathbf{P} \cup \{r\}$, by \cref{np} we have
\[1 \leq v_{p}(N_{E'}) \leq v_{p}(N').\]
Therefore, $N_{E'}$ divides $N'$ and $E'$ is a normalised eigenform for the announced Hecke operators. The rest of the corollary follows from \cref{raise_several}.
\end{pf}
We now state the first main result of this section. It gives for a given $\lambda$, an explicit algorithm to check the reducibility of the representation $\bar{\rho}_{f,\lambda}$.
\begin{thm}\label{red_thm}
Let $f$ be a newform of weight $k \geq 2$, level $N \geq 1$, and character $\varepsilon$. Let $\lambda$ be a prime ideal of $K_{f}$ above a prime number $\ell$. The following assertions are equivalent:
\begin{enumerate}
\item $\bar{\rho}_{f,\lambda}$ is reducible;

\item\label{red_cong} Let $\mathfrak{L}$ be a place of $\bar{\Q}$ above $\lambda$. There exists $(\varepsilon_{1},\varepsilon_{2},m_{1},m_{2}) \in R_{N,k,\varepsilon}(\mathfrak{L})$ (see \cref{RNkeps}) such that the following holds. Let $k'$, $r$, and $N'$ be as in \eqref{kE0}, \eqref{rE} and \eqref{Nf'} respectively. Define
\[a = \left\{\begin{array}{ll}
4 & \text{if } \left|\begin{array}{l}
k \equiv m_{1}+m_{2}+3 \pmod{4},\\
\ell = 3 \text{ and } \forall p \mid N', p \equiv 1 \pmod{9};
\end{array}\right.\\
0 & \text{otherwise},
\end{array}\right. b = \left\{\begin{array}{ll}
3 & \text{if } \ell \mid N';\\
6 & \text{if } \left|\begin{array}{l}
\ell = 3 \text{ and } \forall p \mid N',\\
p \equiv 1 \pmod{9};
\end{array}\right.\\
\ell + 1 & \text{otherwise},
\end{array}\right.\]
\[\tilde{k} = a + b + \max(k,k'+bm_{1}).\]

For every prime $p \leq B := \frac{N'\tilde{k}}{12} \prod\limits_{q \mid N'} \left(1 + \frac{1}{q}\right)$ not dividing $r\ell$, we have
\begin{itemize}
\item $p \nmid N$ and $a_{p}(f) \equiv p^{m_{1}}\varepsilon_{1}(p) + p^{m_{2}}\varepsilon_{2}(p) \pmod{\mathfrak{L}}$;
\item $p \mid N$ and $a_{p}(f) \equiv p^{m_{1}} b_{p} \pmod{\mathfrak{L}}$ for some $b_{p} \in \{0,\varepsilon_{1}(p),p^{m_{2}-m_{1}}\varepsilon_{2}(p)\}$.
\end{itemize}
\end{enumerate}
When this holds, we moreover have $\bar{\rho}_{f,\lambda} \cong \bar{\chi}_{\ell}^{m_{1}}\bar{\varepsilon_{1}} \oplus \bar{\chi}_{\ell}^{m_{2}} \bar{\varepsilon_{2}}$.
\end{thm}
\begin{pf}
Assertion \ref{red_cong}. is weaker than the second part of \cref{bp}. Therefore \ref{red}. implies \ref{red_cong}.

Assume that \ref{red_cong}. holds. Consider again $k'$, $r$, $E$, $N'$, and $f'$ defined in \eqref{kE0}, \eqref{rE} and \eqref{Nf'} respectively. Define $\bm{\mathrm{P}} := \left\{p \text{ prime, } p \mid N, p \nmid \ell r, p \leq B\right\}$, and $\bm{\mathrm{b}} := (b_{p})_{p \in \bm{\mathrm{P}}}$. Finally, with the notation of \cref{raise_several}, consider the form $E' := E_{\bm{\mathrm{P}}}^{\bm{\mathrm{b}}}$.

We wish to apply \cref{Sturm_eigen} with $f = f'$, $g = E'$, $m_{f} = 1$ and $m_{g} = m_{1}+1$. By \cref{E'}, we have $E' \in \Mod_{k'}(N',\varepsilon_{1}\varepsilon_{2})$, it has $\mathfrak{L}$-integral Fourier coefficients and it is an eigenform for all the Hecke operators at level $N'$ of index less than $B$, except maybe at $\ell$. Moreover, from the identity $\bar{\chi}_{\ell}^{m_{1}+m_{2}}\bar{\varepsilon_{1}}\bar{\varepsilon_{2}} = \bar{\chi}_{\ell}^{k-1}\bar{\varepsilon}$, we have
\[\bar{\chi}_{\ell}^{k'+2(m_{1}+1)}\bar{\varepsilon_{1}}\bar{\varepsilon_{2}} = \bar{\chi}_{\ell}^{(m_{2}-m_{1}+1) + 2(m_{1}+1)}\bar{\varepsilon_{1}}\bar{\varepsilon_{2}} = \bar{\chi}_{\ell}^{m_{1}+m_{2}+3}\bar{\varepsilon_{1}}\bar{\varepsilon_{2}} = \bar{\chi}_{\ell}^{k+2}\bar{\varepsilon}.\]
Let $p$ be a prime number less than $B$. If $p \mid \ell r$, then we have $p^{m_{1}+1}a_{p}(E') \equiv 0 \equiv pa_{p}(f') \pmod{\mathfrak{L}}$. Otherwise, by \cref{raise_several} we have $p^{m_{1}+1}a_{p}(E') \equiv p^{m_{1}+1}b_{p} \equiv pa_{p}(f') \pmod{\mathfrak{L}}$. The definitions of $a$, $b$ and $\tilde{k}$ correspond to those of $a$, $b$ and $k$ in \cref{Sturm_theta} (the case $\ell = 2$ and $N' \leq 2$ never occurs as proved in \cref{f'}). By \cref{Sturm_eigen}, we therefore obtain the congruence $na_{n}(f') \equiv n^{m_{1}+1}a_{n}(E') \pmod{\mathfrak{L}}$ for every non-negative integer $n$. By \cref{bp}, we thus have $\bar{\rho}_{f,\lambda} \cong \bar{\chi}_{\ell}^{m_{1}}\bar{\varepsilon_{1}} \oplus \bar{\chi}_{\ell}^{m_{2}}\bar{\varepsilon_{2}}$.
\end{pf}

\begin{Rk}
From this theorem, we can deduce an algorithm that takes a prime ideal $\lambda$ as input and decides whether the representation $\bar{\rho}_{f,\lambda}$ is reducible or not, and computes the representation if it is reducible. In particular, it justifies the reducibility modulo $11$ of the representation treated in \cite[\nopp5.1.2]{Billerey14}. We give further details on how to explicitly do this is PARI/GP in \cref{sec_Num}. Moreover, the theorem extends the case $m = 1$ of \cite[Proposition 2.]{Kraus97}.
\end{Rk}

The previous theorem holds without any restriction on $\ell$, but the result depends on $\ell$ through
\begin{enumerate}
\item the set $R_{N,k,\varepsilon}(\mathfrak{L})$;
\item the integer $B$ that bounds the number of congruences to check.
\end{enumerate}
We first remove the dependency in $\ell$ in the set $R_{N,k,\varepsilon}(\mathfrak{L})$.

\begin{df}\label{RNeps}
Define $R_{N,\varepsilon}$ as the set of pairs $(\varepsilon_{1},\varepsilon_{2})$ of primitive Dirichlet characters such that $\varepsilon_{1}\varepsilon_{2} = \varepsilon$ and for every prime number $p$, we have $v_{p}\left(\frac{N}{\mathfrak{c}_{1}\mathfrak{c}_{2}}\right) \in \{0,1,2\}$, where $\mathfrak{c}_{i}$ is the conductor of $\varepsilon_{i}$.
\end{df}
\begin{prop}\label{bigred_params}
Assume $\ell \geq k-1$ and $\ell \nmid N\varphi(N)$. The representation $\bar{\rho}_{f,\lambda}$ is reducible if and only if there exists $(\varepsilon_{1},\varepsilon_{2}) \in R_{N,\varepsilon}$ such that $\bar{\rho}_{f,\lambda} \cong \bar{\varepsilon_{1}} \oplus \bar{\chi}_{\ell}^{k-1}\bar{\varepsilon_{2}}$. We moreover have $a_{\ell}(f) \equiv \varepsilon_{1}(\ell) + \ell^{k-1}\varepsilon_{2}(\ell) \pmod{\mathfrak{L}}$.
\end{prop}
\begin{pf}
From \cref{red_params}, if $\bar{\rho}_{f,\lambda}$ is reducible, then there exists a quadruple $(\varepsilon_{1},\varepsilon_{2},m_{1},m_{2}) \in R_{N,k,\varepsilon}(\mathfrak{L})$ such that $\bar{\rho}_{f,\lambda} \cong \bar{\chi}_{\ell}^{m_{1}}\bar{\varepsilon_{1}} \oplus \bar{\chi}_{\ell}^{m_{2}}\bar{\varepsilon_{2}}$. By the assumptions $\ell \nmid N$ and $\ell \geq k-1$, together with \cref{DF}, $f$ must be ordinary at $\lambda$ and we have an equality of sets
\[\left\{\mu\left(a_{\ell}(f)\right),\bar{\chi}_{\ell}^{k-1}\mu\left(\frac{\varepsilon(\ell)}{a_{\ell}(f)}\right)\right\} = \left\{\bar{\chi}_{\ell}^{m_{1}}\bar{\varepsilon_{1}},\bar{\chi}_{\ell}^{m_{2}}\bar{\varepsilon_{2}}\right\}.\]
We conclude that $(m_{1},m_{2}) = (0,k-1)$ and that $a_{\ell}(f) \equiv \varepsilon_{1}(\ell) \equiv \varepsilon_{1}(\ell)+\ell^{k-1}\varepsilon_{2}(\ell) \pmod{\mathfrak{L}}$. Finally, the character $\varepsilon(\varepsilon_{1}\varepsilon_{2})^{-1}$ reduces to the trivial character modulo $\mathfrak{L}$, and because $\ell \nmid \varphi(N)$, it must have prime-to-$\ell$ order. Using \cref{Rootmod}, it must be trivial and we get $\varepsilon = \varepsilon_{1}\varepsilon_{2}$.
\end{pf}
The following result will allow us to both remove the dependency in $\ell$ in the bound $B$ of \cref{red_thm}, and bound the set of $\ell$ such that $\bar{\rho}_{f,\lambda}$ is reducible.
\begin{prop}\label{big_cong}
Assume $\ell > k+1$ and $\ell \nmid N\varphi(N)$. The representation $\bar{\rho}_{f,\lambda}$ is reducible if and only if there exist a pair $(\varepsilon_{1},\varepsilon_{2}) \in R_{N,\varepsilon}$, and $\mathbf{b} \in \prod\limits_{p \in \mathbf{P}} \{0,\varepsilon_{1}(p),p^{k-1}\varepsilon_{2}(p)\}$, with $\mathbf{P} := \{p \text{ prime}, p \mid N, p \nmid \ell r\}$, such that $f' \equiv E' \pmod{\mathfrak{L}}$, with $f'$ defined as in \eqref{Nf'} and $E' := E_{\mathbf{P}}^{\mathbf{b}}$.
\end{prop}
\begin{pf}
If we have $f' \equiv E' \pmod{\mathfrak{L}}$, then in particular for all primes $p \nmid N\ell r$, we have
\[a_{p}(f) = a_{p}(f') \equiv a_{p}(E') \equiv \varepsilon_{1}(p) + p^{k-1}\varepsilon_{2}(p) \pmod{\mathfrak{L}}.\]
By \cref{bp}, $\bar{\rho}_{f,\lambda}$ is therefore reducible.

Assume that $\bar{\rho}_{f,\lambda}$ is reducible. The existence of $(\varepsilon_{1},\varepsilon_{2})$ is granted by \cref{bigred_params}. Moreover, by \cref{bp} there exists $\mathbf{b} \in \prod\limits_{p \in \mathbf{P}} \{0,\varepsilon_{1}(p),p^{k-1}\varepsilon_{2}(p)\}$ such that for every prime number $p$, we have a congruence $a_{p}(f') \equiv a_{p}(E') \pmod{\mathfrak{L}}$. By \cref{E'}, $E' \in \Mod_{k}(N',\varepsilon)$ has $\mathfrak{L}$-integral Fourier coefficients and is an eigenform for all the Hecke operators at level $N'$. By \cref{f'}, $f'$ has the same properties and therefore the modular form $f'-E'$ is constant modulo $\mathfrak{L}$. 

By the assumptions $\ell > k+1$ and $\ell \nmid N\varphi(N)$, we have $\ell \geq 5$ and $\ell \nmid N$. Therefore, we know from \cite[Theorem 12.3.7]{Diamond95} that Katz' modular form spaces with coefficients in $\bar{\F}_{\ell}$ are isomorphic to the spaces of reduction modulo $\mathfrak{L}$ of modular forms with $\mathfrak{L}$-integral coefficients. Therefore, from \cite[Corollary 4.4.2]{Katz73}, for $f'-E'$ to be congruent to a non-zero constant we must have $k \equiv 0 \pmod{\ell-1}$. This cannot hold under the assumption $\ell > k+1$ and we get $f' \equiv E' \pmod{\mathfrak{L}}$.
\end{pf}
We now state our second main result. It is analogous to \cref{red_thm} for the prime numbers $\ell > k+1$ and $\ell \nmid N\varphi(N)$.
\begin{thm}\label{bigred_thm}
Let $f$ be a newform of weight $k \geq 2$, level $N \geq 1$, and character $\varepsilon$. Let $\lambda$ be a prime ideal of $K_{f}$ above a prime number $\ell$. Assume $\ell > k+1$ and $\ell \nmid N\varphi(N)$. The following assertions are equivalent.
\begin{enumerate}
\item $\bar{\rho}_{f,\lambda}$ is reducible.
\item \label{bigred_cong} Let $\mathfrak{L}$ be a place of $\bar{\Q}$ above $\lambda$. There exists $(\varepsilon_{1},\varepsilon_{2}) \in R_{N,\varepsilon}$ such that the following holds. Let $r$ be as is \eqref{rE} (recall that $(m_{1},m_{2}) = (0,k-1)$) and let $N'$ be as in \eqref{Nf'}. Define
\[C = \left\{\begin{array}{ll}
0 & \text{if } r > 1 \text{ or } \varepsilon_{1} \neq \mathds{1};\\
-\frac{B_{k,\varepsilon_{0}}}{2k}\prod\limits_{p \mid N} a_{p}(f)(a_{p}(f)-p^{k-1}\varepsilon_{0}(p)) & \text{otherwise},
\end{array}\right.\]
where $\varepsilon_{0}$ is the primitive character associated to $\varepsilon$.

We have $C \equiv 0 \pmod{\mathfrak{L}}$, and for all primes $p \leq B := \frac{N'k}{12} \prod\limits_{q \mid N'} \left(1 + \frac{1}{q}\right)$, we have either $p \mid r$ or
\begin{itemize}
\item $a_{p}(f) \equiv \varepsilon_{1}(p)+p^{k-1}\varepsilon_{2}(p) \pmod{\mathfrak{L}}$, if $p \nmid N$;
\item $a_{p}(f) \equiv b_{p} \pmod{\mathfrak{L}}$ for some $b_{p} \in \{0,\varepsilon_{1}(p),p^{k-1}\varepsilon_{2}(p)\}$, if $p \mid N$.
\end{itemize}
\end{enumerate}
When this holds, we moreover have $\bar{\rho}_{f,\lambda} \cong \bar{\varepsilon_{1}} \oplus \bar{\chi}_{\ell}^{k-1} \bar{\varepsilon_{2}}$.
\end{thm}
\begin{pf}
Assume $\bar{\rho}_{f,\lambda}$ to be reducible. Introduce $f'$ and $E'$ as in \cref{big_cong}. The congruences for $a_{p}(f)$ follow from the congruence $f' \equiv E' \pmod{\mathfrak{L}}$. It only remains to prove that $C \equiv 0 \pmod{\mathfrak{L}}$. Because $f'$ is cuspidal, its constant coefficient at infinity is equal to $0$. Therefore, the one of $E'$ must be congruent to $0$ modulo $\mathfrak{L}$. 

The congruence $C \equiv 0 \pmod{\mathfrak{L}}$ is non-trivial only if $r = 1$ and $\varepsilon_{1} = \mathds{1}$. In this case, we have $\varepsilon_{2} = \varepsilon_{0}$, the set $\bm{\mathrm{P}}$ of \cref{big_cong} is the set of prime divisors of $N$ and for all $p \in \bm{\mathrm{P}}$ we have $b_{p} \equiv a_{p}(f) \pmod{\mathfrak{L}}$. Therefore, the constant coefficient of $E'$ is equal to
\[-\frac{B_{k,\varepsilon_{2}}}{2k} \prod\limits_{p \in \bm{\mathrm{P}}} b_{p}(b_{p}-p^{k-1}\varepsilon_{2}(p)) \equiv -\frac{B_{k,\varepsilon_{0}}}{2k} \prod\limits_{p \mid N} a_{p}(f)(a_{p}(f)-p^{k-1}\varepsilon_{0}(p)) \equiv C \pmod{\mathfrak{L}}.\]
This proves that $1.$ implies $2.$.

Assume now that the second part of the theorem holds. Consider again the modular forms $E$ and $f'$, and define $\bm{\mathrm{P}}_{\leq B} := \left\{p \text{ prime, } p \mid N, p \leq B\right\}$ and $\bm{\mathrm{b}}_{\leq B} := (b_{p})_{p \in \bm{\mathrm{P}}_{\leq B}}$, and let $E' := E_{\bm{\mathrm{P}}_{\leq B}}^{\bm{\mathrm{b}}_{\leq B}}$. By \cref{E'}, we have $E' \in \Mod_{k}(N',\varepsilon)$, it has $\mathfrak{L}$-integral coefficients and it is an eigenform for all the Hecke operators at level $N'$ of index less than $B$. The form $f'$ has moreover the same properties and for all prime numbers $p \leq B$, we have by assumption $a_{p}(f') \equiv a_{p}(E') \pmod{\mathfrak{L}}$. In order to apply \cref{Sturm_eigen} to $f = f'$, $g = E'$, $m_{f} = m_{g} = 0$, we need to have $a_{0}(E') \equiv 0 \pmod{\mathfrak{L}}$. From \cref{cst_coef}, we have $a_{0}(E') = 0$ if $\varepsilon_{1} \neq \mathds{1}$ or $r > 1$. Else, if $\varepsilon_{1} = \mathds{1}$ and $r = 1$, we have $\varepsilon_{2} = \varepsilon_{0}$ and
\[a_{0}(E') = -\frac{B_{k,\varepsilon_{2}}}{2k} \prod\limits_{p \mid N, p \leq B} b_{p}(b_{p}-p^{k-1}\varepsilon_{2}(p)) \equiv -\frac{B_{k,\varepsilon_{0}}}{2k} \prod\limits_{p \mid N, p \leq B} a_{p}(f)(a_{p}(f)-p^{k-1}\varepsilon_{0}(p)) \pmod{\mathfrak{L}}.\]
By the assumption $C \equiv 0 \pmod{\mathfrak{L}}$, we have either $a_{0}(E') \equiv 0 \pmod{\mathfrak{L}}$, or there exists $p_{0} \mid N$, $p_{0} > B$, such that $a_{p_{0}}(f)(a_{p_{0}}(f)-p_{0}^{k-1}\varepsilon_{0}(p_{0})) \equiv 0 \pmod{\mathfrak{L}}$. Define then,
\[E'' := \left\{\begin{array}{ll}
E' & \text{if } a_{0}(E') \equiv 0 \pmod{\mathfrak{L}};\\
{E'}_{p_{0}}^{b_{p_{0}}} & \text{otherwise}.
\end{array}\right.\]
\[\text{with } b_{p_{0}} = \left\{\begin{array}{ll}
0 & \text{if } a_{p_{0}}(f) \equiv 0 \pmod{\mathfrak{L}};\\
p_{0}^{k-1}\varepsilon_{0}(p_{0}) & \text{if } a_{p_{0}}(f) \equiv p_{0}^{k-1}\varepsilon_{0}(p_{0}) \pmod{\mathfrak{L}}.
\end{array}\right.\]
By \cref{E'}, $E''$ still lies in $\Mod_{k}(N',\varepsilon)$, has $\mathfrak{L}$-integral Fourier coefficients, is an eigenform for the Hecke operators at level $N'$ of index less than $B$, for any prime $p \leq B$, we have $a_{p}(E'') \equiv a_{p}(f') \pmod{\mathfrak{L}}$, and its constant Fourier coefficient vanishes modulo $\mathfrak{L}$. By \cref{Sturm_eigen}, we finally get $E'' \equiv f' \pmod{\mathfrak{L}}$ and we therefore have $\bar{\rho}_{f,\lambda} \cong \bar{\varepsilon_{1}} \oplus \bar{\chi}_{\ell}^{k-1}\bar{\varepsilon_{2}}$.
\end{pf}
\begin{Rk}
Notice that we could have always taken $r = 4$ from the start (\textit{i.e.} from \eqref{rE}) without modifying any of the results of \cref{sec_red}. The version of \cref{red_thm} and \cref{bigred_thm} we exposed in the introduction assumed that. The coefficient $C$ is then equal to zero and we get back the results announced previously.
\end{Rk}
From \cref{bigred_thm} we also deduce a bound for the reducible primes in terms of $N$, $k$ and $\varepsilon$ only.
\begin{thm}\label{red_bound}
Assume that $\bar{\rho}_{f,\lambda}$ is reducible, then one of the following condition holds.
\begin{itemize}
\item $\ell \leq k+1$;
\item $\ell \mid N\varphi(N)$;
\item there exists $(\varepsilon_{1}$, $\varepsilon_{2}) \in R_{N,\varepsilon}$ such that $\ell$ divides the algebraic norm of one the following non-zero quantities
\begin{enumerate}
\item $B_{k,(\varepsilon_{1}^{-1}\varepsilon_{2})_{0}}$;
\item $p^{k}-(\varepsilon_{1}\varepsilon_{2}^{-1})_{0}(p)$ for a prime $p$ such that $p \mid \mathfrak{c}_{1}\mathfrak{c}_{2}$, $p \nmid \mathfrak{c}_{0}$ with $\mathfrak{c}_{0}$ the conductor of $(\varepsilon_{1}\varepsilon_{2}^{-1})_{0}$.
\end{enumerate}
\end{itemize}
\end{thm}
\begin{pf}
Assume $\ell > k+1$ and $\ell \nmid N\varphi(N)$. From \cref{big_cong}, if $\bar{\rho}_{f,\lambda}$ is reducible, we have a congruence modulo $\mathfrak{L}$ between the cuspidal modular form $f'$, and the Eisenstein series $E'$. Therefore, by Katz' $q$-expansion principle (see \cite[]{Katz73}) the constant coefficient of $E'$ must be congruent to $0$ modulo $\mathfrak{L}$ at every cusp. By \cref{Upsilon}, the constant coefficient of $E'$ at the cusp $\frac{1}{\mathfrak{c}_{2}}$ divides the quantity
\[-\varepsilon_{1}(-1)\frac{W((\varepsilon_{1}\varepsilon_{2}^{-1})_{0})}{W(\varepsilon_{2}^{-1})}
	\frac{B_{k,(\varepsilon_{1}^{-1}\varepsilon_{2})_{0}}}{2k}  \left(\frac{\mathfrak{c}_{2}}{\mathfrak{c}_{0}}\right)^{k} \prod_{p \mid \mathfrak{c}_{1}\mathfrak{c}_{2}} \left(1 - \frac{(\varepsilon_{1}\varepsilon_{2}^{-1})_{0}(p)}{p^{k}}\right) \prod_{p \mid N'} \left(1 - \frac{\varepsilon_{1}(p)\varepsilon_{2}^{-1}(p)}{p^{k}}\right)\left(1 - \frac{1}{p}\right).\]
Let us look at the prime factors of the norm of this coefficient.
\begin{itemize}
\item The number $-\varepsilon_{1}(-1)$ is a unit. Its norm has no prime factor.

\item By \cref{W_divisors}, the prime factors of the norm of $\frac{W((\varepsilon_{1}\varepsilon_{2}^{-1})_{0})}{W(\bar{\varepsilon_{2}})}\big(\frac{\mathfrak{c}_{2}}{\mathfrak{c}_{0}}\big)^{k}$ are only powers of prime factors of $N$. By assumption, $\ell$ does not divide them.

\item For $p \mid N'$, we have $1 - \frac{1}{p} = \frac{p-1}{p}$. By the assumption $\ell \nmid N\varphi(N)$, this cannot vanish modulo $\mathfrak{L}$.

\item For $p \mid N'$ again, let us prove that the prime factors of the norm of $\left(1 - \frac{\varepsilon_{1}(p)\varepsilon_{2}^{-1}(p)}{p^{k}}\right)$ are redundant with the ones of $N$ and $\left(p^{k}-(\varepsilon_{1}\varepsilon_{2}^{-1})_{0}(p)\right)$. Note that we either have $p = 2$ and $(k,\varepsilon_{1},\varepsilon_{2}) = (2,\mathds{1},\mathds{1})$, or $N' = N$. In the first case we have $1-\frac{\varepsilon_{1}(p)\varepsilon_{2}^{-1}(p)}{p^{k}} = \frac{3}{4}$. This cannot vanish modulo $\mathfrak{L}$ by assumption because $2$ and $3$ are less or equal to $k+1 = 3$. Otherwise, we have $p \mid N$, then either $p \mid \mathfrak{c}_{0}$ and $1-\frac{\varepsilon_{1}(p)\varepsilon_{2}^{-1}(p)}{p^{k}} = 1 \not\equiv 0 \pmod{\mathfrak{L}}$, or $p \nmid \mathfrak{c}_{0}$ and $1-\frac{\varepsilon_{1}(p)\varepsilon_{2}^{-1}(p)}{p^{k}} = \frac{p^{k}-(\varepsilon_{1}\varepsilon_{2}^{-1})_{0}(p)}{p^{k}}$. Therefore, $\ell$ must divide the algebraic norm of $p^{k}-(\varepsilon_{1}\varepsilon_{2}^{-1})_{0}(p)$.

\item Finally $\ell$ divides either the norm of $\frac{B_{k,(\varepsilon_{1}^{-1}\varepsilon_{2})_{0}}}{2k}$ and thus $B_{k,(\varepsilon_{1}^{-1}\varepsilon_{2})_{0}}$ because $2k$ is non-zero modulo $\mathfrak{L}$ by assumption, or $p^{k}-(\varepsilon_{1}\varepsilon_{2}^{-1})_{0}(p)$ for $p \mid \mathfrak{c}_{1}\mathfrak{c}_{2}$. This final quantity contains only prime factors of $N$ if $p \mid \mathfrak{c}_{0}$. We can therefore consider only the primes $p \nmid \mathfrak{c}_{0}$.
\end{itemize}
\end{pf}

\section{Dihedral representations}\label{sec_dihedral}
Let $k \geq 2$, $N \geq 1$ be two integers, and $\varepsilon$ be a Dirichlet character modulo $N$ of conductor $\mathfrak{c}$. Let $f \in \Smod_{k}^{\mathrm{new}}(N,\varepsilon)$ be a newform. We study in this section the case where $\bar{\rho}_{f,\lambda}$ has projective dihedral image. The main result of this section is \cref{Dihthm}. The strategy used to prove this result is similar to the one of \cite[Theorem 3.1]{Billerey14} but several technicalities arise when dealing with forms with non-trivial characters. However, our method is based on a more refined use of the local description of $\bar{\rho}_{f,\lambda}$ at the prime dividing the level, leading to an improved bound in the trivial character case.

\subsection{CM forms}
Let $\varphi$ be a Dirichlet character of modulus $M$. We define the twist of $f$ by $\varphi$ as the only newform $f \otimes \varphi$ such that $a_{p}(f \otimes \varphi) = \varphi(p)a_{p}(f)$ for all but finitely many primes $p$. We have the following result from \cite[§§1-3]{Atkin78}.

\begin{prop}\label{Atkin}
With the above notations we have $f \otimes \varphi \in \Smod_{k}(\lcm(N,M^{2},\mathfrak{c}M),\varphi^{2}\varepsilon)$. For all primes $p \nmid M$, we have $a_{p}(f \otimes \varphi) = \varphi(p)a_{p}(f)$ and the $p$-part of the level of $f \otimes \varphi$ is equal to $p^{v_{p}(M)}$.

Moreover, if $v_{p}(N) = v_{p}(\mathfrak{c})$ and $\varphi_{p} = \varepsilon_{p}^{-1}$, where $\varphi_{p}$ and $\varepsilon_{p}$ denote the $p$-parts of $\varphi$ and $\varepsilon$ respectively, then the $p$-parts of the levels of $f \otimes \varphi$ and $f$ are equal and we have $a_{p}(f \otimes \varphi) = \varepsilon_{p}'(p)\varphi_{p}'(p)\bar{a_{p}(f)}$, with $\bar{a_{p}(f)}$ the complex conjugate of $a_{p}(f)$ and $\varepsilon_{p}'$ and $\varphi_{p}'$ are the prime-to-$p$ part of $\varepsilon$ and $\varphi$ respectively.
\end{prop}
We take this definition of CM forms from \cite[p. 34]{Ribet77}.

\begin{df}\label{CM}
Suppose $\varphi$ is not the trivial character. The form $f$ is said to have complex multiplication by $\varphi$ if $\varphi(p)a_{p}(f) = a_{p}(f)$ for all primes $p$ is a set of primes of density $1$.
\end{df}

\subsection{Study of the ramification of a character}
Let $\ell$ be a prime number and let $\lambda$ be a prime ideal of $\mathcal{O}_{f}$ above $\ell$. We assume for the rest of the section that $\ell \neq 2$. Consider the projectivisation $\P\bar{\rho}_{f,\lambda} : G_{\Q} \overset{\bar{\rho}_{f,\lambda}}{\longrightarrow} \GL_{2}(\F_{\lambda}) \to \PGL_{2}(\F_{\lambda})$,
and assume that $\P\bar{\rho}_{f,\lambda}(G_{\Q})$ is a dihedral group $D_{2n}$ of order $2n$ with $\ell \nmid n$. Let $C$ be the unique cyclic subgroup of order $n$ in $\P\bar{\rho}_{f,\lambda}(G_{\Q})$. It is a normal subgroup, and we get a quadratic character:
\[\theta_{f,\lambda} : G_{\Q} \overset{\P\bar{\rho}_{f,\lambda}}{\longrightarrow} D_{2n} \to D_{2n}/C \cong \Z/2\Z.\]
Recall that the elements in $D_{2n} \setminus C$ are of order $2$ and hence have all trace $0$. Let us focus on the ramification of the character $\theta_{f,\lambda}$. Let $\mathfrak{c}_{f,\lambda}$ be the conductor of $\theta_{f,\lambda}$ seen as a Dirichlet character. Because $\bar{\rho}_{f,\lambda}$ is unramified outside $N\ell$, the prime factors of $\mathfrak{c}_{f,\lambda}$ are among the ones of $N\ell$. Let $K_{f,\lambda}$ be the number field fixed by the kernel of $\theta_{f,\lambda}$. It is a quadratic extension of $\Q$, and we can write $K_{f,\lambda} = \Q(\sqrt{d_{f,\lambda}})$ with $d_{f,\lambda} \in \Z$ square-free. By the conductor-discriminant formula \cite[VII. (11.9)]{Neukirch99}, we have
\[\mathfrak{c}_{f,\lambda} = \left\{\begin{array}{ll}
|d_{f,\lambda}| & \text{if } d_{f,\lambda} \equiv 1 \pmod{4}\\
4|d_{f,\lambda}| & \text{if } d_{f,\lambda} \equiv 2,3 \pmod{4}
\end{array}\right..\]
We thus have the following result.
\begin{prop}\label{ramif2}
Let $p$ be a prime number dividing $\mathfrak{c}_{f,\lambda}$. Then we have either $p \neq 2$, $p \mid N\ell$ and $v_{p}(\mathfrak{c}_{f,\lambda}) = 1$, or $p = 2$, $p \mid N$ and $v_{p}(\mathfrak{c}_{f,\lambda}) \in \{2,3\}$.
\end{prop}
We now prove the following lemma.

\begin{lm}\label{lmimage}
Let $G$ be a finite group and let $K$ be a field a positive characteristic $\ell$. Let $\rho : G \to \GL_{2}(K)$ be a morphism with values in the subgroup of $\GL_{2}(K)$ of upper triangular matrices. Put $\rho = \begin{pmatrix}
\chi_{1} & c\\
0 & \chi_{2}
\end{pmatrix}$. If the order of $\P\rho$ is prime to $\ell$, then $\Im(\P\rho) \cong \Im(\chi_{1}\chi_{2}^{-1})$.
\end{lm}
\begin{pf}
From the isomorphisms $G/\ker(\chi_{1}\chi_{2}^{-1}) \cong \Im(\chi_{1}\chi_{2}^{-1})$ and $G/\ker(\P\rho) \cong \Im(\P\rho)$, it suffices to prove that $\ker(\P\rho) = \ker(\chi_{1}\chi_{2}^{-1})$.

Let $\sigma \in \ker(\P\rho)$. There exists $\lambda \in K$ such that $\rho(\sigma) = \lambda I_{2}$, and we immediately have $\chi_{1}(\sigma) = \chi_{2}(\sigma)$. Therefore, $\sigma \in \ker(\chi_{1}\chi_{2}^{-1})$.

Let $\sigma \in \ker(\chi_{1}\chi_{2}^{-1})$ and let $\lambda := \chi_{1}(\sigma) = \chi_{2}(\sigma)$. We then have
\[\rho(\sigma) = \begin{pmatrix}
\lambda & c(\sigma)\\
0 & \lambda
\end{pmatrix}.\]
The order of $\rho(\sigma)$ in $\PGL_{2}(K)$ is equal to $1$ if $c(\sigma) = 0$ and to $\ell$ if not. By assumption, the image of $\P\rho$ is of order prime to $\ell$. Thus $c(\sigma) = 0$ and $\sigma \in \ker(\P\rho)$.
\end{pf}
We first study the ramification of $\theta_{f,\lambda}$ at $\ell$.

\begin{prop}\label{ramif_l}
Assume $\ell \nmid N$ and $\ell > k$. Then, with the terminology of \cref{DF} we have the following.
\begin{enumerate}
\item If $f$ is ordinary at $\lambda$ and $\ell \neq 2k - 1$, then $\theta_{f,\lambda}$ is unramified at $\ell$;
\item If $f$ is not ordinary at $\lambda$ and $\ell \neq 2k - 3$, then $\theta_{f,\lambda}$ is unramified at $\ell$.
\end{enumerate}
\end{prop}
\begin{pf}
We are under the hypotheses of \cref{DF}. 
\begin{itemize}
\item If $f$ is ordinary at $\lambda$, then we have
\[\bar{\rho}_{f,\lambda|I_{\ell}} \cong \begin{pmatrix}
\bar{\chi}_{\ell}^{k-1} & \star\\
0 & 1
\end{pmatrix}.\]
The order of the image of $\P\bar{\rho}_{f,\lambda}$ is prime to $\ell$. Therefore, so is the one of $\bar{\rho}_{f,\lambda|I_{\ell}}$. By \cref{lmimage}, we have $\Im(\P\bar{\rho}_{f,\lambda}) \cong \Im({\bar{\chi}_{\ell}^{k-1}})$. The image of $\bar{\chi}_{\ell}^{k-1}$ is a cyclic group of order $\frac{\ell - 1}{\gcd(\ell-1,k-1)}$. If $\ell > k$ and $\ell \neq 2k-1$, it is greater than $2$. Therefore, the image of $\P\bar{\rho}_{f,\lambda}$ is necessarily included in $C$ and $\theta_{f,\lambda}$ is unramified at $\ell$.

\item If $f$ is not ordinary at $\lambda$, then we have
\[\bar{\rho}_{f,\lambda|I_{\ell}} \cong \begin{pmatrix}
\psi^{k-1} & 0\\
0 & \psi^{\ell(k-1)}
\end{pmatrix}.\]
By \cref{lmimage}, $\Im(\P\bar{\rho}_{f,\lambda})$ is isomorphic to the image of $\psi^{(k-1)(\ell-1)}$. It is of order $\frac{\ell+1}{\gcd(\ell+1,k-1)}$ which is again greater than $2$ if we assume $\ell > k$ and $\ell \neq 2k - 3$. We conclude as before.
\end{itemize}
\end{pf}
We now look at the ramification at the primes $p \mid N$, $p \neq \ell$.

\begin{prop}\label{ramif_N}
Let $p$ be a prime number dividing $N$ and different from $\ell$.
\begin{enumerate}
\item If $v_{p}(N) = 1$ and $v_{p}(\mathfrak{c}) = 0$, then $\theta_{f,\lambda}$ is unramified at $p$.
\item Assume $v_{p}(N) = v_{p}(\mathfrak{c})$. If $\theta_{f,\lambda}$ is ramified at $p$, then the $p$-parts of $\varepsilon$ and $\theta_{f,\lambda}$ are equal modulo $\lambda$ (in particular, $\varepsilon_{p}$ has order $2$ modulo $\lambda$).
\end{enumerate}
\end{prop}

\begin{pf}
\begin{enumerate}
\item By \cref{LW} the restriction of $\bar{\rho}_{f,\lambda}$ at an inertia subgroup $I_{p}$ at $p$, is given by
\[\bar{\rho}_{f,\lambda|I_{p}} \cong \begin{pmatrix}
\mathds{1} & \star\\
0 & \mathds{1}
\end{pmatrix}.\]
By assumption, the order of image of $\P\bar{\rho}_{f,\lambda}$ is prime to $\ell$. Therefore, by \cref{lmimage} $\Im(\P\bar{\rho}_{f,\lambda})$ is trivial, and $\P\bar{\rho}_{f,\lambda}$ is unramified at $p$, as well as $\theta_{f,\lambda}$.

\item By \cref{LW}, we have
\[\bar{\rho}_{f,\lambda|I_{p}} \cong \mathds{1} \oplus \bar{\varepsilon}_{|I_{p}}.\]
As a character of $G_{\Q}$, $\varepsilon$ factors through the group $\Gal(\Q(\zeta_{N})/\Q)$. Moreover, an inertia subgroup at $p$ in $\Gal(\Q(\zeta_{N})/\Q)$ is given by the group $\Gal(\Q(\zeta_{N})/\Q(\zeta_{p^{v_{p}(N)}})) \cong \left(\Z/p^{v_{p}(N)}\Z\right)^{\times}$. Thus, as a Dirichlet character, the restriction of $\varepsilon$ to $I_{p}$ is $\varepsilon_{p}$. By \cref{lmimage}, we thus have $\P\bar{\rho}_{f,\lambda}(I_{p}) \cong \bar{\varepsilon_{p}}(I_{p})$. If $\theta_{f,\lambda}$ is ramified at $p$, then its image must have a non-trivial intersection with $D_{2n} \setminus C$. Therefore, it is isomorphic to $\Z/2\Z$. We deduce that the reduction of $\varepsilon_{p}$ modulo $\lambda$ is of order $2$, and that the $p$-part of $\theta_{f,\lambda}$ corresponds to the reduction of $\varepsilon_{p}$.
\end{enumerate}
\end{pf}

\subsection{Proof of the main result}
In order to prove the main result of this section, consider the form
\[g = f \otimes \theta_{f,\lambda}.\]
It is by construction a newform of weight $k$, character $\varepsilon$, and we have
\begin{prop}\label{twist_isom}
The Galois representations $\bar{\rho}_{f,\lambda}$ and $\bar{\rho}_{g,\lambda}$ are isomorphic.
\end{prop}
\begin{pf}
Because $f$ and $g$ have same weight and character, the determinants of $\bar{\rho}_{f,\lambda}$ and $\bar{\rho}_{g,\lambda}$ are the same. Let $p$ be a prime not dividing $N\ell$. By \cref{Atkin} and \cref{ramif2}, the representations are both unramified at $p$ and their traces at a Frobenius element at $p$ are the reductions modulo $\lambda$ of $a_{p}(f)$ and $a_{p}(g) = \theta_{f,\lambda}(p)a_{p}(f)$ respectively. If $\theta_{f,\lambda}(p) = 1$, then the traces agree. Else, if $\theta_{f,\lambda}(p) = -1$, then by definition of $\theta_{f,\lambda}$ we have $\Tr(\bar{\rho}_{f,\lambda}(\mathrm{Frob}_{p})) = 0$ and the traces agree again. By Brauer-Nesbitt theorem \cite[Lemme 3.2.]{Deligne74}, we have $\bar{\rho}_{f,\lambda} \cong \bar{\rho}_{g,\lambda}$.
\end{pf}
From now on, let assume that $\ell > k$, $\ell \notin \{2k-1,2k-3\}$. From \cref{ramif_l}, it means that the character $\theta_{f,\lambda}$ is unramified at $\ell$. The following proposition contains much information about the level and the Fourier coefficients of $g$.
\begin{prop}\label{g}
Assume that $\ell \nmid N\varphi(N)$. Write $N_{g}$ for the level of $g$ and let $p$ be prime number.
\begin{enumerate}[(i)]
\item If either $p \nmid N$, or $v_{p}(N) = 1$ and $v_{p}(\mathfrak{c}) = 0$, then $v_{p}(N_{g}) = v_{p}(N)$ and $a_{p}(g) = \theta_{f,\lambda}(p)a_{p}(f)$.

\item If $v_{p}(N) = v_{p}(\mathfrak{c}) \geq 1$, then $v_{p}(N_{g}) = v_{p}(N)$ and $a_{p}(g) = \theta_{f,\lambda}(p)a_{p}(f)$ if $\theta_{f,\lambda}$ is unramified at $p$, and $a_{p}(g) = \bar{a_{p}(f)}\varepsilon_{p}'(p)\left(\theta_{f,\lambda}\right)_{p}'(p)$ otherwise, where $\varepsilon_{p}'$ and $\left(\theta_{f,\lambda}\right)_{p}'$ are the prime-to-$p$ parts of $\varepsilon$ and $\theta_{f,\lambda}$ respectively.

\item If $v_{p}(N) \geq 2$ and $v_{p}(\mathfrak{c}) < v_{p}(N)$, then $v_{p}(N_{g}) \leq v_{p}(N) + 2\min(v_{p}(N),v_{p}(2))$.
\end{enumerate}
In particular, we have $N_{g} \mid N\gcd(N,2)^{2}$.
\end{prop}
\begin{pf}
\begin{enumerate}[(i)]
\item If either $p \nmid N$, or $v_{p}(N) = 1$ and $v_{p}(\mathfrak{c}) = 0$, then by \cref{ramif_N}, $\theta_{f,\lambda}$ is unramified at $p$. By \cref{Atkin}, we have $v_{p}(N_{g}) = v_{p}(N)$ and $a_{p}(g) = \theta_{f,\lambda}(p)a_{p}(f)$.

\item If $v_{p}(N) = v_{p}(\mathfrak{c})$. Then either $\theta_{f,\lambda}$ is unramified at $p$, and we again have $v_{p}(N_{g}) = v_{p}(N)$ and $a_{p}(g) = \theta_{f,\lambda}(p)a_{p}(f)$, or it is ramified. In this latter case, by \cref{ramif_N}, we have $\varepsilon_{p} \equiv \left(\theta_{f,\lambda}\right)_{p} \pmod{\lambda}$. As we assumed $\ell \nmid \varphi(N)$, by \cref{Rootmod}, we get that $\varepsilon_{p}$ has order $2$ and that $\left(\theta_{f,\lambda}\right)_{p} = \varepsilon_{p} = \varepsilon_{p}^{-1}$. Therefore by \cref{Atkin}, we have $v_{p}(N_{g}) = v_{p}(N)$ and $a_{p}(g) = \bar{a_{p}(f)}\varepsilon_{p}'(p)\left(\theta_{f,\lambda}\right)_{p}'(p)$.

\item Finally assume $v_{p}(N) \geq 2$ and $v_{p}(\mathfrak{c}) < v_{p}(N)$. If $p \neq 2$, then we have $v_{p}(\mathfrak{c}_{f,\lambda}) = 1$ and by \cref{Atkin} we find that $v_{p}(N_{g}) \leq v_{p}(N)$. Assume that $p = 2$ and let us look at the Artin conductors of $\bar{\rho}_{f,\lambda}$ and $\bar{\rho}_{g,\lambda}$. By \cref{twist_isom}, they are equal. Moreover, by \cref{CL} the difference between the $2$-adic valuations of $N$ and $N(\bar{\rho}_{f,\lambda})$, and of $N_{g}$ and $N(\bar{\rho}_{g,\lambda})$ respectively, cannot be greater than $2$. Therefore, we have
\[v_{2}(N{_g}) \leq v_{2}(N(\bar{\rho}_{g,\lambda})) + 2 = v_{2}(N(\bar{\rho}_{f,\lambda})) + 2 \leq v_{2}(N)+2.\]
\end{enumerate}
\end{pf}
Finally consider the set $\bm{\mathrm{P}} = \{p \text{ prime}, v_{p}(N) \geq 2 \text{ and } v_{p}(\mathfrak{c}) < v_{p}(N)\}$, and with the notations of \cref{raise_several}, define
\[h := g_{\bm{\mathrm{P}}}^{(0)_{p \in \bm{\mathrm{P}}}}.\]

\begin{prop}\label{h}
Assume that $\ell \nmid N\varphi(N)$. The form $h$ is of weight $k$, character $\varepsilon$, and level dividing $N\gcd(N,2)^{2}$. Moreover, it is a normalised eigenform for all the Hecke operators at level $N\gcd(N,2)^{2}$ and for all prime numbers $p$, we have
\[a_{p}(h) = \left\{\begin{array}{ll}
a_{p}(g) & \text{if } p \notin \bm{\mathrm{P}}\\
0 & \text{if } p \in \bm{\mathrm{P}}
\end{array}\right..\]
\end{prop}
\begin{pf}
The computation of the weight and the character of $h$ follows from \cref{raise_several} as well as Fourier coefficients of $h$. Write $N_{h}$ for the level of $h$ and $N_{g}$ for the level of $g$. With the notations of \cref{raise_several}, we have $N_{h} = N_{g}\prod\limits_{p \in \bm{\mathrm{P}}} p^{n_{p}}$ and $h$ is an eigenform for all the Hecke operators at level $N_{h}$. Let us prove that $N_{h}$ divides $N\gcd(N,2)^{2}$.

If $n_{p} = 0$, there is nothing to do. If $n_{p} = 2$, then $p$ does not divide the level of $g$. By definition of $\bm{\mathrm{P}}$ we have $v_{p}(N_{h}) = v_{p}(N_{g})+2 = 2 \leq v_{p}(N)$. If $n_{p} = 1$, then we have $p \mid N_{g}$ and $a_{p}(g) \neq 0$. Moreover, the character of $g$ is $\varepsilon$. Therefore, by \cite[Theorem 4.6.17]{Miyake06} we cannot have $v_{p}(N_{g}) \geq 2$ and $v_{p}(N_{g}) > v_{p}(\mathfrak{c})$. However, by assumption we have $v_{p}(N) \geq 2$ and $v_{p}(N) > v_{p}(\mathfrak{c})$. We conclude that either $v_{p}(N_{g}) = v_{p}(\mathfrak{c}) < v_{p}(N)$, or $v_{p}(N_{g}) \leq 1 < v_{p}(N)$. In every case we have $v_{p}(N_{h}) = v_{p}(N_{g}) + 1 \leq v_{p}(N)$.

Finally, to prove that $h$ is an eigenform for all the Hecke operators at level $N\gcd(N,2)^{2}$, we just have to prove that the prime divisors of $N_{h}$ and $N\gcd(N,2)^{2}$ are the same. This is clear from \cref{g} for the prime numbers outside $\bm{\mathrm{P}}$. If $p \in \bm{\mathrm{P}}$, then the $p$-adic valuation of the level of $h$ is also positive by \cref{raise_several}. This finishes the proof.
\end{pf}
We now prove our main result.
\begin{thm}\label{Dihthm}
Assume $\bar{\rho}_{f,\lambda}$ has dihedral projective image. If $N = 1$, then we have $\ell \leq k$ or $\ell \in \{2k - 1,2k-3\}$. Else, if $N \geq 2$ and $f$ does not have complex multiplication, then we have
\[\ell \leq \left(2N^{\frac{k-1}{2}}\right)^{[K_{f}:\Q]} \times \max\left(
    \left(\frac{k}{3}\left(2\log\log(N)+2.4\right)\right)^{\frac{k-1}{2}},
    \left(\frac{5}{2}N^{\frac{k-1}{2}}\right)
\right)^{[K_{f}:\Q]}.\]
\end{thm}
\begin{pf}
Assume that $\ell > k$, $\ell \notin \{2k-1;2k-3\}$ and $\ell \nmid N\varphi(N)$. If $N = 1$, then $\theta_{f,\lambda}$ is unramified everywhere and thus trivial, which is a contradiction. Assume now that $N \geq 2$. Because $f$ does not have complex multiplication, there exists some prime number $p$ such that $a_{p}(f) \neq a_{p}(h)$. Moreover, from \cref{h}, $f$ and $h$ are both modular forms of weight $k$, level $N\gcd(N,2)^{2}$, and character $\varepsilon$, that are eigenforms for all the Hecke operators at this level. Therefore, from \cite{Murty97}, $p$ must be less than $B = \frac{kN\gcd(N,2)^{2}}{12} \prod\limits_{\underset{p \text{ prime}}{p \mid N}} \left(1+\frac{1}{p}\right)$. We treat various cases for $p$.
\begin{itemize}
\item If $p = \ell$, then we have $\ell \leq B$.

\item If $p \nmid N\ell$, then from \cref{g} and \cref{h}, we have $a_{p}(h) = a_{p}(g) = \theta_{f,\lambda}(p)a_{p}(f)$. This implies that $\theta_{f,\lambda}(p) = -1$ and by \cref{twist_isom}, we get $a_{p}(f) \equiv 0 \pmod{\lambda}$.

\item If $v_{p}(N) \geq 2$ and $v_{p}(N) > v_{p}(\mathfrak{c})$, then by \cite[Theorem 4.6.17]{Miyake06} and \cref{h}, we have $a_{p}(f) = 0 = a_{p}(h)$. We cannot have $a_{p}(f) \neq a_{p}(h)$ in this situation.

\item If $v_{p}(N) = 1$ and $v_{p}(\mathfrak{c}) = 0$, then by \cref{g} and \cref{h}, we have $a_{p}(h) = a_{p}(g) = \theta_{f,\lambda}(p)a_{p}(f)$ and $v_{p}(N_{g}) = 1$. Therefore, we have $\theta_{f,\lambda}(p) = -1$ and by \cref{LW} and \cref{twist_isom}, we get
\[\begin{pmatrix}
pa_{p}(f) & \star\\
0 & a_{p}(f)
\end{pmatrix} \cong \begin{pmatrix}
-pa_{p}(f) & \star\\
0 & -a_{p}(f)
\end{pmatrix}.\]
Hence, we either have $a_{p}(f) \equiv 0 \pmod{\lambda}$, or $p+1 \equiv 0 \pmod{\ell}$.

\item Finally assume that $v_{p}(N) = v_{p}(\mathfrak{c})$. By \cref{g} and \cref{h}, we have $a_{p}(h) = a_{p}(g)$ and $v_{p}(N_{g}) = v_{p}(N)$. Therefore, by \cref{twist_isom} and \cref{LW}, we have an isomorphism of representations of $G_{p}$:
\[\mu(a_{p}(f)) \oplus \mu(a_{p}(f))\bar{\chi}_{\ell}^{k-1}\bar{\varepsilon_{|G_{p}}} \cong \mu(a_{p}(h)) \oplus \mu(a_{p}(h))\bar{\chi}_{\ell}^{k-1}\bar{\varepsilon_{|G_{p}}}.\]
The second character of each representation is ramified at $p$, while the first is not. We deduce that we have $a_{p}(f) \equiv a_{p}(h) \pmod{\mathfrak{L}}$.

If $\theta_{f,\lambda}$ is unramified at $p$, we have by \cref{g}, $a_{p}(h) = \theta_{f,\lambda}(p)a_{p}(f)$, and therefore $\theta_{f,\lambda}(p) = -1$ and $a_{p}(f) \equiv 0 \pmod{\lambda}$.

If $\theta_{f,\lambda}$ is ramified at $p$, we have from \cref{g} and \cref{h},
\[a_{p}(h) = \varepsilon_{p}'(p)\left(\theta_{f,\lambda}\right)_{p}'(p)\bar{a_{p}(f)} \equiv a_{p}(f) \pmod{\lambda}.\]
We moreover know from \cite[Theorem 4.6.17]{Miyake06} that $a_{p}(f)\bar{a_{p}(f)} = p^{k-1}\varepsilon_{p}'(p)$. Therefore we get that $a_{p}(f)^{2} \equiv \varepsilon_{p}'(p)^{2}\left(\theta_{f,\lambda}\right)_{p}'(p)p^{k-1} \pmod{\lambda}$.
\end{itemize}
To sum up, we have either $\ell \leq B$, or $a_{p}(f) \equiv 0 \pmod{\lambda}$, or $p \mid N$ and $p+1 \equiv 0 \pmod{\ell}$, or $p \mid N$ and $a_{p}(f)^{2} \equiv \varepsilon'_{p}(p)^{2}\left(\theta_{f,\lambda}\right)_{p}'(p)p^{k-1} \pmod{\lambda}$. Using Deligne's bounds for the coefficients of a newform and \cref{Sturm_upperbound}, this means that either $\ell \leq N+1$, or
\[\ell \mid |\mathrm{Norm}(a_{p}(f))| = \prod\limits_{\sigma: K_{f} \hookrightarrow \C} |\sigma(a_{p}(f))| \leq \left(2p^{\frac{k-1}{2}}\right)^{[K_{f}:\Q]} \leq \left(2N^{\frac{k-1}{2}}\right)^{[K_{f}:\Q]} \left(\frac{k}{3}\left(2\log\log(N) + 2.4\right)\right)^{\frac{k-1}{2}[K_{f}:\Q]},\]
or
\begin{align*}
\ell \mid \left|\mathrm{Norm}\left(a_{p}(f)^{2} -\varepsilon_{p}'(p)^{2}\left(\theta_{f,\lambda}\right)_{p}'(p)p^{k-1}\right)\right| &\leq \prod\limits_{\sigma: K_{f} \hookrightarrow \C} \left(\left|\sigma\left(a_{p}(f)^{2}\right)\right|+p^{k-1}\right)\\
    &\leq \left(5p^{k-1}\right)^{[K_{f}:\Q]}\\
    &\leq \left(5N^{k-1}\right)^{[K_{f}:\Q]} = \left(2N^{\frac{k-1}{2}}\right)^{[K_{f}:\Q]} \times \left(\frac{5}{2}N^{\frac{k-1}{2}}\right)^{[K_{f}:\Q]}.
\end{align*}
We therefore get the wanted result.
\end{pf}

\section{Numerical applications}\label{sec_Num}

\subsection{Checking the reducibility}
We explain here how to use \cref{red_thm} and \cref{bigred_thm} to explicitly compute the prime ideals $\lambda$ for which the representation $\bar{\rho}_{f,\lambda}$ is reducible. We begin by discussing the dependency of the set $R_{N,k,\varepsilon}(\mathfrak{L})$ (see \cref{RNkeps}) in the place $\mathfrak{L}$.

\begin{prop}\label{RNkeps_dep}
Let $N \geq 1$ and $k \geq 2$ be integers, and let $\varepsilon$ be a Dirichlet character modulo $N$. Let $\ell$ be a prime number and let $\mathfrak{L}$ be a place of $\bar{\Q}$ above $\ell$. The set $R_{N,k,\varepsilon}(\mathfrak{L})$ depends only on $\mathfrak{L} \cap \Q(\varepsilon)$ (and on $N$, $k$ and $\varepsilon$).
\end{prop}
\begin{pf}
Write $\pi_{\mathfrak{L}}$ for the projection modulo $\mathfrak{L}$, and $T_{\mathfrak{L}}$ for the associated Teichmüller lift (see \cref{sec_Tlifts}). Recall that for $x \in \bar{\F}_{\ell}^{\times}$, $T_{\mathfrak{L}}(x)$ is the only root of unity of order prime to $\ell$ and such that $\pi_{\mathfrak{L}}(T_{\mathfrak{L}}(x)) = x$.

We first prove that the map $T_{\mathfrak{L}} \circ \pi_{\mathfrak{L}}$ depends only on $\ell$. Let $\zeta$ be a root of unity of order $n = \ell^{m}q$ with $m \geq 0$ and $\ell \nmid q$. We can then write $\zeta = \zeta^{\ell^{m}a} \cdot \zeta^{qb}$, with $\ell \nmid b$ and $a$ prime to $q$. Because $\zeta$ is a root of unity of order $n$, $\zeta^{\ell^{m}a}$ is a root of unity of order $q$ and $\zeta^{qb}$ is a root of unity of order $\ell^{m}$. From \cref{Rootmod}, we get $\zeta \equiv \zeta^{\ell^{m}a} \pmod{\mathfrak{L}}$, and $T_{\mathfrak{L}}\circ \pi_{\mathfrak{L}}(\zeta) = \zeta^{\ell^{m}a}$. Therefore, it depends only on $\ell$.

Let $(\varepsilon_{1},\varepsilon_{2},m_{1},m_{2}) \in R_{N,k,\varepsilon}(\mathfrak{L})$. The only dependency on the place $\mathfrak{L}$ is the congruence
\[\bar{\chi}_{\ell}^{k-1}\varepsilon \equiv \bar{\chi}_{\ell}^{m_{1}+m_{2}}\varepsilon_{1}\varepsilon_{2} \pmod{\mathfrak{L}}.\]
Decompose $\varepsilon$ as $\varepsilon_{\ell} \varepsilon'$, where $\varepsilon_{\ell}$ is the $\ell$-part of $\varepsilon$, and $\varepsilon'$ is unramified at $\ell$. Looking at the $\ell$-part of the congruence in one hand, and at the prime-to-$\ell$ part in another hand, the congruence is equivalent to
\begin{equation}
\bar{\chi}_{\ell}^{k-1}\varepsilon_{\ell} \equiv \bar{\chi}_{\ell}^{m_{1}+m_{2}} \pmod{\mathfrak{L}} \quad \text{and} \quad \varepsilon' \equiv \varepsilon_{1}\varepsilon_{2} \pmod{\mathfrak{L}}.
\end{equation}
Applying $T_{\mathfrak{L}}$ to the second equation, we get $T_{\mathfrak{L}}\circ\pi_{\mathfrak{L}}(\varepsilon') = T_{\mathfrak{L}}\circ\pi_{\mathfrak{L}}(\varepsilon_{1}\varepsilon_{2})$. We have seen that this depends only on $\ell$. Let us look at the first equation. The projection of $\varepsilon_{\ell}$ modulo $\mathfrak{L}$ depends only on $\mathfrak{L} \cap \Q(\varepsilon)$. Moreover, $\pi_{\mathfrak{L}}(\varepsilon_{\ell})$ is a character modulo $\mathfrak{L}$ of conductor $\ell$. Therefore, there exists an integer $k_{\ell} \in \ent{0}{\ell-1}$, depending only on $\mathfrak{L}\cap\Q(\varepsilon)$, such that $\pi_{\mathfrak{L}}(\varepsilon_{\ell}) = \bar{\chi}_{\ell}^{k_{\ell}}$. The equation $\bar{\chi}_{\ell}^{k-1}\varepsilon_{\ell} \equiv \bar{\chi}_{\ell}^{m_{1}+m_{2}} \pmod{\mathfrak{L}}$ is therefore equivalent to $k+k_{\ell}-1 \equiv m_{1}+m_{2} \pmod{\ell-1}$ and depends only on $\mathfrak{L}\cap\Q(\varepsilon)$.
\end{pf}
Notice that we have in fact proved that $R_{N,k,\varepsilon}(\mathfrak{L})$ depends only on $\mathfrak{L} \cap \Q(T_{\mathfrak{L}}\circ\pi_{\mathfrak{L}}(\varepsilon_{\ell}))$ but we will only use what we have stated. A practical application of this result is that we can compute the set $R_{N,k,\varepsilon}(\mathfrak{L})$ while knowing only a prime ideal $\lambda$ below $\mathfrak{L}$ in a finite extension of $\Q(\varepsilon)$, like $K_{f}$ for example. For $\lambda$ a prime ideal in an extension of $\Q(\varepsilon)$, we will freely write $R_{N,k,\varepsilon}(\lambda)$ for the set $R_{N,k,\varepsilon}(\mathfrak{L})$ for any place $\mathfrak{L}$ above $\lambda$. We also deduce from \cref{RNkeps_dep}, a procedure to compute $R_{N,k,\varepsilon}(\lambda)$:

\begin{algo}\label{algo_RNkeps}
\textbf{Input}: Two integers $N \geq 1$, $k \geq 2$, a Dirichlet character $\varepsilon : \left(\Z/N\Z\right)^{\times} \to \C^{\times}$, and a prime ideal $\lambda$ in an finite extension of $\Q(\varepsilon)$ above a prime number $\ell$.

\noindent \textbf{Output}: The set $R_{N,k,\varepsilon}(\lambda)$.
\begin{enumerate}
\item Compute $\varepsilon_{\ell}$ and $\varepsilon'$, the $\ell$-part and prime-to-$\ell$ part of $\varepsilon$ respectively.
\item Compute the unique Dirichlet character $\varepsilon''$ modulo $N$ such that $\varepsilon''$ has prime-to-$\ell$ order, is unramified at $\ell$ and $\varepsilon'\varepsilon''^{-1}$ has order a power of $\ell$. This corresponds to the character $T_{\mathfrak{L}}\circ\pi_{\mathfrak{L}}(\varepsilon')$ for any place $\mathfrak{L}$ above $\lambda$.
\item Compute the integer $k_{\ell} \in \ent{0}{\ell-2}$ such that for all $n \in \ent{1}{N}$ prime to $N$, $\varepsilon_{\ell}(n) \equiv n^{k_{\ell}} \pmod{\lambda}$. We then have $\varepsilon_{\ell} \equiv \bar{\chi}_{\ell}^{k_{\ell}} \pmod{\lambda}$.
\item Compute the set $M_{N,k,\varepsilon}(\lambda)$ of pairs of integers $(m_{1},m_{2})$ such that $0 \leq m_{1} \leq m_{2} < \ell-1$ and $m_{1}+m_{2} \equiv k+k_{\ell}-1 \pmod{\ell-1}$.
\item Compute the set $E_{N,k,\varepsilon}(\lambda)$ of pairs of Dirichlet characters $(\varepsilon_{1},\varepsilon_{2})$ of conductor $(\mathfrak{c}_{1},\mathfrak{c}_{2})$ and such that $\varepsilon_{1}$ and $\varepsilon_{2}$ have prime-to-$\ell$ order, are unramified at $\ell$, satisfy $\varepsilon_{1}\varepsilon_{2} = \varepsilon''$ and for all primes $p \neq \ell$, we have $v_{p}\left(\frac{N}{\mathfrak{c}_{1}\mathfrak{c}_{2}}\right) \in \{0,1,2\}$.
\item Return the set $E_{N,k,\varepsilon}(\lambda) \times M_{N,k,\varepsilon}(\lambda) = R_{N,k,\varepsilon}(\lambda)$.
\end{enumerate}
\end{algo}

We now give the two main algorithm that follows from \cref{bigred_thm} et \cref{red_thm} respectively. The first algorithm computes the prime ideals $\lambda$ of $\mathcal{O}_{f}$, of residual characteristic $\ell$ such that $\ell > k+1$ and $\ell \nmid N\varphi(N)$, for which $\bar{\rho}_{f,\lambda}$ is reducible, together with the description of $\bar{\rho}_{f,\lambda}$. The correctness of the algorithm is granted by \cref{bigred_thm}.

\begin{algo}\label{algo_big}
\textbf{Input}: A newform $f$, described by its Fourier coefficients $(a_{n}(f))_{n \geq 0}$ as elements of the number field $K_{f}$, together with its level $N$, weight $k$, and character $\varepsilon$.

\noindent \textbf{Output}: The set of prime ideals $\lambda$ of $\mathcal{O}_{f}$ of residual characteristic $\ell$ such that $\ell > k+1$ and $\ell \nmid N\varphi(N)$, for which $\bar{\rho}_{f,\lambda}$ is reducible, together with the shape of $\bar{\rho}_{f,\lambda}$.

\begin{enumerate}
\item Set $\mathrm{Red}(f) = \emptyset$.
\item Compute the set $R_{N,\varepsilon}$ (see \cref{RNeps}).
\item For $(\varepsilon_{1},\varepsilon_{2}) \in R_{N,\varepsilon}$, 
\begin{enumerate}
    \item Compute $r$, $C$, and $B$ defined in \eqref{rE}, and \cref{bigred_thm} respectively.
    \item Compute the set $P(\varepsilon_{1},\varepsilon_{2})$ of prime divisors of the $\gcd$ of the algebraic norms of
    \begin{itemize}
    \item $C$;
    \item $a_{p}(f)-\varepsilon_{1}(p)-p^{k-1}\varepsilon_{2}(p)$, for $p \nmid Nr$, $p \leq B$;
    \item and $a_{p}(f)\left(a_{p}(f)-\varepsilon_{1}(p)\right)\left(a_{p}(f)-p^{k-1}\varepsilon_{2}(p)\right)$, for $p \mid N$, $p \nmid r$, $p \leq B$,
    \end{itemize}
    that are bigger than $k+1$ and do not divide $N\varphi(N)$. By \cref{bigred_thm}, these are the only prime numbers bigger than $k+1$ and not dividing $N\varphi(N)$ for which $\bar{\rho}_{f,\lambda}$ can be reducible.
\end{enumerate}
\item For $(\varepsilon_{1},\varepsilon_{2}) \in R_{N,\varepsilon}$ and for $\ell \in P(\varepsilon_{1},\varepsilon_{2})$,
\begin{enumerate}
    \item Compute the prime ideals $\lambda$ of $\mathcal{O}_{f}$ above $\ell$.
    \item For each such $\lambda$, compute a prime ideal $\mathfrak{L}$ in the ring of integers of $K_{f}(\varepsilon_{1},\varepsilon_{2})$ above $\lambda$.
    \item For each such $\mathfrak{L}$, check the following congruences.
    \begin{itemize}
    \item $C \equiv 0 \pmod{\mathfrak{L}}$;
    \item $a_{p}(f) \equiv \varepsilon_{1}(p)+p^{k-1}\varepsilon_{2}(p) \pmod{\mathfrak{L}}$ for all $p \nmid Nr$, $p \leq B$;
    \item $a_{p}(f)(a_{p}(f)-\varepsilon_{1}(p))(a_{p}(f)-p^{k-1}\varepsilon_{2}(p)) \equiv 0 \pmod{\mathfrak{L}}$ for all $p \mid N$, $p \nmid r$, $p \leq B$.
    \end{itemize}
    If they all hold, add $(\lambda,\varepsilon_{1},\varepsilon_{2},0,k-1)$ to $\mathrm{Red}(f)$. By \cref{bigred_thm}, $\bar{\rho}_{f,\lambda}$ is reducible and we have $\bar{\rho}_{f,\lambda} \cong \bar{\varepsilon_{1}} \oplus \bar{\chi}_{\ell}^{k-1} \bar{\varepsilon_{2}}$.
\end{enumerate}
\item Return $\mathrm{Red}(f)$.
\end{enumerate}
\end{algo}
We now turn to the computation of the reducible primes of residue characteristic $\ell$ such that $\ell \leq k+1$ or $\ell \mid N\varphi(N)$. The correctness of the following algorithm follows by \cref{red_thm}.

\begin{algo}\label{algo_small}
\textbf{Input}: A newform $f$, described by its Fourier coefficients $(a_{n}(f))_{n \geq 0}$ as elements of the number field $K_{f}$, together with its level $N$, weight $k$, and character $\varepsilon$.

\noindent \textbf{Output}: The set of prime ideals $\lambda$ of $\mathcal{O}_{f}$ of residual characteristic $\ell$ such that $\ell \leq k+1$ or $\ell \mid N\varphi(N)$, for which $\bar{\rho}_{f,\lambda}$ is reducible, together with the shape of $\bar{\rho}_{f,\lambda}$.

\begin{enumerate}
\item Set $\mathrm{Red}(f) = \emptyset$.
\item Compute the set $P$ of prime numbers $\ell$ such that $\ell \leq k+1$ or $\ell \mid N\varphi(N)$.
\item For each $\ell \in P$, compute the set $P(\ell)$ of prime ideals $\lambda$ in $\mathcal{O}_{f}$ above $\ell$.
\item For each $\ell \in P$ and for each $\lambda \in P(\ell)$, compute $R_{N,k,\varepsilon}(\lambda)$ using \cref{algo_RNkeps}. We can do this because we have $\Q(\varepsilon) \subset K_{f}$.
\item For each $\ell \in P$, for each $\lambda \in P(\ell)$, and for each $(\varepsilon_{1},\varepsilon_{2},m_{1},m_{2}) \in R_{N,k,\varepsilon}(\lambda)$,
\begin{enumerate}
    \item Compute a prime ideal $\mathfrak{L}$ in the ring of integers of $K_{f}(\varepsilon_{1},\varepsilon_{2})$ above $\lambda$.
    \item Compute $r$ and $B$ defined in \eqref{rE} and \cref{red_thm} respectively.
    \item Check the following congruences.
    \begin{itemize}
    \item $a_{p}(f) \equiv p^{m_{1}}\varepsilon_{1}(p) + p^{m_{2}}\varepsilon_{2}(p) \pmod{\mathfrak{L}}$ for all $p \nmid Nr$, $p \leq B$;
    \item $a_{p}(f)(a_{p}(f)-p^{m_{1}}\varepsilon_{1}(p))(a_{p}(f)-p^{m_{2}}\varepsilon_{2}(p)) \equiv 0 \pmod{\mathfrak{L}}$ for all $p \mid N$, $p \nmid r$, $p \leq B$.
    \end{itemize}
    If they all hold, add $(\lambda,\varepsilon_{1},\varepsilon_{2},m_{1},m_{2})$ to $\mathrm{Red}(f)$. By \cref{red_thm}, $\bar{\rho}_{f,\lambda}$ is reducible and we have $\bar{\rho}_{f,\lambda} \cong \bar{\chi}_{\ell}^{m_{1}} \bar{\varepsilon_{1}} \oplus \bar{\chi}_{\ell}^{m_{2}}\bar{\varepsilon_{2}}$.
\end{enumerate}
\item Return $\mathrm{Red}(f)$.
\end{enumerate}
\end{algo}

The correctness of \cref{algo_big} and \cref{algo_small} follows directly from \cref{bigred_thm} and \cref{red_thm} respectively. The most time-consuming computation is step $3(b)$ of \cref{algo_big}. This depends on the size of the “big” reducible primes. We have implemented these algorithms in PARI/GP \cite{PARI} and we have been able to execute them as long as the degree of $K_{f}$ keeps reasonable (say $[K_{f}:\Q] \leq 20$). The second limiting factor being the weight $k$ that controls the size of the Fourier coefficients of $f$. Our code is available on GitHub at the following address:
\begin{center}
\url{https://github.com/bpeaucelle/mfreducible}
\end{center}

\subsection{Numerical Examples}
We present here some examples of applications of the algorithms described above to compute the reducible primes of a given newform. Throughout this section we use the Conrey representation $(\varepsilon_{a}(b))_{b \wedge a = 1}$ for the Dirichlet characters of modulus $a$. This is the way they are described in the LMFDB for example (and in some extent in PARI/GP). Notice that there would be no possible confusion with the characters $\varepsilon_{1}$, $\varepsilon_{2}$ from the algorithms above.

\subsubsection{A concrete example}
Consider the space $\Smod_{7}^{\mathrm{new}}(7,\varepsilon_{7}(3))$. It has dimension $6$ over $\C$ and is generated by $2$ newforms, $f_{1}$ and $f_{2}$, up to conjugation by $\Gal(\bar{\Q}/\Q(\varepsilon_{7}(3)))$. We have $K_{f_{1}} = \Q[t]/(t^{2}-t+1)$ and $K_{f_{2}} = \Q[x]/(x^{4}+2x^{2}+4)$. Notice that $\left(1,x,\frac{x^{2}}{2},\frac{x^{3}}{2}\right)$ is a integer basis of $\mathcal{O}_{f_{2}}$, and that $\varepsilon_{7}(3)$ sends $3$ to $t$ in $K_{f_{1}}$ and to $-\frac{x^{2}}{2}$ in $K_{f_{2}}$. The $q$-expansions of $f_{1}$ and $f_{2}$ are given by
\[f_{1} = q + 12tq^{2} + (-7t-7)q^{3} + (80t-80)q^{4} + (-105t+210)q^{5} + (-168t+84)q^{6} - 343q^{7} + O(q^{8}),\]
\begin{align*}
f_{2} = q + \left(\frac{3}{2}x^{3}+2x^{2}+3x\right)q^{2} &+ \left(13x^{3}-\frac{3}{2}x^{2} + 13x + 3\right)q^{3} + (15x^{2}-24x+30)q^{4}\\
    &+ \left(-25x^{3}-\frac{25}{2}x^{2}+50x-50\right)q^{5} + O(q^{6}).
\end{align*}
The set of prime numbers less than $k+1 = 8$ and dividing $N\varphi(N) = 42$ is therefore $\{2,3,5,7\}$. We treat those primes separately below.

\begin{itemize}
\item $\ell = 2$: The ideal $2\mathcal{O}_{f_{1}}$ is prime and we have $2\mathcal{O}_{f_{2}} = \left((x,2)\mathcal{O}_{f_{2}}\right)^{2}$. Because, the ideal generated by $2$ in $\Z[\varepsilon_{7}(3)]$ is prime, we have 
\[R_{7,7,\varepsilon_{7}(3)}(2\mathcal{O}_{f_{1}}) = R_{7,7,\varepsilon_{7}(3)}((2,x)\mathcal{O}_{f_{2}}) = \{(\mathds{1},\varepsilon_{7}(4),0,0);(\varepsilon_{7}(4),\mathds{1},0,0)\}.\]
In fact, the two elements of this set give rise to the same representation. We only study the case $(\varepsilon_{1},\varepsilon_{2},m_{1},m_{2}) = (\varepsilon_{7}(4),\mathds{1},0,0)$. According to \cref{red_thm}, we have
\[k' = 2, \quad r = 1, \quad N' = 7, \quad a = 0, \quad b = 3, \quad \tilde{k} = 10, \quad B = 6 + \frac{2}{3}.\]
To check the reducibility of $\bar{\rho}_{f_{1},(2)}$ and $\bar{\rho}_{f_{2},(2,x)}$ we only have to check the third and fifth coefficients of $f_{1}$ and $f_{2}$. The table below shows the reduction modulo the prime ideals above of $a_{p}(f_{i})-1-\varepsilon_{7}(4)(p)$ for $i = 1$, $2$, and $p = 3$, $5$. From \cref{red_thm}, we know that $\bar{\rho}_{f_{i},\lambda}$ is reducible if and only if the row corresponding to $f_{i}$ contains only zeros.

\begin{minipage}{.52\linewidth}
\begin{center}
\begin{tabular}{|c|c|c|}
\hline
$p$ & $3$ & $5$\\
\hline
$a_{p}(f_{1})-\left(1+\varepsilon_{7}(4)(p)\right)\pmod{2}$ & $0$ & $0$\\
\hhline{|=|=|=|}
$a_{p}(f_{2})-\left(1+\varepsilon_{7}(4)(p)\right)\pmod{(2,x)}$ & $0$ & $0$\\
\hline
\end{tabular}
\end{center}
\end{minipage}
\begin{minipage}{.43\linewidth}
Therefore, we have $\bar{\rho}_{f_{1},(2)} \cong \mathds{1} \oplus \bar{\varepsilon_{7}(4)}$ and $\bar{\rho}_{f_{2},(2,x)} \cong \mathds{1} \oplus \bar{\varepsilon_{7}(4)}$.
\end{minipage}

\item $\ell = 3$: We have $3\mathcal{O}_{f_{1}} = \left((3,t+1)\mathcal{O}_{f_{1}}\right)^{2}$ and $3\mathcal{O}_{f_{2}} = \left((3,x^2+1)\mathcal{O}_{f_{2}}\right)^{2}$. As for $\ell = 2$, the ideal generated by $3$ in $\Z[\varepsilon_{7}(3)]$ is prime. Therefore, we have
\[R_{7,7,\varepsilon_{7}(3)}((3,t+1)\mathcal{O}_{f_{1}}) = R_{7,7,\varepsilon_{7}(3)}((3,x^{2}+1)\mathcal{O}_{f_{2}}) = \{(\mathds{1},\varepsilon_{7}(6));(\varepsilon_{7}(6),\mathds{1})\} \times \{(0,0);(1,1)\}.\]
As above, we treat only the cases $(\varepsilon_{1},\varepsilon_{2},m_{1},m_{2}) = (\varepsilon_{7}(6),\mathds{1},0,0)$ and $(\varepsilon_{7}(6),\mathds{1},1,1)$. According to \cref{red_thm}, we have in both cases
\[k' = 1, \quad r = 1, \quad N' = 7, \quad 
a = 0, \quad b = 4, \quad \tilde{k} = 11, \quad B = 7 + \frac{1}{3}.\]
We have to look at the second, fifth and seventh coefficients of $f_{1}$ and $f_{2}$. Let look at the second and fifth first.

\begin{center}
\begin{tabular}{|c|c|c|}
\hline
$p$ & $2$ & $5$\\
\hline
$a_{p}(f_{1}) - \left(1+\varepsilon_{7}(6)(p)\right) \pmod{(3,t+1)}$ & $1$ & $0$\\
\hline
$a_{p}(f_{1}) - \left(p+p\varepsilon_{7}(6)(p)\right) \pmod{(3,t+1)}$ & $2$ & $0$\\
\hhline{|=|=|=|}
$a_{p}(f_{2}) -\left(1+\varepsilon_{7}(6)(p)\right) \pmod{(3,x^{2}+1)}$ & $2$ & $0$\\
\hline
$a_{p}(f_{2}) - \left(p+p\varepsilon_{7}(6)(p)\right) \pmod{(3,x^{2}+1)}$ & $0$ & $0$\\
\hline
\end{tabular}
\end{center}

From these computations, we deduce that the only representation that can be reducible is $\bar{\rho}_{f_{2},(3,x^{2}+1)}$, and that it can only isomorphic to $\bar{\chi}_{3} \oplus \bar{\chi}_{3}\bar{\varepsilon_{7}(6)}$. To confirm this isomorphism, we finally have to check that there exists some $b_{7} \in \{0,7,7\varepsilon_{7}(6)(7)\} = \{0,7\}$ such that $a_{7}(f_{2}) \equiv 7b_{7} \pmod{(3,x^{2}+1)}$. We find that we have $a_{7}(f_{2}) \equiv 7 \pmod{(3,x^{2}+1)}$. Therefore, the representation $\bar{\rho}_{f_{1},(3,t+1)}$ is irreducible and we have $\bar{\rho}_{f_{2},(3,x^{2}+1)} \cong \bar{\chi}_{3} \oplus \bar{\chi}_{3}\bar{\varepsilon_{7}(6)}$.

\item $\ell = 5$: We have that $5$ is prime in $\mathcal{O}_{f_{1}}$ and $5\mathcal{O}_{f_{2}} = (5,x^{2}-2x-2)(5,x^{2}+2x-2)$. There is again only one prime ideal above $5$ in $\Z[\varepsilon_{7}(3)]$ and we have
\[R_{7,7,\varepsilon_{7}(3)}(5\mathcal{O}_{f_{1}}) = R_{7,7,\varepsilon_{7}(3)}(5,x^{2}\pm2x-2) = \{(\mathds{1},\varepsilon_{7}(3));(\varepsilon_{7}(3),\mathds{1})\} \times \{(0,2);(1,1)\}.\]
Looking at the congruences at $p = 3$ for $f_{1}$ and $p = 2$ for $f_{2}$, we have
\[a_{3}(f_{1})-\left(1+3^{2}\varepsilon_{7}(3)(3)\right) \equiv 4t+2 \pmod{5}, \quad a_{3}(f_{1})-\left(\varepsilon_{7}(3)(3) + 3^{2}\right) \equiv 2t+4 \pmod{5},\]
\[a_{2}(f_{2})-\left(1+2^{2}\varepsilon_{7}(3)(2)\right) \equiv \left\{\begin{array}{l}
2 \pmod{(5,x^{2}-2x-2});\\
4x \pmod{(5,x^{2}+2x-2)},
\end{array}\right.\]
\[a_{2}(f_{2}) - \left(\varepsilon_{7}(3)(2)+2^{2}\right) \equiv \left\{\begin{array}{l}
2x+3 \pmod{(5,x^{2}-2x-2)};\\
2x+1 \pmod{(5,x^{2}+2x-2)},
\end{array}\right.\]
\[a_{2}(f_{2}) - \left(2+2\varepsilon_{7}(3)(2)\right) \equiv \left\{\begin{array}{l}
3x+2 \pmod{(5,x^{2}-2x-2)};\\
x \pmod{(5,x^{2}+2x-2)}.
\end{array}\right.\]
The only candidate remaining is $(\varepsilon_{7}(3),\mathds{1},1,1)$ for $\bar{\rho}_{f_{1},(5)}$ (and $(\mathds{1},\varepsilon_{7}(3),1,1)$ which gives the same representation). We have
\[k' = 1, \quad r = 1, \quad N' = 7, \quad a = 0, \quad b = 6, \quad \tilde{k} = 13, \quad B = 8 + \frac{2}{3}.\]
We check the second, third and seventh coefficients and we get
\[a_{2}(f_{1}) \equiv 2 + 2\varepsilon_{7}(3)(2) \pmod{5}, \quad a_{3}(f_{1}) \equiv 3 + 3\varepsilon_{7}(3)(3) \pmod{5}, \quad a_{7}(f_{1}) \equiv 7 \pmod{5}.\] Therefore, the representations $\bar{\rho}_{f_{2},(5,x^{2}-2x-2)}$ and $\bar{\rho}_{f_{2},(5,x^{2}+2x-2)}$ are irreducible, and we have an isomorphism $\bar{\rho}_{f_{1},(5)} \cong \bar{\chi}_{5} \oplus \bar{\chi}_{5}\bar{\varepsilon_{7}(3)}$.

\item $\ell = 7$: We have $7\mathcal{O}_{f_{1}} = (7,t-5)(7,t-3)$ and $7\mathcal{O}_{f_{2}} = (7,x-1)(7,x-2)(7,x+2)(7,x+1)$. This time $7$ decomposes in $\Z[\varepsilon_{7}(3)]$ and we have
\begin{align*}
&R_{7,7,\varepsilon_{7}(3)}(7,t-3) = R_{7,7,\varepsilon_{7}(3)}(7,x \pm 1) = \{(\mathds{1},\mathds{1})\} \times \{(0,1) ; (2,5) ; (3,4)\}\\
\text{and } &R_{7,7,\varepsilon_{7}(3)}(7,t-5) = R_{7,7,\varepsilon_{7}(3)}(7,x \pm 2) = \{(\mathds{1},\mathds{1})\} \times \{(0,5) ; (1,4) ; (2,3)\}.
\end{align*}
For $f_{1}$, looking at $p = 2$ leaves us only with $(\varepsilon_{1},\varepsilon_{2},m_{1},m_{2}) = (\mathds{1},\mathds{1},2,5)$ for the ideal $(7,t-3)$ and $(\mathds{1},\mathds{1},1,4)$ for $(7,t-5)$. In both cases we have to look at congruences up to $p = 5$, and we get
\[\bar{\rho}_{f_{1},(7,t-3)} \cong \bar{\chi}_{7}^{2} \oplus \bar{\chi}_{7}^{5} \quad \text{and} \quad \bar{\rho}_{f_{1},(7,t-5)} \cong \bar{\chi}_{7} \oplus \bar{\chi}_{7}^{4}.\]
For $f_{2}$, looking at $p = 3$ leaves us with $(\mathds{1},\mathds{1},2,5)$ for $(7,x+1)$, $(\mathds{1},\mathds{1},1,4)$ for $(7,x+2)$, $(\mathds{1},\mathds{1},2,3)$ for $(7,x-2)$, and $(\mathds{1},\mathds{1},3,4)$ for $(7,x-1)$. In the two first cases we have $r = 1$ and we have to look at congruences up to $p = 5$. In the two last cases we have $r = 4$ and we have to check congruences up to $p = 53$ and $p = 67$ respectively. In every case, \cref{red_thm} shows that the corresponding representation is reducible. To sum up we have
\[\bar{\rho}_{f_{2},(7,x-1)} \cong \bar{\chi}_{7}^{3} \oplus \bar{\chi}_{7}^{4}, \quad \bar{\rho}_{f_{2},(7,x+1)} \cong \bar{\chi}_{7}^{2} \oplus \bar{\chi}_{7}^{5},\]
\[\bar{\rho}_{f_{2},(7,x-2)} \cong \bar{\chi}_{7}^{2} \oplus \bar{\chi}_{7}^{3}, \quad \bar{\rho}_{f_{2},(7,x+2)} \cong \bar{\chi}_{7} \oplus \bar{\chi}_{7}^{4}.\]
\end{itemize}
We finally look at the prime numbers $\ell > 7$. We have $R_{7,\varepsilon_{7}(3)} = \{(\mathds{1},\varepsilon_{7}(3)),(\varepsilon_{7}(3),\mathds{1})\}$. Let $(\varepsilon_{1},\varepsilon_{2}) \in R_{7,\varepsilon_{7}(3)}$. We have $r = 1$, $N' = 1$, $B = 4 + \frac{1}{3}$, and
\[C(\varepsilon_{1},\varepsilon_{2}) = \left\{\begin{array}{ll}
0 & \text{if } (\varepsilon_{1},\varepsilon_{2}) =  (\varepsilon_{7}(3),\mathds{1});\\
-\frac{B_{7,\varepsilon_{7}(3)}}{14} a_{7}\left(f_{i}\right)^{2} & \text{if } (\varepsilon_{1},\varepsilon_{2}) =  (\mathds{1},\varepsilon_{7}(3)).
\end{array}\right.\]
We first look at $f_{1}$. We find that $43$ is the only prime factor greater than $7$ of the $\gcd$ of the algebraic norms of $C(\varepsilon_{1},\varepsilon_{2})$ and $a_{p}(f_{1})-\varepsilon_{1}(p) - p^{6}\varepsilon_{2}(p)$, for $p = 2$, $3$. In $\mathcal{O}_{f_{1}}$ we have $43\mathcal{O}_{f_{1}} = (43,t-7)(43,t+6)$ and we get the following table.
\begin{center}
\begin{tabular}{|c|c|c|c||c|c|}
\hline
$(\varepsilon_{1},\varepsilon_{2})$ & \multicolumn{3}{c||}{$(\mathds{1},\varepsilon_{7}(3))$} & \multicolumn{2}{c|}{$(\varepsilon_{7}(3),\mathds{1})$}\\
\hline
 & \multirow{2}{*}{$C(\mathds{1},\varepsilon_{7}(3))$} & \multicolumn{2}{c||}{$a_{p}(f_{1})-1-p^{6}\varepsilon_{7}(3)(p)$} & \multicolumn{2}{c|}{$a_{p}(f_{1})-\varepsilon_{7}(3)(p)-p^{6}$}\\
\cline{3-6}
 & & $p = 2$ & $p = 3$ & $p = 2$ & $p = 3$\\
\hline
$(43,t-7)$ & $0$ & $0$ & $0$ & $14$ & $25$\\
\hline
$(43,t+6)$ & $40$ & $31$ & $22$ & $0$ & $0$\\
\hline
\end{tabular}
\end{center}
Therefore, we get $\bar{\rho}_{f_{1},(43,t-7)} \cong \mathds{1} \oplus \bar{\chi}_{43}^{6}\bar{\varepsilon_{7}(3)}$ and $\bar{\rho}_{f_{1},(43,t+6)} \cong \bar{\varepsilon_{7}(3)} \oplus \bar{\chi}_{43}^{6}$.

We now turn to $f_{2}$. Computing again the $\gcd$ of the algebraic norm of $C(\varepsilon_{1},\varepsilon_{2})$ and $a_{p}(f_{2})-\varepsilon_{1}(p)-p^{7}\varepsilon_{2}(p)$, for $p = 2$, $3$, we find that the only possible residue characteristics are $97$ and $3919$. We have the following decompositions in $\mathcal{O}_{f_{2}}$:
\[97\mathcal{O}_{f_{2}} = (97,x-19)(97,x-5)(97,x+5)(97,x+19),\]
\[3919\mathcal{O}_{f_{2}} = (3919,x-934)(3919,x-621)(3919,x+621)(3919,x+934).\]

\begin{center}
\begin{tabular}{|c|c|c|c||c|c|}
\hline
$(\varepsilon_{1},\varepsilon_{2})$ & \multicolumn{3}{c||}{$(\mathds{1},\varepsilon_{7}(3))$} & \multicolumn{2}{c|}{$(\varepsilon_{7}(3),\mathds{1})$}\\
\hline
 & \multirow{2}{*}{$C(\mathds{1},\varepsilon_{7}(3))$} & \multicolumn{2}{c||}{$a_{p}(f_{2})-1-p^{6}\varepsilon_{7}(3)(p)$} & \multicolumn{2}{c|}{$a_{p}(f_{2})-\varepsilon_{7}(3)(p)-p^{6}$}\\
\cline{3-6}
 & & $p = 2$ & $p = 3$ & $p = 2$ & $p = 3$\\
\hhline{|=|=|=|=||=|=|}
$(97,x-19)$ & $9$ & $33$ & $75$ & $30$ & $57$\\
\hline
$(97,x-5)$ & $0$ & $0$ & $0$ & $8$ & $66$ \\
\hline
$(97,x+5)$ & $0$ & $80$ & $15$ & $88$ & $81$ \\
\hline
$(97,x+19)$ & $11$ & $3$ & $18$ & $0$ & $0$ \\
\hhline{|=|=|=|=||=|=|}
$(3919,x-934)$ & $3160$ & $3231$ & $1337$ & $0$ & $0$ \\
\hline
$(3919,x-621)$ & $0$ & $0$ & $0$ & $3042$ & $609$ \\
\hline
$(3919,x+621)$ & $0$ & $1685$ & $2010$ & $808$ & $2619$ \\
\hline
$(3919,x+934)$ & $1455$ & $3038$ & $3047$ & $3726$ & $1710$ \\
 \hline
\end{tabular}
\end{center}
Therefore, the representations $\bar{\rho}_{f_{2},(97,x-19)}$, $\bar{\rho}_{f_{2},(97,x+5)}$, $\bar{\rho}_{f_{2},(3919,x + 621)}$ and $\bar{\rho}_{f_{2},(3919,x+934)}$ are irreducible and we have
\[\bar{\rho}_{f_{2},(97,x-5)} \cong \mathds{1} \oplus \bar{\chi}_{97}^{6} \bar{\varepsilon_{7}(3)}, \quad \bar{\rho}_{f_{2},(97,x+19)} \cong \bar{\varepsilon_{7}(3)} \oplus \bar{\chi}_{97}^{6},\]
\[\bar{\rho}_{f_{2},(3919,x-934)} \cong \bar{\varepsilon_{7}(3)} \oplus \bar{\chi}_{3919}^{6}, \quad \bar{\rho}_{f_{2},(3919,x-621)} \cong \mathds{1} \oplus \bar{\chi}_{3919}^{6}\bar{\varepsilon_{7}(3)}.\]
The following table sums up all the cases for which $\bar{\rho}_{f_{i},\lambda}$ is reducible.

\begin{center}
\begin{tabular}{|c|c|c||c|c|c|c|}
\hline
$\ell$ & \multicolumn{2}{c||}{$f_{1}$} & \multicolumn{4}{c|}{$f_{2}$}\\
\hline
\multirow{2}{*}{$2$} & \multicolumn{2}{c||}{$(2)$} & \multicolumn{4}{c|}{$(2,x)$}\\
 & \multicolumn{2}{c||}{$\mathds{1} \oplus \bar{\varepsilon_{7}(4)}$} & \multicolumn{4}{c|}{$\mathds{1} \oplus \bar{\varepsilon_{7}(4)}$}\\
\hline
\multirow{2}{*}{$3$} & \multicolumn{2}{c||}{\multirow{2}{*}{Irreducible}} & \multicolumn{4}{c|}{$(3,x^{2}+1)$}\\
 & \multicolumn{2}{c||}{} & \multicolumn{4}{c|}{$\bar{\chi}_{3} \oplus \bar{\chi}_{3}\bar{\varepsilon_{7}(6)}$}\\
\hline
\multirow{2}{*}{$5$} & \multicolumn{2}{c||}{(5)} & \multicolumn{4}{c|}{\multirow{2}{*}{Irreducible}}\\
 & \multicolumn{2}{c||}{$\bar{\chi}_{5} \oplus \bar{\chi}_{5} \bar{\varepsilon_{7}(3)}$} & \multicolumn{4}{c|}{}\\
\hline
\multirow{2}{*}{$7$} & $(7,t-3)$ & $(7,t-5)$ & $(7,x-2)$ & $(7,x-1)$ & $(7,x+1)$ & $(7,x+2)$\\
 & $\bar{\chi}_{7}^{2} \oplus \bar{\chi}_{7}^{5}$ & $\bar{\chi}_{7} \oplus \bar{\chi}_{7}^{4}$ & $\bar{\chi}_{7}^{2} \oplus \bar{\chi}_{7}^{3}$ & $\bar{\chi}_{7}^{3} \oplus \bar{\chi}_{7}^{4}$ & $\bar{\chi}_{7}^{2} \oplus \bar{\chi}_{7}^{5}$ & $\bar{\chi}_{7} \oplus \bar{\chi}_{7}^{4}$\\
 \hline
\multirow{2}{*}{$\ell > k+1$} & & & \multicolumn{2}{c|}{$(97,x-5)$} & \multicolumn{2}{c|}{$(97,x+19)$}\\
 & $(43,t-7)$ & $(43,t+6)$ & \multicolumn{2}{c|}{$\mathds{1} \oplus \bar{\chi}_{97}^{6}\bar{\varepsilon_{7}(3)}$} & \multicolumn{2}{c|}{$\bar{\varepsilon_{7}(3)} \oplus \bar{\chi}_{97}^{6}$}\\
\cline{4-7}
\multirow{2}{*}{and $\ell \nmid N\varphi(N)$} & $\mathds{1} \oplus \bar{\chi}_{47}^{6} \bar{\varepsilon_{7}(3)}$ & $\bar{\varepsilon_{7}(3)} \oplus \bar{\chi}_{43}^{6}$ & \multicolumn{2}{c|}{$(3919,x-934)$} & \multicolumn{2}{c|}{$(3919,x-621)$}\\
 & & & \multicolumn{2}{c|}{$\bar{\varepsilon_{7}(3)} \oplus \bar{\chi}_{3919}^{6}$} & \multicolumn{2}{c|}{$\mathds{1} \oplus \bar{\chi}_{3919}^{6}\bar{\varepsilon_{7}(3)}$}\\
\hline
\end{tabular}
\end{center}

\subsubsection{Irreducible everywhere representation}

We present an example of a form which all residual representations are irreducible. Fix $(N,k,\varepsilon) = (35,4,\mathds{1})$. The space $\Smod_{4}^{\mathrm{new}}(35,\mathds{1})$ has dimension $6$ over $\C$ and contains $3$ newforms up to conjugation by $\Gal(\bar{\Q}/\Q)$. Let $f$ be the newform of this space which $q$-expansion is
\[f = q + (y+4)q^{2} + (1-4y)q^{3} + O\left(q^{4}\right),\]
where $y$ is a root of $X^{2}-2$, and generates $K_{f}$. We have in this case
\[R_{35,\mathds{1}} = \{(\mathds{1},\mathds{1})\}.\]
Therefore, by \cref{red_bound} the only prime ideals $\lambda$ of $\mathcal{O}_{f}$ for which $\bar{\rho}_{f,\lambda}$ can be reducible are of residue characteristic $\ell \in \{2,3,5,7\}$ (because we have $B_{4,\mathds{1}} = -\frac{1}{30}$). Let us look at each of these cases.
\begin{itemize}
\item $\ell = 2$: We have $2\mathcal{O}_{f} = (2,y)^{2}$ and $R_{35,4,\mathds{1}}(2,y) = \{(\mathds{1},\mathds{1},0,0)\}$. However, we have
\[\mathrm{Tr}\left(\bar{\rho}_{f,(2,y)}(\mathrm{Frob}_{3})\right) \equiv a_{3}(f) \equiv 1 \pmod{(2,y)}\]
\[\text{and}\]
\[\mathrm{Tr}\left((\mathds{1}\oplus\mathds{1})(\mathrm{Frob}_{3})\right) \equiv 0 \pmod{(2,y)}.\]
Therefore $\bar{\rho}_{f,(2,y)}$ is irreducible.

\item $\ell = 3$: The ideal $3\mathcal{O}_{f}$ is prime and we have $R_{35,4,\mathds{1}}(3) = \{(\mathds{1},\mathds{1},0,1)\}$. However, we have
\[\mathrm{Tr}\left(\bar{\rho}_{f,(3)}\left(\mathrm{Frob}_{2}\right)\right) \equiv a_{2}(f) \equiv y+1 \pmod{3}\]
\[\text{and}\]
\[\mathrm{Tr}\left(\left(\mathds{1}\oplus\bar{\chi}_{3}\right)\left(\mathrm{Frob}_{2}\right)\right) \equiv 0 \pmod{3}.\]
Therefore, $\bar{\rho}_{f,(3)}$ is irreducible.

\item $\ell = 5$: Again, $5$ is prime in $\mathcal{O}_{f}$ and we have $R_{35,4,\mathds{1}}(5) = \{(\mathds{1},\mathds{1},0,3);(\mathds{1},\mathds{1},1,2)\}$. Looking at a Frobenius element at $2$, we have
\[\mathrm{Tr}\left(\bar{\rho}_{f,(5)}\left(\mathrm{Frob}_{2}\right)\right) \equiv a_{2}(f) \equiv y+4 \pmod{5}\]
\[\text{and}\]
\[\mathrm{Tr}\left(\left(\mathds{1}\oplus\bar{\chi}_{5}^{3}\right)\left(\mathrm{Frob}_{2}\right)\right) \equiv 4 \pmod{5}, \quad \mathrm{Tr}\left(\left(\bar{\chi}_{5}\oplus\bar{\chi}_{5}^{2}\right)\left(\mathrm{Frob}_{2}\right)\right) \equiv 0 \pmod{5}.\]
Therefore, $\bar{\rho}_{f,(5)}$ is irreducible.

\item $\ell = 7$: We have $7\mathcal{O}_{f} = (7,y-3)(7,y+3)$ and $R_{35,7,\mathds{1}}(7,y\pm 3) = \{(\mathds{1},\mathds{1})\} \times \{(0,3);(1,2);(4,5)\}$. However, for a Frobenius element at $7$ we have
\[\left\{\begin{array}{l} \mathrm{Tr}\left(\bar{\rho}_{f,(7,y-3)}\left(\mathrm{Frob}_{3}\right)\right) \equiv 
3 \pmod{(7,y-3)};\\
\mathrm{Tr}\left(\bar{\rho}_{f,(7,y+3)}\left(\mathrm{Frob}_{3}\right)\right) \equiv 6 \pmod{(7,y+3)},
\end{array}\right.\]
\[\text{and}\]
\[\mathrm{Tr}\left(\left(\mathds{1}\oplus\bar{\chi}_{7}^{3}\right)\left(\mathrm{Frob}_{3}\right)\right) \equiv 0 \pmod{7}, \quad \mathrm{Tr}\left(\left(\bar{\chi}_{7}\oplus\bar{\chi}_{7}^{2}\right)\left(\mathrm{Frob}_{3}\right)\right) \equiv 5 \pmod{7},\]
\[\mathrm{Tr}\left(\left(\bar{\chi}_{7}^{4}\oplus\bar{\chi}_{7}^{5}\right)\left(\mathrm{Frob}_{3}\right)\right) \equiv 2 \pmod{7}.\]
Therefore, the representations $\bar{\rho}_{f,(7,y-3)}$ and $\bar{\rho}_{f,(7,y+3)}$ are irreducible.
\end{itemize}
Thus, for all prime ideals $\lambda$ in $\mathcal{O}_{f}$, the representation $\bar{\rho}_{f,\lambda}$ is irreducible.

\printbibliography
\end{document}